\documentclass[11pt]{amsart}
\usepackage{amsmath,amsfonts,amssymb}
\def\wh{\widehat}
\def\wt{\widetilde}
\def\R{\mathbb R}
\def\C{\mathbb C}
\def\Z{\mathbb Z}

\def\N{\mathbb N}
\def\A{\mathcal A}
\def\B{\mathcal B}
\def\u{\mathbf u}
\def\e{\mathbf e}
\def\f{\mathbf f}
\def\x{\mathbf x}
\def\v{\mathbf v}
\def\a{\mathbf a}
\def\s{\mathbf s}
\def\p{\mathbf p}
\def\r{\mathbf r}
\def\b{\mathbf b}
\def\c{\mathbf c}

\def\w{\mathbf w}
\def\z{\mathbf z}
\def\h{\mathbf h}

\def\D{\mathbf D}
\def\FF{\mathbf F}
\def\1{\mathbf 1}

\def\H{\mathfrak H}

\def\O{\mathcal O}

\def\bxi{\boldsymbol \xi}
\def\bnu{\boldsymbol \nu}

\def\eeta{\boldsymbol \eta}
\def\bphi{\boldsymbol \phi}
\def\bpsi{\boldsymbol \psi}
\def\brho{\boldsymbol \rho}

\def\bPhi{\boldsymbol \Phi}
\def\bvarphi{\boldsymbol \varphi}
\def\1{\bold 1}

\def\div{\mathrm{div}\,}

\def\rank{\mathrm{rank}\,}

\def\eps{\varepsilon}

\def\Ker{\mathrm{Ker}\,}
\def\meas{\mathrm{meas}\,}

\def\leq{\leqslant}
\def\le{\leqslant}
\def\geq{\geqslant}
\def\ge{\geqslant}

\begin{document}

\title[Homogenization of the Neumann problem for elliptic systems]{Homogenization of the Neumann problem
 for elliptic systems with periodic coefficients}

\author{T.~A.~Suslina}

\address{St.~Petersburg State University, Department of Physics,
Ul'yanovskaya 3, Petrodvorets, St.~Petersburg, 198504, RUSSIA}

\email{suslina@list.ru}


\thanks{Work done during a visit to the Mittag-Leffler Institute (Djursholm, Sweden).}



\begin{abstract}
Let $\O \subset \R^d$ be a bounded domain with the boundary of class $C^{1,1}$.
In $L_2(\O;\C^n)$, a matrix elliptic second order differential operator
$\A_{N,\eps}$ with the Neumann boundary condition is considered. Here \hbox{$\eps>0$} is a small parameter,
the coefficients of $\A_{N,\eps}$ are periodic
and depend on $\x/\eps$. There are no regularity assumptions on the coefficients.
It is shown that the resolvent \hbox{$(\A_{N,\eps}+\lambda I)^{-1}$}
converges in the $L_2(\O;\C^n)$-operator norm to the resolvent of the effective operator $\A_N^0$ with constant
coefficients, as $\eps \to 0$.
A sharp order error estimate
$\|(\A_{N,\eps}+\lambda I)^{-1} - (\A_{N}^0 +\lambda I)^{-1}\|_{L_2\to L_2} \le C\eps$
is obtained.
Approximation for the resolvent
$(\A_{N,\eps}+\lambda I)^{-1}$ in the norm of operators acting from $L_2(\O;\C^n)$ to the Sobolev space
$H^1(\O;\C^n)$ with an error $O(\sqrt{\eps})$ is found.
Approximation is given by the sum of the operator $(\A^0_N +\lambda I)^{-1}$
and the first order corrector.
In a strictly interior subdomain $\O'$ a similar approximation with an error $O(\eps)$ is obtained.
\end{abstract}

\maketitle
\section*{Introduction}

The paper concerns homogenization theory for periodic differential operators (DO's).
A broad literature is devoted to homogenization problems in the small period limit.
First of all, we mention the books [BeLPa], [BaPan], [ZhKO].

\noindent\textbf{0.1. Operator-theoretic approach to homogenization problems.}
In a series of papers [BSu1-4] by M.~Sh.~Birman and T.~A.~Suslina an
operator-theoretic (spectral) approach to homogenization problems was suggested and developed.
By this approach, the so-called operator error estimates in homogenization problems for elliptic DO's
were obtained.
Matrix strongly elliptic DO's acting in $L_2(\R^d;\C^n)$
and admitting a factorization of the form
$$
\A_\eps = b(\D)^* g(\x/\eps)b(\D),\ \  \eps>0,
\eqno(0.1)
$$
were considered.
Here $g(\x)$ is a periodic bounded and positive definite matrix-valued function,
and $b(\D)$ is a first order DO.
The precise assumptions on $g(\x)$ and $b(\D)$ are given below in  \S 1.

In [BSu1-4], the equation ${\A}_\eps \u_\eps + \u_\eps = \FF$ with $\FF \in L_2(\R^d;\C^n)$ was considered.
The behavior of the solution $\u_\eps$ for small $\eps$ was studied.
The solution $\u_\eps$ converges in $L_2(\R^d;\C^n)$ to the solution $\u_0$ of the
"homogenized" equation ${\A}^0 \u_0 + \u_0 = \FF$, as $\eps \to 0$.
Here ${\A}^0 = b(\D)^* g^0 b(\D)$ is the \textit{effective operator} with the constant
 effective matrix $g^0$. In [BSu1,2], it was proved that
$$
\| \u_\eps - \u_0 \|_{L_2(\R^d)} \leq C \eps \| \FF \|_{L_2(\R^d)}.
$$
In operator terms, it means that the resolvent $({\A}_\eps +I)^{-1}$
converges in the $L_2(\R^d;\C^n)$-operator norm to the resolvent of the effective operator, as $\eps \to 0$,
and we have
$$
\| ({\A}_\eps +I)^{-1} - ({\A}^0 +I)^{-1} \|_{L_2(\R^d) \to L_2(\R^d)}
\le C \eps.
\eqno(0.2)
$$

In [BSu3], more accurate approximation (including the corrector) for the resolvent
$({\A}_\eps +I)^{-1}$ in the $L_2(\R^d;\C^n)$-operator norm with an error term
$O(\eps^2)$ was obtained. (Here we do not go in details.)

In [BSu4], approximation of the resolvent of ${\A}_\eps$
in the norm of operators acting from $L_2(\R^d;\C^n)$ to the Sobolev space $H^1(\R^d;\C^n)$ was found:
$$
\| ({\A}_\eps +I)^{-1} - ({\A}^0 +I)^{-1} - \eps K(\eps) \|_{L_2(\R^d) \to H^1(\R^d)}
\le C \eps;
\eqno(0.3)
$$
this corresponds to approximation of $\u_\eps$ in the energy norm.
Here $K(\eps)$ is a corrector. It contains rapidly oscillating factors, and so depends on $\eps$.

 Estimates of the form (0.2), (0.3) are called \textit{operator error estimates}.
 They are order-sharp; the constants in estimates are controlled explicitly in terms of the
 problem data. The method of [BSu1--4] is based on the scaling transformation, the Floquet-Bloch theory
 and the analytic perturbation theory.

\smallskip\noindent\textbf{0.2. A different approach} to operator error estimates in homogenization
problems was suggested by  V.~V.~Zhikov. In [Zh1, Zh2, ZhPas, Pas], the scalar elliptic operator
$- \div g(\x/\eps) \nabla$ (where $g(\x)$ is a matrix with real entries)
and the system of elasticity theory were studied. Estimates of the form
(0.2) and (0.3) for the corresponding problems in $\R^d$ were obtained.
The method was based on analysis of the first order approximation to the solution and introduction
 of an additional parameter.
Besides the problems in $\R^d$, homogenization problems in a bounded domain $\O \subset \R^d$
with the Dirichlet or Neumann boundary condition were studied.
Approximation of the solutions in $H^1(\O)$ was deduced from the corresponding results in $\R^d$.
Due to the boundary influence, estimates in a bounded domain become worse and the error term is
$O(\eps^{1/2})$:
$$
\| \A_{\flat,\eps}^{-1} - (\A^0_\flat)^{-1} - \eps K_\flat(\eps)\|_{L_2(\O)\to H^1(\O)} \le C \eps^{1/2},\ \ \flat=D,N.
\eqno(0.4)
$$
Here $\A_{D,\eps}$ and $\A_{N,\eps}$ are operators with the Dirichlet or Neumann boundary conditions,
$\A^0_D$, $\A_N^0$ and $K_D(\eps)$, $K_N(\eps)$ are the corresponding effective operators and correctors,
respectively.

In [ZhPas], an estimate
$\| \A_{\flat,\eps}^{-1} - (\A^0_\flat)^{-1}\|_{L_2\to L_2} \le C \eps^{1/2}$ was obtained as a (rough)
consequence of (0.4).
Refinement of this estimate is a natural problem.
In [ZhPas], for the case of the scalar elliptic operator
$- \div g(\x/\eps) \nabla$ (where $g(\x)$ is a matrix with real entries)
with the Dirichlet boundary condition an estimate for
$\|\A_{D,\eps}^{-1} - (\A^0_D)^{-1} \|_{L_2\to L_2}$ of order $\eps^{\frac{d}{2d-2}}$
for $d \geq 3$ and of order $\eps |\log \eps|$ for $d=2$ was obtained.
The proof essentially relies on the maximum principle which is specific for scalar elliptic equations.

Similar results for the operator $- \div g(\x/\eps) \nabla$ in a bounded domain with the Dirichlet or Neumann
boundary conditions were obtained in the papers [Gr1, Gr2] by G.~Griso by the
"unfolding" method. In [Gr1] estimate of the form (0.4) was obtained, and in
[Gr2] the sharp order estimate
$$
\|\A_{\flat,\eps}^{-1} - (\A^0_\flat)^{-1} \|_{L_2(\O)\to L_2(\O)} \leq C \eps,
\quad \flat = D,N,
\eqno(0.5)
$$
was proved for the first time (for the same scalar elliptic operator).

\smallskip\noindent\textbf{0.3. Operator error estimates for matrix elliptic operators
in a bounded domain.} In the recent papers [PSu1,2], [Su1,2],
the Dirichlet problem for the equation $\A_{D,\eps}\u_\eps = \FF$ with $\FF \in L_2(\O;\C^n)$
in a bounded domain $\O \subset \R^d$ with the boundary of class $C^{1,1}$ was studied.
Here $\A_{D,\eps}$ is a matrix operator of the form (0.1) with the Dirichlet condition on $\partial \O$.
In [PSu1,2], estimate of the form (0.4) for the operator $\A_{D,\eps}$
was obtained. The method was based on using estimates (0.2), (0.3) for the problem in $\R^d$
and on estimates for the discrepancy (of the boundary layer type);
some technical tricks were borrowed from [ZhPas].

In [Su1,2], the author succeeded in proving a sharp order
 error estimate of the form (0.5) for the same matrix operator $\A_{D,\eps}$.
The problem was to obtain an estimate of order
$\eps$ for the $L_2$-norm of the discrepancy;
for this purpose, estimate of the form (0.4) was applied
and the duality arguments were used.

In a recent paper [KeLiS], homogenization problems
for uniformly elliptic systems in a bounded domain with the Dirichlet or Neumann
boundary conditions have been studied.
Under the assumptions that the coefficients are real-valued and H\"older continuous
the authors obtained estimate of the form (0.5) in the corresponding problems.

Note that the class of strongly elliptic operators (0.1) that we consider is wider
than the class of operators studied in [KeLiS].
Moreover, we do not impose any regularity assumptions on the coefficients,
which extends the area of possible applications.

\smallskip\noindent\textbf{0.4. Main results.}
In the present paper, analogs of the results from [PSu1,2], [Su1,2] for the Neumann problem are obtained.
Note that the homogenization problem with the Neumann boundary condition is more complicated than
 that with the Dirichlet condition.
Main difficulty is related to the fact that the boundary conditions in the initial and homogenized problems
are different. Indeed,
the conormal derivative of the solution contains rapidly oscillating coefficients,
while in the homogenized problem the conormal derivative has constant coefficients.
(For the Dirichlet problem, boundary conditions in the initial and homogenized problems
are one and the same.)

We study matrix DO's $\A_{N,\eps}$ in a bounded domain $\O\subset \R^d$ with the boundary of class
$C^{1,1}$. The operator $\A_{N,\eps}$ is given by the differential expression
(0.1) with the Neumann condition on $\partial \O$.
The effective operator $\A_N^0$ is given by the expression $b(\D)^* g^0 b(\D)$ with the Neumann condition.
The behavior of the solutions of the equation $(\A_{N,\eps} + \lambda I)\u_\eps =\FF$ for small $\eps$ is studied.
Here $\FF \in L_2(\O;\C^n)$, and $\lambda$ is subject to the restriction which ensures
that the operator $\A_{N,\eps} + \lambda I$ is positive definite.
The case where $\lambda=0$, which is important for applications, is studied separately in \S 9.
In this case, additional orthogonality conditions on $\FF$ and $\u_\eps$ are imposed.

In operator terms, the following estimates are obtained:
$$
\| (\A_{N,\eps} + \lambda I)^{-1} - (\A^0_N + \lambda I)^{-1} \|_{L_2(\O) \to L_2(\O)}
\leq C \eps,
\eqno(0.6)
$$
$$
\| (\A_{N,\eps} + \lambda I)^{-1} - (\A^0_N + \lambda I)^{-1} - \eps K_N(\eps) \|_{L_2(\O) \to H^1(\O)}
\leq C \eps^{1/2}.
\eqno(0.7)
$$
Here $K_N(\eps)$ is the corresponding corrector.
The form of the corrector depends on the properties of the periodic solution
$\Lambda(\x)$ of the auxiliary problem (1.7).
In general case, the corrector contains an auxiliary smoothing operator.
If $\Lambda$ is bounded, the corrector has a standard form.
Besides approximation of the solution $\u_\eps$ in $H^1(\O;\C^n)$,
we also obtain approximation of the flux $\p_\eps = g^\eps b(\D)\u_\eps$ in $L_2(\O;\C^m)$.
For a strictly interior subdomain $\O'$, the following sharp order estimate is proved:
$$
\| (\A_{N,\eps} + \lambda I)^{-1} - (\A^0_N + \lambda I)^{-1} - \eps K_N(\eps) \|_{L_2(\O) \to H^1(\O')}
\leq C \eps.
\eqno(0.8)
$$

The author considers estimate (0.6)
as the main achievement of the paper.

\smallskip\noindent\textbf{0.5. The method} is based on using estimates
(0.2), (0.3) for homogenization problem in $\R^d$ and taking the boundary influence into account.
Main difficulties are related to estimation of the "boundary layer correction term" $\w_\eps$ which is
solution of the Neumann problem for the homogeneous equation
\hbox{$\A_\eps \w_\eps + \lambda \w_\eps =0$} in $\O$ with the boundary condition
$\partial^\eps_{\bnu} \w_\eps = \brho_\eps$ on $\partial \O$.
Here $\partial^\eps_{\bnu} \w_\eps$ is the conormal derivative of $\w_\eps$,
$\brho_\eps$ is some function defined in terms of
the solution $\u_0$ of the homogenized problem and containing rapidly oscillating coefficients
(see \S 4).

Some technical tricks, in particular, using the Steklov smoothing operator,
are borrowed from [ZhPas].

First we prove estimate (0.7).
For this, it is necessary to estimate the $H^1$-norm of
$\w_\eps$ by $C \eps^{1/2} \|\FF\|_{L_2(\O)}$.
Next we obtain (0.6) by estimating the
$L_2$-norm of $\w_\eps$ in terms of $C \eps \|\FF\|_{L_2(\O)}$.
For this, we rely on the (already proved) estimate
(0.7) and use the duality arguments.

\smallskip\noindent\textbf{0.6. Plan of the paper.} The paper consists of nine sections.
In \S 1, we introduce the class of operators in $L_2(\R^d;\C^n)$
that we consider, describe the effective operator and the corrector, and give
the results from [BSu2,4] and [PSu2] needed in what follows. Precisely, we deduce theorems about approximation
of the resolvent $(\A_\eps + \lambda I)^{-1}$ for arbitrary $\lambda >0$
from the known results for the operator $(\A_\eps +I)^{-1}$.
In \S 2, the statement of the Neumann problem in a bounded domain is discussed,
some auxiliary material is given, and the homogenized problem is described.
In \S 3, different auxiliary statements are collected.
  Main results (Theorems 4.1 and 4.2) are formulated in \S 4.
  Also, in \S 4 the first part of the proof is given:
  the boundary layer correction term is introduced,
  and the proofs of Theorems 4.1 and 4.2 are reduced to the estimates for this correction term.
  \S 5 contains the proof of Theorem 4.2, and \S 6 contains the proof of Theorem 4.1.
  In \S 7 the case where $\Lambda \in L_\infty$ is studied;
in this case it is possible to get rid of the smoothing operator in the corrector.
In \S 8, estimates in a strictly interior subdomain
$\O'$ of the domain $\O$ are obtained:
it is shown that estimate (0.6) and the results in $\R^d$ imply estimate (0.8).
Finally, \S 9 is devoted to the case where $\lambda =0$: the corresponding
results are deduced from the theorems for the case where $\lambda >0$.

\smallskip\noindent\textbf{0.7. Notation.} Let $\H$ and $\H_*$ be complex separable Hilbert spaces.
The symbols $(\cdot,\cdot)_{\H}$ and $\|\cdot\|_{\H}$ stand for the inner product and the norm in
$\H$; the symbol $\|\cdot\|_{\H \to \H_*}$ stands for the norm of a linear continuous operator
acting from $\H$ to $\H_*$.
By $I=I_\H$ we denote the identity operator in $\H$.

The symbols $\langle \cdot, \cdot \rangle$ and $|\cdot|$ stand for the inner product and the norm in
$\C^n$; $\1 = \1_n$ is the identity $(n\times n)$-matrix.
If $a$ is an $(n\times n)$-matrix, the symbol $|a|$ denotes the norm of $a$
as an operator in $\C^n$. We denote $\x = (x_1,\dots,x_d)\in \R^d$, $iD_j = \partial_j = \partial/\partial x_j$,
$j=1,\dots,d$, $\D = -i \nabla = (D_1,\dots,D_d)$.
The $L_p$-classes of $\C^n$-valued functions
in a domain ${\mathcal O} \subset \R^d$ are denoted by
$L_p({\mathcal O};\C^n)$, $1 \le p \leq \infty$.
The Sobolev classes of $\C^n$-valued functions in a domain ${\mathcal O} \subset \R^d$
are denoted by $H^s({\mathcal O};\C^n)$.
Next, $H^1_0(\O;\C^n)$ is the closure of $C_0^\infty(\O;\C^n)$ in $H^1(\O;\C^n)$.
If $n=1$, we write simply $L_p({\mathcal O})$, $H^s({\mathcal O})$, etc., but sometimes
we use such abbreviated notation also for the spaces of vector-valued or matrix-valued functions.

\smallskip\noindent\textbf{0.8. Acknowledgements.} The results of the paper have been obtained and
the paper was written during the author's visit to the Mittag-Leffler Institute (Sweden)
in September and October 2012.
The author is grateful to the director of the Institute A.~Laptev, to the organizers of the program
G.~Rozenblum and G.~Raikov, to the team of the Institute
for warm hospitality and creation of excellent conditions for scientific research.

\section*{\S 1. Homogenization problem for a periodic elliptic operator in
 $L_2(\mathbb{R}^d;\mathbb{\C}^n)$}

 In this section, we describe the class of matrix elliptic operators
 that we consider and give the results for homogenization problem in $\R^d$ obtained in [BSu2], [BSu4], and also in  [PSu2]. Precisely, from the known results
 about approximation of the resolvent $(\A_\eps + I)^{-1}$ we deduce
theorems about approximation of the resolvent $(\A_\eps +\lambda I)^{-1}$ for any $\lambda>0$.

\smallskip\noindent\textbf{1.1. Lattices in $\R^d$.}
Let $\Gamma \subset \R^d$ be a lattice generated by the basis $\a_1,\dots,\a_d \in \R^d$:
$$
\Gamma = \{ \a \in \R^d:\ \a = \sum_{j=1}^d \nu_j \a_j,\ \ \nu_j \in \Z\},
$$
and let $\Omega$ be the (elementary) cell of the lattice $\Gamma$:
$$
\Omega := \{\x \in \R^d:\ \x = \sum_{j=1}^d \tau_j \a_j,\ \ -\frac{1}{2} < \tau_j < \frac{1}{2} \}.
$$
We use the notation $|\Omega| = \meas \Omega$.
The basis $\b_1,\dots, \b_d$ in $\R^d$ dual to $\a_1,\dots, \a_d$ is defined
by the relations $\langle \b_i, \a_j\rangle  = 2 \pi \delta_{ij}$.
This basis generates a lattice $\wt{\Gamma}$ dual to the lattice $\Gamma$.
Below we use the notation
$$
r_0 = \frac{1}{2} \min_{0 \ne \b \in \wt{\Gamma}} |\b|,\quad r_1 = \frac{1}{2} \text{diam}\,\Omega.
$$

By $\wt{H}^1(\Omega)$ we denote the subspace of all functions in
$H^1(\Omega)$ whose $\Gamma$-periodic extension to $\R^d$ belongs to $H^1_{\text{loc}}(\R^d)$.
If $\varphi(\x)$ is a $\Gamma$-periodic function in $\R^d$, we denote
$$
\varphi^\eps(\x) := \varphi(\eps^{-1}\x),\quad \eps >0.
$$

\smallskip\noindent\textbf{1.2. The class of operators.}
In $L_2(\R^d;\C^n)$, consider a second order DO $\A_\eps$ formally given by the differential expression
$$
\mathcal{A}_\eps = b(\D)^* g^\eps(\x) b (\D),\ \ \eps >0.
\eqno(1.1)
$$
Here $g(\x)$ is a measurable Hermitian \hbox{$(m \times m)$}-matrix-valued
function (in general, with complex entries). We assume that $g(\x)$ is periodic with respect to
the lattice $\Gamma$, bounded and uniformly positive definite, i.~e.,
$$
c \1_m \le g(\x) \le \wt{c} \1_m,\ \ \x \in \R^d;\ \ 0< c \le \wt{c} <\infty.
$$
Next, $b(\D)$ is a homogeneous $(m\times n)$-matrix first order DO
with constant coefficients:
$$
b(\D)=\sum_{l=1}^d b_l D_l.
\eqno(1.2)
$$
Here $b_l$ are constant matrices (in general, with complex entries).
The symbol $b(\boldsymbol{\xi}) = \sum_{l=1}^d b_l \xi_l$, $\boldsymbol{\xi} \in \mathbb{R}^d$,
corresponds to the operator $b(\D)$.
It is assumed that $m \ge n$ and
$$
\textrm{rank} \, b (\boldsymbol{\xi}) = n, \ 0 \ne \boldsymbol{\xi} \in \R^d.
\eqno(1.3)
$$
This condition is equivalent to the following inequalities
$$
\alpha_0 \mathbf{1}_n \leq b(\boldsymbol{\theta})^* b(\boldsymbol{\theta}) \leq \alpha_1 \mathbf{1}_n,
\quad \boldsymbol{\theta} \in \mathbb{S}^{d-1}, \quad 0 < \alpha_0 \leq \alpha_1 < \infty,
\eqno(1.4)
$$
with some positive constants $\alpha_0$ and $\alpha_1$.
From (1.4) it follows that
$$
|b_l| \le \alpha_1^{1/2},\quad l=1,\dots,d.
\eqno(1.5)
$$

The precise definition of the operator $\mathcal{A}_\eps$ is given in terms of
the corresponding quadratic form
$$
a_\eps[\u,\u] = \int_{\mathbb{R}^d} \left\langle g^\eps (\x) b (\D) \u, b(\D) \u \right\rangle \,
d \x, \quad \u \in H^1 (\mathbb{R}^{d}; \mathbb{C}^{n}).
$$
Under the above assumptions this form is closed in $L_2(\R^d;\C^n)$ and nonnegative.
Using the Fourier transformation and condition (1.4), it is easy to check that
$$
c_0 \int_{\R^d} |\D \u|^2\, d\x \le a_\eps[\u, \u] \le c_1 \int_{\R^d} |\D \u|^2\, d\x,
\ \ \u \in H^1(\R^d;\C^n),
\eqno(1.6)
$$
where $c_0 = \alpha_0 \|g^{-1}\|^{-1}_{L_\infty}$, $c_1 = \alpha_1 \|g\|_{L_\infty}$.

The simplest example of the operator (1.1) is the scalar elliptic operator
$
\A_\eps = -\div g^\eps(\x) \nabla = \D^* g^\eps(\x)\D.
$
In this case we have $n=1$, $m=d$, $b(\D)=\D$. Obviously, condition (1.4) is valid with $\alpha_0=\alpha_1=1$.
Another example is the operator of elasticity theory which can be written in the form (1.1)
with $n=d$ and $m=d(d+1)/2$. These and other examples are discussed in [BSu2] in detail.

\smallskip\noindent\textbf{1.3. The effective operator.}
In order to formulate the results, we need to introduce the effective operator $\A^0$.

Let an $(n\times m)$-matrix-valued function $\Lambda(\x)$ be the (weak) $\Gamma$-periodic solution
of the problem
$$
b(\D)^* g (\x)\left( b(\D) \Lambda(\x) + \mathbf{1}_m \right) = 0,
\quad \int_{\Omega} \Lambda(\x) \, d \x = 0.
\eqno(1.7)
$$
In other words, for the columns $\v_j(\x)$, $j=1,\dots,m,$ of the matrix $\Lambda(\x)$
the following is true: $\v_j \in \wt{H}^1(\Omega;\C^n)$, we have
$$
\int_\Omega \langle g(\x) (b(\D)\v_j(\x) + \wh{\e}_j), b(\D) \boldsymbol{\eta}(\x)\rangle\,d\x =0,
\quad \forall \boldsymbol{\eta} \in \wt{H}^1(\Omega;\C^n),
$$
and $\int_\Omega \v_j(\x)\,d\x=0$. Here $\wh{\e}_1,\dots,\wh{\e}_m$ is the standard orthonormal basis in $\C^m$.

The so-called \emph{effective matrix} $g^0$ of size $m\times m$ is defined as follows:
$$
g^0 = |\Omega|^{-1}
\int_{\Omega} \wt{g}(\x)\, d \x,
\eqno(1.8)
$$
where
$$
\wt{g}(\x):= g (\x) \left( b (\D) \Lambda (\x) + \mathbf{1}_m \right).
\eqno(1.9)
$$
It turns out that the matrix $g^0$ is positive definite.
The \textit{effective operator} $\A^0$ for the operator (1.1) is given by the differential expression
$$
\A^0 = b (\D)^* g^0 b (\D)
$$
on the domain $H^2(\R^d;\C^n)$.

Below we need the following estimates for $\Lambda(\x)$ proved in [BSu3, (6.28) and Subsection 7.3]:
$$
\|\D \Lambda\|_{L_2(\Omega)} \le |\Omega|^{1/2} m^{1/2} \alpha_0^{-1/2} \|g\|_{L_\infty}^{1/2}
\|g^{-1}\|_{L_\infty}^{1/2},
\eqno(1.10)
$$
$$
\| \Lambda\|_{L_2(\Omega)} \le |\Omega|^{1/2} m^{1/2} (2r_0)^{-1} \alpha_0^{-1/2} \|g\|_{L_\infty}^{1/2}
\|g^{-1}\|_{L_\infty}^{1/2}.
\eqno(1.11)
$$

\smallskip\noindent\textbf{1.4. Properties of the effective matrix}.
The following properties of the effective matrix are proved in [BSu2, Chapter 3, Theorem 1.5].

\smallskip\noindent\textbf{Proposition 1.1.} \textit{The effective matrix satisfies the estimates}
$$
\underline{g} \le g^0 \le \overline{g}.
\eqno(1.12)
$$
\textit{Here}
$$
\overline{g}= |\Omega|^{-1} \int_\Omega g(\x)\,d\x, \quad
\underline{g}= \left(|\Omega|^{-1} \int_\Omega g(\x)^{-1}\,d\x\right)^{-1}.
$$
\textit{If} $m = n$, \textit{the effective matrix} $g^0$ \textit{coincides with} $\underline{g}$.

\smallskip
In homogenization theory for specific DO's, estimates (1.12) are known as the Voight-Reuss bracketing.
We distinguish the cases where one of the inequalities in (1.12) becomes an identity.
The following statements were obtained in [BSu2, Chapter 3, Propositions 1.6 and 1.7].

\smallskip\noindent\textbf{Proposition 1.2.} \textit{The identity} $g^0 = \overline{g}$
\textit{is equivalent to the relations}
$$
b(\D)^* {\mathbf g}_k(\x) = 0, \ \ k = 1,\dots,m,
\eqno(1.13)
$$
\textit{where} ${\mathbf g}_k(\x)$, $k = 1,\dots, m$, \textit{are the columns of the matrix} $g(\x)$.

\smallskip\noindent\textbf{Proposition 1.3.} \textit{The identity} $g^0 = \underline{g}$
\textit{is equivalent to the relations}
$$
 {\mathbf l}_k(\x) = {\mathbf l}_k^0 + b(\D) \w_k, \ \ {\mathbf l}_k^0 \in \C^m,
 \ \ \w_k \in \wt{H}^1(\Omega;\C^n),\ \ k = 1,\dots,m,
\eqno(1.14)
$$
 \textit{where} ${\mathbf l}_k(\x)$, $k= 1,\dots, m$, \textit{are the columns of the matrix} $g(\x)^{-1}$.

\smallskip
Obviously, (1.12) implies the following estimates for the norms of the matrices $g^0$ and $(g^0)^{-1}$:
$$
| g^0| \le \|g\|_{L_\infty},\ \ | (g^0)^{-1}| \le \|g^{-1}\|_{L_\infty}.
\eqno(1.15)
$$
Note that, by (1.4) and (1.15), the symbol of the effective operator $\A^0$
satisfies the following inequality:
$$
b(\bxi)^* g^0 b(\bxi) \ge c_0 |\bxi|^2 \1_n,\quad \bxi \in \R^d,\quad c_0 = \alpha_0 \|g^{-1}\|^{-1}_{\infty}.
\eqno(1.16)
$$

\smallskip\noindent\textbf{1.5. Smoothing in Steklov's sense.}
Let $S_\eps$ be the operator in $L_2(\R^d;\C^m)$ given by
$$
(S_\eps \u)(\x) = |\Omega|^{-1} \int_\Omega \u(\x - \eps \z)\, d\z
\eqno(1.17)
$$
and called the \textit{Steklov smoothing operator}.
Note that
$$
\|S_\eps\|_{L_2(\R^d)\to L_2(\R^d)} \leq 1.
\eqno(1.18)
$$
Obviously, $D^\alpha S_\eps \u = S_\eps D^\alpha \u$ for
$\u \in H^s(\R^d;\C^m)$ and $|\alpha| \leq s$.
Therefore,
$$
\|S_\eps\|_{H^s(\R^d)\to H^s(\R^d)} \leq 1,\ \ s \in {\mathbb N}.
\eqno(1.19)
$$

We mention some properties of the operator (1.17), see [ZhPas, Lemmas 1.1 and 1.2]
or [PSu2, Propositions 3.1, 3.2].

\smallskip\noindent\textbf{Proposition 1.4.}
\textit{For any} $\u \in H^1(\R^d;\C^m)$ \textit{we have}
$$
\| S_\eps \u - \u\|_{L_2(\R^d)} \leq \eps r_1 \| \D \u \|_{L_2(\R^d)}, \ \ \eps>0.
$$

\smallskip\noindent\textbf{Proposition 1.5.} \textit{Let} $f(\x)$
\textit{be a} $\Gamma$-\textit{periodic function in} $\R^d$ \textit{such that} $f \in L_2(\Omega)$.
\textit{Let $[f^\eps]$ be the operator of multiplication by the function $f^\eps(\x)$. Then the operator} $[f^\eps]S_\eps$ \textit{is continuous in} $L_2(\R^d;\C^m)$, \textit{and}
$$
\| [f^\eps]S_\eps \|_{L_2(\R^d)\to L_2(\R^d)}
\leq |\Omega|^{-1/2} \| f \|_{L_2(\Omega)},\ \ \eps>0.
$$

\smallskip\noindent\textbf{1.6. The results for homogenization problem in $\R^d$}.
Consider the following elliptic equation in $\R^d$:
$$
\A_\varepsilon \u_{\varepsilon} + \lambda \u_{\varepsilon} = \FF,
\eqno(1.20)
$$
where $\FF \in L_2(\R^d;\C^n)$, and $\lambda >0$ is a parameter.
It is known that the solution $\u_{\varepsilon}$ converges in
$L_2 (\R^d; \C^n)$ to the solution of the "homogenized" equation
$$
{\A}^0 \u_0 + \lambda \u_0 = \FF,
\eqno(1.21)
$$
as $\varepsilon \to 0$.

\smallskip\noindent\textbf{Theorem 1.6.} \textit{Let} $\u_\eps$ \textit{be the solution of the equation}
(1.20), \textit{and let} $\u_0$ \textit{be the solution of the equation} (1.21). \textit{Then}
$$
 \| \u_{\varepsilon} - \u_0 \|_{L_2 (\mathbb{R}^{d}; \mathbb{C}^{n})}
 \leq C_1(\lambda) \varepsilon  \| \FF \|_{L_2 (\mathbb{R}^{d}; \mathbb{C}^{n})},
 \quad \varepsilon >0,
 $$
\textit{or, in operator terms},
$$
\| (\A_{\varepsilon} + \lambda I)^{-1} - (\A^0 + \lambda I)^{-1}
\|_{L_2 (\mathbb{R}^{d}; \mathbb{C}^{n}) \rightarrow L_2 (\mathbb{R}^{d}; \mathbb{C}^{n})}
\leq \ C_1(\lambda) \varepsilon, \quad \varepsilon >0.
\eqno(1.22)
$$
\textit{Here $C_1(\lambda)=\check{C}_1\lambda^{-1/2}$,
and the constant $\check{C}_1$ depends only on the norms}
$\| g \|_{L_\infty}$, $\| g^{-1} \|_{L_\infty}$, \textit{on the constants} $\alpha_0$, $\alpha_1$
\textit{from} (1.4) \textit{and the parameters of the lattice} $\Gamma$.

\smallskip
\noindent\textbf{Proof.}
In [BSu2, Chapter 4, Theorem 2.1], estimate (1.22) was proved in the case where
$\lambda =1$ for $0< \eps \le 1$.
We only have to explain how this estimate is carried over to the general case.

Note that for $\lambda=1$ and $\eps >1$ the left-hand side of (1.22) is obviously less than
2, and then also less than 2$\eps$. Therefore, we start with the estimate
$$
\| (\A_{\varepsilon} + I)^{-1} - (\A^0 + I)^{-1}
\|_{L_2 (\mathbb{R}^{d}) \rightarrow L_2 (\mathbb{R}^{d})}
\leq \ \check{C}_1 \varepsilon, \quad \varepsilon >0.
\eqno(1.23)
$$
Here the constant $\check{C}_1$ depends only on
$\| g \|_{L_\infty}$, $\| g^{-1} \|_{L_\infty}$, $\alpha_0$, $\alpha_1$, and the parameters of the lattice $\Gamma$.

Next, by the scaling transformation, (1.23) is equivalent to the inequality
$$
\| (\A + \eps^2 I)^{-1} - (\A^0 + \eps^2 I)^{-1}
\|_{L_2 (\mathbb{R}^{d}) \rightarrow L_2 (\mathbb{R}^{d})}
\leq \ \check{C}_1  \varepsilon^{-1}, \quad \eps >0.
\eqno(1.24)
$$
Here $\A = b(\D)^* g(\x)b(\D)$. Using (1.24) with $\eps$ replaced by $\eps \lambda^{1/2}$
and applying the inverse transformation, we arrive at (1.22) with the constant
$C_1(\lambda)=\check{C}_1\lambda^{-1/2}$.
$\ \bullet$

\smallskip
In order to approximate $\u_\eps$ in $H^1(\R^d;\C^n)$ it is necessary to take
the first order corrector into account. We put
$$
K_\lambda(\eps) = [\Lambda^{\varepsilon}] S_{\varepsilon} b (\D) (\A^0 + \lambda I)^{-1}.
\eqno(1.25)
$$
Here $[\Lambda^{\varepsilon}]$
is the operator of multiplication by the matrix-valued function $\Lambda(\eps^{-1}\x)$, and
$S_\eps$ is the smoothing operator defined by (1.17).
The operator (1.25) is a continuous mapping of $L_2(\R^d;\C^n)$ into $H^1(\R^d;\C^n)$.
Indeed, the operator $b (\D) (\A^0 + \lambda I)^{-1}$ is continuous from
$L_2(\R^d;\C^n)$ to $H^1(\R^d;\C^m)$. The operator $[\Lambda^\eps]S_\eps$
is a continuous mapping of $H^1(\R^d;\C^m)$ into $H^1(\R^d;\C^n)$.
It can be easily checked by using Proposition 1.5 and relation $\Lambda \in \wt{H}^1(\Omega)$.
We have $\eps \|K_\lambda(\eps)\|_{L_2 \to H^1} =O(1)$ for small $\eps$.
(See the proof of Theorem 1.7 below, where this is checked for $\lambda=1$.)

"The first order approximation" to the solution $\u_\eps$ is given by
$$
\v_{\varepsilon} = \u_0 + \varepsilon \Lambda^{\varepsilon} S_{\varepsilon} b (\D) \u_0=
(\A^0+\lambda I)^{-1} \FF + \eps K_\lambda(\eps) \FF.
\eqno(1.26)
$$

\smallskip\noindent\textbf{Theorem 1.7.} \textit{Let} $\u_\eps$ \textit{be the solution of the equation}
(1.20), \textit{and let} $\u_0$ \textit{be the solution of the equation} (1.21). \textit{Suppose that}
$\v_\eps$ \textit{is defined by} (1.26). \textit{Then}
$$
 \| \u_{\varepsilon} - \v_\eps \|_{H^1 (\mathbb{R}^{d}; \mathbb{C}^{n})}
 \leq C_2(\lambda) \varepsilon  \| \FF \|_{L_2 (\mathbb{R}^{d}; \mathbb{C}^{n})},
 \quad  \varepsilon >0,
 \eqno(1.27)
 $$
\textit{or, in operator terms},
$$
\| (\A_{\varepsilon} + \lambda I)^{-1} - (\A^0 + \lambda I)^{-1} -
\varepsilon K_\lambda (\varepsilon) \|_{L_2 (\mathbb{R}^{d}; \mathbb{C}^{n})
\to H^1 (\mathbb{R}^{d}; \mathbb{C}^{n})} \leq \ C_2(\lambda) \varepsilon,\ \eps>0.
\eqno(1.28)
$$
\textit{Here $C_2(\lambda)= \wh{C}_2 (\lambda^{-1/2}+1)$,
and the constant $\wh{C}_2$ depends only on}
$d, m, \alpha_0, \alpha_1, \| g \|_{L_\infty}, \| g^{-1} \|_{L_\infty}$,
\textit{and the parameters of the lattice} $\Gamma$.

\smallskip
\noindent\textbf{Proof.}
In [BSu4, Theorem 10.6], a similar result was proved for $\lambda=1$, but with a different smoothing
operator instead of $S_\eps$. In [PSu2, Theorem 3.3], it was shown that it is possible to
pass to the smoothing operator $S_\eps$,
and (1.28) was proved for $\lambda=1$ and $0< \eps \le 1$.
We only have to explain how this estimate can be carried over to the general case.

Thus, we start with the estimate
$$
\| (\A_{\varepsilon} + I)^{-1} - (\A^0 + I)^{-1} -
\varepsilon K_1 (\varepsilon) \|_{L_2 (\mathbb{R}^{d})
\to H^1 (\mathbb{R}^{d})} \leq \ \check{C}_2 \varepsilon,\ 0< \eps\le 1.
\eqno(1.29)
$$
The constant $\check{C}_2$ depends only on
$d, m, \alpha_0, \alpha_1, \| g \|_{L_\infty}, \| g^{-1} \|_{L_\infty}$, and the parameters of the lattice $\Gamma$.

For $\eps >1$ estimates are trivial: each operator under the norm sign in (1.29)
is estimated separately. From the lower inequality (1.6) it follows that
$$
\min \{c_0,1\} \|(\A_\eps +I)^{-1}\v\|^2_{H^1(\R^d)} \le ((\A_\eps +I)^{-1}\v,\v)_{L_2(\R^d)}
\le \|\v\|^2_{L_2(\R^d)}
$$
for $\v \in L_2(\R^d;\C^n)$, whence
$$
\|(\A_\eps +I)^{-1}\|_{L_2(\R^d) \to H^1(\R^d)} \le \max \{1,c_0^{-1/2}\} =
\max \{1, \alpha_0^{-1/2} \|g^{-1}\|_{L_\infty}^{1/2}\}.
\eqno(1.30)
$$
By (1.15), the norm of $(\A^0 + I)^{-1}$ satisfies the same estimate:
$$
\|(\A^0 +I)^{-1}\|_{L_2(\R^d) \to H^1(\R^d)} \le
\max \{1, c_0^{-1/2} \}.
\eqno(1.31)
$$

Now we estimate the $(L_2 \to H^1)$-norm of the operator $\eps K_1(\eps)=
\eps [\Lambda^\eps]S_\eps b(\D)(\A^0+I)^{-1}$.
Let $\FF \in L_2(\R^d;\C^n)$. Then, by Proposition 1.5,
$$
\| \eps K_1(\eps) \FF \|_{L_2(\R^d)} \le \eps |\Omega|^{-1/2} \|\Lambda\|_{L_2(\Omega)}
\| b(\D) (\A^0 +I)^{-1}\FF\|_{L_2(\R^d)}.
$$
Using the Fourier transformation and (1.4), (1.16), we obtain:
$$
\begin{aligned}
&\| b(\D) (\A^0 +I)^{-1}\FF\|^2_{L_2(\R^d)}
\le
\int_{\R^d} |b(\bxi)|^2 |(b(\bxi)^* g^0 b(\bxi) + \1)^{-1}|^2 |\wh{\FF}(\bxi)|^2\,d\bxi
\cr
&\le \alpha_1 \int_{\R^d} |\bxi|^2 (c_0 |\bxi|^2 + 1)^{-2} |\wh{\FF}(\bxi)|^2\,d\bxi
\le \alpha_1 (2c_0)^{-1} \|\FF\|^2_{L_2(\R^d)}.
\end{aligned}
\eqno(1.32)
$$
Here $\wh{\FF}(\bxi)$ is the Fourier image of $\FF(\x)$.
Hence,
$$
\| \eps K_1(\eps) \FF \|_{L_2(\R^d)} \le \eps |\Omega|^{-1/2} \|\Lambda\|_{L_2(\Omega)}
\alpha_1^{1/2} (2c_0)^{-1/2} \|\FF\|_{L_2(\R^d)}.
\eqno(1.33)
$$
Consider the derivatives
$$
\eps \partial_j (K_1(\eps) \FF) = [(\partial_j \Lambda)^\eps] S_\eps  b(\D) (\A^0 +I)^{-1}\FF +
\eps [\Lambda^\eps] S_\eps  b(\D) \partial_j (\A^0 +I)^{-1}\FF.
$$
By Proposition 1.5, we have:
$$
\begin{aligned}
\sum_{j=1}^d &\|\eps \partial_j (K_1(\eps) \FF)\|^2_{L_2(\R^d)}
\le 2 |\Omega|^{-1} \| \D \Lambda\|^2_{L_2(\Omega)}
\| b(\D) (\A^0 +I)^{-1}\FF\|^2_{L_2(\R^d)}
\cr
&+ 2 \eps^2 |\Omega|^{-1} \| \Lambda\|^2_{L_2(\Omega)}
\sum_{j=1}^d \| b(\D) \partial_j (\A^0 +I)^{-1}\FF\|^2_{L_2(\R^d)}.
\end{aligned}
\eqno(1.34)
$$
Similarly to (1.32),
$$
\begin{aligned}
\sum_{j=1}^d &\| b(\D) \partial_j (\A^0 +I)^{-1}\FF\|^2_{L_2(\R^d)}
\cr
&\le \alpha_1 \int_{\R^d} |\bxi|^4 (c_0 |\bxi|^2 + 1)^{-2} |\wh{\FF}(\bxi)|^2\,d\bxi
\le \alpha_1 c_0^{-2} \|\FF\|^2_{L_2(\R^d)}.
\end{aligned}
\eqno(1.35)
$$
As a result, from (1.32), (1.34), and (1.35) it follows that
$$
\begin{aligned}
\sum_{j=1}^d \|\eps \partial_j (K_1(\eps) \FF)\|^2_{L_2(\R^d)}
&\le
 |\Omega|^{-1} \| \D \Lambda\|^2_{L_2(\Omega)} \alpha_1 c_0^{-1} \| \FF\|^2_{L_2(\R^d)}
\cr
&+ 2 \eps^2 |\Omega|^{-1} \| \Lambda\|^2_{L_2(\Omega)} \alpha_1 c_0^{-2} \| \FF\|^2_{L_2(\R^d)}.
\end{aligned}
\eqno(1.36)
$$
By (1.33) and (1.36),
$$
\begin{aligned}
\| \eps K_1(\eps)\|_{L_2(\R^d) \to H^1(\R^d)}
&\le
|\Omega|^{-1/2} \alpha_1^{1/2}  c_0^{-1/2} \| \D \Lambda\|_{L_2(\Omega)}
\cr
&+  \eps |\Omega|^{-1/2} \alpha_1^{1/2} \| \Lambda\|_{L_2(\Omega)}
(2c_0)^{-1/2} (1 + 4 c_0^{-1})^{1/2}.
\end{aligned}
$$
Combining this with (1.10) and (1.11), we obtain
$$
\| \eps K_1(\eps)\|_{L_2(\R^d) \to H^1(\R^d)}
\le C_3 + C_4 \eps,\quad \eps >0,
\eqno(1.37)
$$
where
$$
\begin{aligned}
C_3 &= m^{1/2} \alpha_1^{1/2} \alpha_0^{-1} \|g\|^{1/2}_{L_\infty} \|g^{-1}\|_{L_\infty},
\cr
C_4 &= 2^{-3/2}
m^{1/2} r_0^{-1} \alpha_1^{1/2} \alpha_0^{-1} \|g\|^{1/2}_{L_\infty} \|g^{-1}\|_{L_\infty}
\left( 1+ 4 \alpha_0^{-1} \|g^{-1}\|_{L_\infty}\right)^{1/2}.
\end{aligned}
$$

Relations (1.30), (1.31), and (1.37) imply that
$$
\begin{aligned}
\|& (\A_{\varepsilon} + I)^{-1} - (\A^0 + I)^{-1} -
\varepsilon K_1 (\varepsilon) \|_{L_2 (\mathbb{R}^{d})
\to H^1 (\mathbb{R}^{d})}
\cr
&\leq 2 \max\{ 1, \alpha_0^{-1/2}\|g^{-1}\|_{L_\infty}^{1/2} \} + C_3 + C_4 \eps,\ \  \eps >0.
\end{aligned}
$$
Obviously, for $\eps >1$ the right-hand side does not exceed $\wt{C}_2 \eps$, where
$\wt{C}_2 = 2 \max\{ 1, \alpha_0^{-1/2}\|g^{-1}\|_{L_\infty}^{1/2}\} + C_3 + C_4$.
Combining this with (1.29), we obtain estimate of the form (1.29) for all $\eps >0$:
$$
\| (\A_{\varepsilon} + I)^{-1} - (\A^0 + I)^{-1} -
\varepsilon K_1 (\varepsilon) \|_{L_2 (\mathbb{R}^{d})
\to H^1 (\mathbb{R}^{d})} \leq  \wh{C}_2 \eps,\ \  \eps >0,
\eqno(1.38)
$$
where $\wh{C}_2 = \max \{\check{C}_2, \wt{C}_2\}$.

A similar approximation for the operator
$(\A_\eps +\lambda I)^{-1}$ with arbitrary $\lambda >0$ can be easily deduced from (1.38).
Indeed, (1.38) is equivalent to
$$
\| (-\Delta + I)^{1/2} \left((\A_{\varepsilon} + I)^{-1} - (\A^0 +I)^{-1} -
\varepsilon K_1 (\varepsilon)\right) \|_{L_2 (\mathbb{R}^{d})
\to L_2 (\mathbb{R}^{d})} \leq  \wh{C}_2 \varepsilon
\eqno(1.39)
$$
for $\eps>0$. By the scaling transformation, (1.39) is equivalent to
$$
\begin{aligned}
\|& (-\Delta + \eps^2 I)^{1/2} \left((\A + \eps^2 I)^{-1} - (\A^0 + \eps^2 I)^{-1} -
\Lambda S b(\D) (\A^0 + \eps^2 I)^{-1} \right) \|_{L_2 \to L_2}
\cr
&\leq \wh{C}_2,\quad \eps>0.
\end{aligned}
\eqno(1.40)
$$
Here $S=S_1$ is the operator (1.17) with $\eps =1$.
From (1.40) with $\eps$ replaced by $\eps \lambda^{1/2}$, applying the inverse transformation, we obtain:
$$
\| (-\Delta + \lambda I)^{1/2} \left((\A_{\varepsilon} + \lambda I)^{-1} - (\A^0 + \lambda I)^{-1} -
\varepsilon K_\lambda (\varepsilon)\right) \|_{L_2 (\mathbb{R}^{d})
\to L_2 (\mathbb{R}^{d})} \leq \wh{C}_2 \varepsilon
$$
for $\eps>0$. This implies (1.28) with
$C_2(\lambda) = \wh{C}_2 (\lambda^{-1/2}+ 1)$. $\ \bullet$

\smallskip
Theorem 1.7 allows one to obtain approximation for the flux
$\p_\eps := g^\eps b(\D) \u_\eps$ in $L_2(\R^d;\C^m)$.

\smallskip\noindent\textbf{Theorem 1.8.} \textit{Let} $\u_\eps$ \textit{be the solution of the equation}
(1.20), \textit{and let} $\u_0$ \textit{be the solution of the equation} (1.21). \textit{Let
$\p_\eps := g^\eps b(\D) \u_\eps$. Then we have}
$$
 \| \p_{\varepsilon} - \wt{g}^\eps S_\eps b(\D)\u_0 \|_{L_2(\R^{d}; \mathbb{C}^{m})}
 \leq C_5(\lambda) \varepsilon  \| \FF \|_{L_2 (\mathbb{R}^{d}; \mathbb{C}^{n})},
 \quad  \varepsilon >0.
 \eqno(1.41)
 $$
\textit{Here $\wt{g}(\x)$ is the matrix} (1.9),
$C_5(\lambda)= C_5' \lambda^{-1/2}+C_5''$, \textit{and the constants $C_5'$, $C_5''$ depend only on}
$d, m, \alpha_0, \alpha_1, \| g \|_{L_\infty}, \| g^{-1} \|_{L_\infty}$,
\textit{and the parameters of the lattice} $\Gamma$.

\smallskip
\noindent\textbf{Proof.}
By (1.4) and (1.27),
$$
 \| \p_{\varepsilon} - {g}^\eps b(\D)\v_\eps \|_{L_2(\R^{d})}
 \leq \eps \alpha_1^{1/2} \|g\|_{L_\infty} C_2(\lambda)
 \| \FF \|_{L_2 (\mathbb{R}^{d})},
 \quad  \varepsilon >0.
 \eqno(1.42)
 $$
From (1.2) and (1.26) it follows that
$$
{g}^\eps b(\D)\v_\eps =
g^\eps b(\D) \u_0 + g^\eps(b(\D)\Lambda)^\eps S_\eps b(\D) \u_0 +
\eps \sum_{l=1}^d g^\eps b_l \Lambda^\eps S_\eps b(\D) D_l \u_0.
\eqno(1.43)
$$
Using (1.5) and Proposition 1.5, we estimate the last term in the right-hand side of (1.43):
$$
\begin{aligned}
&\eps \left\| \sum_{l=1}^d g^\eps b_l \Lambda^\eps S_\eps b(\D) D_l \u_0 \right\|_{L_2(\R^d)}
\cr
&\le \eps \|g\|_{L_\infty} \alpha_1^{1/2} |\Omega|^{-1/2} \| \Lambda \|_{L_2(\Omega)}
d^{1/2} \| \D b(\D) \u_0\|_{L_2(\R^d)}.
\end{aligned}
\eqno(1.44)
$$
Next, by Proposition 1.4,
$$
\| g^\eps b(\D) \u_0 -  g^\eps S_\eps b(\D) \u_0 \|_{L_2(\R^d)}
\le \eps \|g\|_{L_\infty} r_1 \|\D b(\D)\u_0\|_{L_2(\R^d)}.
\eqno(1.45)
$$
By (1.9), we have
$$
g^\eps S_\eps b(\D)\u_0 +g^\eps (b(\D)\Lambda)^\eps S_\eps b(\D)\u_0=
\wt{g}^\eps S_\eps b(\D)\u_0.
\eqno(1.46)
$$
Combining (1.43)--(1.46) and (1.11), we arrive at
$$
\| g^\eps b(\D) \v_\eps - \wt{g}^\eps S_\eps b(\D) \u_0 \|_{L_2(\R^d)}
\le C_6 \eps \|\D b(\D) \u_0 \|_{L_2(\R^d)},\quad \eps >0,
\eqno(1.47)
$$
where
$C_6 = (dm)^{1/2} (2r_0)^{-1} \alpha_1^{1/2} \alpha_0^{-1/2}
\|g\|^{3/2}_{L_\infty} \|g^{-1}\|^{1/2}_{L_\infty}
+  r_1 \|g\|_{L_\infty}$.

Similarly to (1.35),
$$
\| \D b(\D)\u_0\|^2_{L_2(\R^d)}
\le \alpha_1 \int_{\R^d} |\bxi|^4 (c_0 |\bxi|^2 + \lambda)^{-2} |\wh{\FF}(\bxi)|^2\,d\bxi
\le \alpha_1 c_0^{-2} \|\FF\|^2_{L_2(\R^d)}.
\eqno(1.48)
$$

Combining (1.42) and (1.47), (1.48), we obtain (1.41) with the constant
$C_5(\lambda) = C_5' \lambda^{-1/2} + C_5''$, where
$C_5' = \wh{C}_2 \alpha_1^{1/2} \|g\|_{L_\infty}$,
$C_5'' = \wh{C}_2 \alpha_1^{1/2} \|g\|_{L_\infty} + C_6 \alpha_1^{1/2} \alpha_0^{-1}\|g^{-1}\|_{L_\infty}$.
$\ \bullet$

\smallskip
Now we distinguish the case where the corrector is equal to zero.
Next statement follows from Theorem 1.7, Proposition 1.2, and equation (1.7).

\smallskip\noindent\textbf{Proposition 1.9.}
\textit{Let} $\u_\eps$ \textit{be the solution of the equation}
(1.20), \textit{and let} $\u_0$ \textit{be the solution of the equation} (1.21).
\textit{If $g^0 = \overline{g}$,
 i.~e., relations}  (1.13) \textit{are satisfied, then $\Lambda =0$, $K_\lambda(\eps)=0$, and we have}
 $$
 \| \u_\eps - \u_0\|_{H^1(\R^d;\C^n)} \le C_2(\lambda) \eps \|\FF\|_{L_2(\R^d)},\quad \eps >0.
 $$

\smallskip\noindent\textbf{1.7. The results for homogenization problem in $\R^d$
in the case where $\Lambda \in L_\infty$.}
It turns out that under some additional assumptions on the solution of the problem (1.7)
the smoothing operator $S_{\varepsilon}$ in (1.25) can be removed
(replaced by the identity operator). We impose the following condition.

\smallskip\noindent\textbf{Condition 1.10.}
\textit{Suppose that the} $\Gamma$-\textit{periodic solution}
$\Lambda(\x)$ \textit{of the problem} (1.7) \textit{is bounded}: $\Lambda \in L_\infty$.

We need the following
multiplicative property of $\Lambda$, see [PSu2, Corollary 2.4].

\smallskip\noindent\textbf{Proposition 1.11.}
\textit{Under Condition} 1.10 \textit{for any function $u \in H^1(\R^d)$ and
$\eps>0$ we have}
$$
\int_{\R^d} |(\D \Lambda)^\eps(\x)|^2 |u|^2\,d\x \le
\beta_1 \|u\|^2_{L_2(\R^d)} + \beta_2 \|\Lambda\|^2_{L_\infty} \eps^2 \int_{\R^d} |\D u|^2 d\x.
$$
\textit{The constants $\beta_1$, $\beta_2$ are given by}
$$
\begin{aligned}
\beta_1 &= 16 m \alpha_0^{-1} \|g\|_{L_\infty} \|g^{-1}\|_{L_\infty},
\cr
\beta_2 &= 2(1+ 2d \alpha_0^{-1} \alpha_1 + 20 d \alpha_0^{-1} \alpha_1 \|g\|_{L_\infty} \|g^{-1}\|_{L_\infty}).
\end{aligned}
$$

We put
$$
K^0_\lambda (\varepsilon) = [\Lambda^{\varepsilon}] b (\D) (\A^0 + \lambda I)^{-1}.
$$
By Proposition 1.11, it is easily seen that under Condition 1.10
the operator $K^0_\lambda (\varepsilon)$ is a continuous mapping of
$L_2 (\mathbb{R}^{d}; \mathbb{C}^{n})$ to $H^1 (\mathbb{R}^{d}; \mathbb{C}^{n})$.
Instead of (1.26), consider another first order approximation of $\u_{\varepsilon}$:
$$
\check{\v}_{\varepsilon} = \u_0 + \varepsilon \Lambda^{\varepsilon} b (\D) \u_0
= (\A^0 + \lambda I)^{-1} \FF + \varepsilon K^0_\lambda (\varepsilon) \FF.
\eqno(1.49)
$$

\smallskip\noindent\textbf{Theorem 1.12.} \textit{Suppose that Condition} 1.10
\textit{is satisfied. Let} $\u_\eps$ \textit{be the solution of the equation}
(1.20), \textit{and let} $\u_0$ \textit{be the solution of the equation} (1.21).
\textit{Let} $\check{\v}_\eps$ \textit{be defined by} (1.49). \textit{Then}
$$
 \| \u_{\varepsilon} - \check{\v}_\eps \|_{H^1 (\mathbb{R}^{d}; \mathbb{C}^{n})}
 \leq C_7(\lambda) \varepsilon  \| \FF \|_{L_2 (\mathbb{R}^{d}; \mathbb{C}^{n})},
 \quad  \varepsilon >0,
 \eqno(1.50)
 $$
\textit{or, in operator terms},
$$
\| (\A_{\varepsilon} + \lambda I)^{-1} - (\A^0 + \lambda I)^{-1} -
\varepsilon K^0_\lambda (\varepsilon) \|_{L_2 (\mathbb{R}^{d}; \mathbb{C}^{n})
\to H^1 (\mathbb{R}^{d}; \mathbb{C}^{n})} \leq \ C_7(\lambda) \varepsilon,\ \eps>0.
$$
\textit{For the flux $\p_\eps = g^\eps b(\D)\u_\eps$ we have}
$$
 \| \p_{\varepsilon} - \wt{g}^\eps b(\D) \u_0 \|_{L_2 (\mathbb{R}^{d}; \mathbb{C}^{m})}
 \leq C_8(\lambda) \varepsilon  \| \FF \|_{L_2 (\mathbb{R}^{d}; \mathbb{C}^{n})},
 \quad  \varepsilon >0,
 \eqno(1.51)
$$
\textit{where $\wt{g}(\x)$ is the matrix} (1.9),
$C_7(\lambda)=C_7'\lambda^{-1/2} + C_7''$,
$C_8(\lambda)=C_8'\lambda^{-1/2} + C_8''$,
\textit{and the constants $C_7'$, $C_7''$, $C_8'$, $C_8''$ depend only on}
$d, m, \alpha_0, \alpha_1, \| g \|_{L_\infty}, \| g^{-1} \|_{L_\infty}$,
\textit{the parameters of the lattice} $\Gamma$,
\textit{and the norm} $\|\Lambda\|_{L_\infty}$.

\smallskip\noindent\textbf{Proof.}
For the proof of (1.50) we need to estimate the $H^1$-norm of
the function $\eps \Lambda^\eps (I-S_\eps) b(\D)\u_0$.
We start with the $L_2$-norm. By Condition 1.10 and (1.18),
$$
\eps \|\Lambda^\eps (I-S_\eps) b(\D)\u_0\|_{L_2(\R^d)}
\le 2 \eps \|\Lambda\|_{L_\infty} \|b(\D)\u_0\|_{L_2(\R^d)}.
\eqno(1.52)
$$
Consider the derivatives
$$
\partial_l(\eps \Lambda^\eps (I-S_\eps) b(\D)\u_0) =
(\partial_l \Lambda)^\eps (I-S_\eps) b(\D)\u_0
+ \eps \Lambda^\eps (I-S_\eps) b(\D) \partial_l \u_0.
\eqno(1.53)
$$
The second term in the right-hand side of (1.53)
is estimated with the help of Condition 1.10 and (1.18):
$$
\sum_{l=1}^d \| \eps \Lambda^\eps (I-S_\eps) b(\D) \partial_l \u_0 \|^2_{L_2(\R^d)} \le
4 \eps^2 \|\Lambda\|^2_{L_\infty} \|\D b(\D)\u_0\|^2_{L_2(\R^d)}.
\eqno(1.54)
$$
Using Proposition 1.11,
we estimate the first term in the right-hand side of (1.53):
$$
\begin{aligned}
&\sum_{l=1}^d \|(\partial_l \Lambda)^\eps (I-S_\eps) b(\D)\u_0\|^2_{L_2(\R^d)}
\cr
&\le \beta_1 \| (I-S_\eps) b(\D)\u_0\|^2_{L_2(\R^d)} +
\beta_2 \eps^2 \|\Lambda\|_{L_\infty}^2 \sum_{l=1}^d \| (I-S_\eps) b(\D) D_l \u_0\|^2_{L_2(\R^d)}.
\end{aligned}
\eqno(1.55)
$$
By Proposition 1.4,
$$
\| (I-S_\eps) b(\D)\u_0\|_{L_2(\R^d)} \le \eps r_1 \|\D b(\D) \u_0 \|.
\eqno(1.56)
$$
Now from (1.55), (1.56), and (1.18) it follows that
$$
\sum_{l=1}^d \|(\partial_l \Lambda)^\eps (I-S_\eps) b(\D)\u_0\|^2_{L_2(\R^d)}
\le \eps^2 \left(\beta_1  r_1^2  + 4 \beta_2 \|\Lambda\|_{L_\infty}^2\right) \|\D b(\D) \u_0 \|^2_{L_2(\R^d)}.
\eqno(1.57)
$$
As a result, combining (1.53), (1.54), and (1.57), we arrive at
$$
\begin{aligned}
&\sum_{l=1}^d \| \partial_l(\eps \Lambda^\eps (I-S_\eps) b(\D)\u_0) \|^2_{L_2(\R^d)}
\cr
&\le 2\eps^2 \left( 4(1+\beta_2)\|\Lambda\|_{L_\infty}^2 + \beta_1 r_1^2
\right)\| \D b(\D) \u_0\|^2_{L_2(\R^d)}.
\end{aligned}
\eqno(1.58)
$$

Similarly to (1.32), we have
$$
\| b(\D) \u_0\|^2_{L_2(\R^d)}
\le \alpha_1 \int_{\R^d} |\bxi|^2 (c_0 |\bxi|^2 + \lambda)^{-2} |\wh{\FF}(\bxi)|^2\,d\bxi
\le \alpha_1 (2\lambda c_0)^{-1} \|\FF\|^2_{L_2(\R^d)}.
$$
Together with (1.52), (1.58), and (1.48) this yields
$$
\| \eps \Lambda^\eps (I-S_\eps) b(\D)\u_0 \|_{H^1(\R^d)}
\le C_9(\lambda) \eps \|\FF\|_{L_2(\R^d)},
\eqno(1.59)
$$
where $C_9(\lambda) = C_9' \lambda^{-1/2} + C_9''$,
$C_9'= \sqrt{2} \|\Lambda \|_{L_\infty} \alpha_1^{1/2} \alpha_0^{-1/2} \|g^{-1}\|_{L_\infty}^{1/2}$,
$$
C_9'' = \sqrt{2}\alpha_1^{1/2} \alpha_0^{-1} \|g^{-1}\|_{L_\infty}
\left( 2  (1+\beta_2)^{1/2}\|\Lambda\|_{L_\infty} +
\beta_1^{1/2} r_1 \right).
$$

As a result, the desired estimate (1.50) follows from
(1.26), (1.27), (1.49), and (1.59); the constant $C_7(\lambda)$ is given by
$C_7(\lambda) = C_7'\lambda^{-1/2} + C_7''$, where
$C_7'= \wh{C}_2 + C_9'$, $C_7'' = \wh{C}_2 + C_9''$.

To prove (1.51), note that from (1.50) and (1.4) it follows that
$$
\|\p_\eps - g^\eps b(\D) \check{\v}_\eps \|_{L_2(\R^d)} \le
\eps C_7(\lambda) \|g\|_{L_\infty} \alpha_1^{1/2}\|\FF\|_{L_2(\R^d)}.
\eqno(1.60)
$$
Next, by (1.2), (1.9), and (1.49)
$$
g^\eps b(\D) \check{\v}_\eps = \wt{g}^\eps b(\D)\u_0 +
\eps \sum_{l=1}^d g^\eps b_l \Lambda^\eps b(\D) D_l \u_0.
\eqno(1.61)
$$
The last term in the right-hand side of (1.61) is estimated with the help of
(1.5), (1.48), and Condition 1.10:
$$
\left\| \eps \sum_{l=1}^d g^\eps b_l \Lambda^\eps b(\D) D_l \u_0 \right\|_{L_2(\R^d)} \le
\eps \|g\|_{L_\infty} \|\Lambda\|_{L_\infty} \alpha_1 c_0^{-1} d^{1/2} \|\FF\|_{L_2(\R^d)}.
\eqno(1.62)
$$
Now relations (1.60)--(1.62) imply (1.51) with the constant
$C_8(\lambda)= C_8'\lambda^{-1/2} + C_8''$, where
$C_8'= C_7' \|g\|_{L_\infty} \alpha_1^{1/2}$,
$$
C_8'' = C_7'' \|g\|_{L_\infty} \alpha_1^{1/2}
+ d^{1/2}\|\Lambda\|_{L_\infty} \alpha_1 \alpha_0^{-1}\|g\|_{L_\infty} \|g^{-1}\|_{L_\infty}. \ \ \bullet
$$

\smallskip
In some cases Condition 1.10 is valid automatically.
Next statement was proved in [BSu4, Lemma 8.7].

\smallskip\noindent\textbf{Proposition 1.13.} \textit{Condition} 1.10
\textit{is a fortiori valid if at least one of the following assumptions is fulfilled}:

$1^\circ$. $d \leq 2$;

$2^\circ$. \textit{$d\ge 1$ and}
$\A_\eps = \D^* g^\eps(\x) \D$,
\textit{where} $g(\x)$ \textit{has real entries};

$3^\circ$. \textit{dimension is arbitrary, and} $g^0 = \underline{g}$, \textit{i.~e., relations}
(1.14) \textit{are satisfied}.

\smallskip
Note that Condition 1.10 is also ensured if $g(\x)$ is smooth enough.

We distinguish the special case where $g^0=\underline{g}$. In this case the matrix
(1.9) is constant: $\wt{g}(\x)=g^0 =\underline{g}$;
moreover, Condition 1.10 is satisfied.
Applying the statement of Theorem 1.12
concerning the fluxes, we arrive at the following statement.

\smallskip\noindent\textbf{Proposition 1.14.}
\textit{Suppose that} $\u_\eps$ \textit{is the solution of the equation}
(1.20), $\u_0$ \textit{is the solution of the equation} (1.21), \textit{and}
$\p_\eps = g^\eps b(\D)\u_\eps$.
\textit{Let $g^0 = \underline{g}$, i.~e., relations} (1.14) \textit{are satisfied. Then we have}
$$
 \| \p_{\varepsilon} - {g}^0 b(\D) \u_0 \|_{L_2 (\mathbb{R}^{d}; \mathbb{C}^{m})}
 \leq C_8(\lambda) \varepsilon  \| \FF \|_{L_2 (\mathbb{R}^{d}; \mathbb{C}^{n})},
 \quad  \varepsilon >0.
$$

\section*{\S 2. The Neumann problem in a bounded domain: preliminaries}

\noindent\textbf{2.1. Coercivity.}
Let $\mathcal{O} \subset \mathbb{R}^{d}$ be a bounded domain with the boundary of class $C^{1,1}$.
We impose an additional condition on the symbol
$b(\bxi)=\sum_{l=1}^d b_l \xi_l$ for $\bxi \in \C^d$.

\smallskip\noindent\textbf{Condition 2.1.} \textit{The matrix-valued function
$b(\bxi)$, $\bxi \in \C^d$, is such that}
$$
\rank b(\bxi) =n, \quad 0 \ne \bxi \in \C^d.
\eqno(2.1)
$$

\smallskip
Note that Condition 2.1 is more restrictive than (1.3).
According to [Ne] (see Theorem 7.8 in \S 3.7),
Condition 2.1 \textit{is necessary and sufficient for coercivity
of the form} $\|b(\D)\u\|^2_{L_2(\O)}$ \textit{on $H^1(\O;\C^n)$}
(moreover, this is true in any bounded domain $\O$
with the Lipschitz boundary).

\smallskip\noindent\textbf{Proposition 2.2. [Ne]}
\textit{Condition} 2.1 \textit{is necessary and sufficient for existence of constants
${\mathcal C}_1 >0$ and ${\mathcal C}_2 \ge 0$ such that the G\"arding type inequality}
$$
\|b(\D)\u \|_{L_2(\O)}^2 + {\mathcal C}_2 \|\u\|^2_{L_2(\O)} \ge {\mathcal C}_1 \|\D \u\|^2_{L_2(\O)},\quad
\u \in H^1(\O;\C^n),
\eqno(2.2)
$$
\textit{is satisfied.}

\smallskip\noindent\textbf{Remark 2.3.} The constants ${\mathcal C}_1$ and ${\mathcal C}_2$
depend on the matrix $b(\bxi)$ and the domain $\O$, but in general case it is difficult to control
these constants explicitly.
However, often for particular operators they are known.
Therefore, in what follows we indicate the dependence of other constants on
${\mathcal C}_1$ and ${\mathcal C}_2$.

\smallskip

In what follows, \textit{we assume that Condition} 2.1 \textit{is satisfied}.

\smallskip\noindent\textbf{2.2. Statement of the problem.}
In $L_2 (\mathcal{O}; \mathbb{C}^{n})$, consider the operator
$\mathcal{A}_{N,\varepsilon}$ formally given by the differential expression
$b (\mathbf{D})^* g^{\varepsilon} (\mathbf{x}) b (\mathbf{D})$
with the Neumann condition on $\partial \mathcal{O}$.
The precise definition of $\mathcal{A}_{N,\varepsilon}$ is given in terms of the
quadratic form
$$
a_{N,\eps}[\u,\u] := \int_{\mathcal{O}} \left\langle g^{\varepsilon} (\mathbf{x}) b (\mathbf{D}) \mathbf{u},
b(\mathbf{D}) \mathbf{u} \right\rangle \, d \mathbf{x},
\quad \mathbf{u} \in H^1 (\mathcal{O}; \mathbb{C}^{n}).
\eqno(2.3)
$$

By (1.2) and (1.5),
$$
a_{N,\eps}[\u,\u] \le d \alpha_1 \|g\|_{L_\infty} \|\D\u\|^2_{L_2(\O)},\quad \u \in H^1(\O;\C^n).
\eqno(2.4)
$$
From (2.2) it follows that
$$
a_{N,\eps}[\u,\u] \ge \|g^{-1}\|^{-1}_{L_\infty} \left(
{\mathcal C}_1 \|\D \u\|^2_{L_2(\O)} - {\mathcal C}_2 \| \u\|^2_{L_2(\O)}\right),
\quad \u \in H^1(\O;\C^n).
\eqno(2.5)
$$
By (2.4) and (2.5), the form (2.3) is closed and lower semibounded.
\textit{By definition, $\mathcal{A}_{N,\varepsilon}$ is the selfadjoint operator in
$L_2(\O;\C^n)$ generated by the quadratic form} (2.3).

\textit{Let $\lambda >0$ be a parameter subject to the following restriction}:
$$
\lambda > {\mathcal C}_2 \| g^{-1}\|^{-1}_{L_\infty}.
\eqno(2.6)
$$
Then (2.5) implies that
$$
\begin{aligned}
&a_{N,\eps}[\u,\u] + \lambda \|\u\|^2_{L_2(\O)} \ge c_\lambda \|\u\|^2_{H^1(\O)},\quad \u \in H^1(\O;\C^n),
\cr
&c_\lambda := \min \{ {\mathcal C}_1 \|g^{-1}\|_{L_\infty}^{-1},
\lambda - {\mathcal C}_2 \|g^{-1}\|_{L_\infty}^{-1}\}.
\end{aligned}
\eqno(2.7)
$$
Hence, the operator
$\mathcal{A}_{N,\varepsilon}+ \lambda I$ is positive definite.

\textit{Our goal} is to approximate the inverse operator
$(\mathcal{A}_{N,\varepsilon}+ \lambda I)^{-1}$ for small $\eps$
in the operator norm in $L_2(\O;\C^n)$
and in the norm of operators acting from $L_2(\O;\C^n)$ to $H^1(\O;\C^n)$.
In terms of solutions, we are interested in the behavior of the generalized solution
$\u_\eps \in H^1(\O;\C^n)$ of the Neumann problem
$$
b (\mathbf{D})^* g^{\varepsilon} (\mathbf{x}) b (\mathbf{D}) \mathbf{u}_{\varepsilon}(\x)
+ \lambda \u_\eps(\x) = \mathbf{F}(\x),
\ \ \x \in \O;
\quad \partial^\eps_{\bnu} \u_\eps\vert_{\partial \O}=0,
\eqno(2.8)
$$
where $\mathbf{F} \in L_2 (\mathcal{O}; \mathbb{C}^{n})$.
Then $\u_\eps = (\mathcal{A}_{N,\varepsilon}+\lambda I)^{-1} \FF$.

Here the notation
$\partial^\eps_{\bnu}$ for the "conormal derivative" was used.
Let $\bnu(\x)=\sum_{j=1}^d \nu_j(\x) \e_j$ be the unit outer normal vector to
$\partial \O$ at the point $\x \in \partial \O$.
Here $\e_1,\dots,\e_d$ is the standard orthonormal basis in $\R^d$.
Introduce the matrix
$b^*(\bnu(\x))= \sum_{l=1}^d b_l^* \nu_l(\x)$. Then formally
the conormal derivative is given by
$$
\partial^\eps_{\bnu} \u(\x) := b^*(\bnu(\x)) g^\eps(\x) b(\partial) \u(\x)=
\sum_{l,j=1}^d  b_l^* g^\eps(\x)  b_j \nu_l(\x) \partial_j \u(\x),\ \ \x \in \partial \O.
$$

Note that, under our assumptions, (2.8) has formal sense.
\textit{By definition, we say that $\u_\eps \in H^1(\O;\C^n)$ is the generalized
solution of the problem} (2.8) \textit{if $\u_\eps$ satisfies the following integral identity}
$$
\int_\O \left( \langle g^\eps b(\D)\u_\eps, b(\D) \eeta \rangle +
\lambda \langle \u_\eps, \eeta \rangle \right)\,d\x = \int_\O \langle \FF,\eeta \rangle\,d\x,
\quad \forall \eeta \in H^1(\O;\C^n).
\eqno(2.9)
$$
Taking (2.4) and (2.7) into account, we can view the form $a_{N,\eps}[\u,\eeta]+\lambda (\u,\eeta)_{L_2(\O)}$
as a (new) inner product in $H^1(\O;\C^n)$.
The right-hand side of (2.9) is an antilinear continuous functional of $\eeta \in H^1(\O;\C^n)$.
By the Riss theorem, there exists unique solution $\u_\eps$;
by (2.7), it satisfies the estimate
$$
\| \u_\eps \|_{H^1(\O)} \le c_\lambda^{-1} \|\FF \|_{L_2(\O)}.
\eqno(2.10)
$$
In operator terms, (2.10) means that
$$
\|(\A_{N,\eps} + \lambda I)^{-1}\|_{L_2(\O)\to H^1(\O)} \le c_\lambda^{-1}.
\eqno(2.11)
$$

\smallskip\noindent\textbf{2.3. The trace and extension operators.}
It is well known that, under our assumption on the domain $\O$
(that the boundary is of class $C^{1,1}$),
the trace operator $\gamma$ taking a function $\u$ in
$\O$ into its trace on the boundary $\partial \O$ is defined correctly as a linear continuous
operator
$$
\gamma: H^s(\O;\C^n) \to H^{s-1/2}(\O;\C^n), \quad s=1, 2.
\eqno(2.12)
$$
Herewith,
$$
\begin{aligned}
&\|\gamma\|_{H^1(\O)\to H^{1/2}(\partial \O)} \le \wh{c}_1,
\cr
&\|\gamma\|_{H^2(\O)\to H^{3/2}(\partial \O)} \le \wh{c}_2,
\end{aligned}
\eqno(2.13)
$$
where the constants $\wh{c}_1$ and $\wh{c}_2$ depend only on the domain $\O$.

There exists a linear continuous extension operator
taking a function at the boundary into
its extension to the domain $\O$ (this is the right inverse to the operator (2.12)).
Obviously, such an operator is not unique. It is convenient to choose
this operator in the following way.
Consider the Dirichlet problem
$$
-\Delta \u + \u = 0 \ \ \text{in}\ \O, \quad \u\vert_{\partial \O} = \bvarphi,
\eqno(2.14)
$$
where $\bvarphi \in H^{1/2}(\partial \O;\C^n)$.
It is well known that there exists a unique solution $\u \in H^1(\O;\C^n)$, and it
satisfies the estimate
$
\|\u \|_{H^1(\O)} \le \wt{c}_1 \|\bvarphi\|_{H^{1/2}(\partial \O)}.
$
The constant $\wt{c}_1$ depends only on the domain $\O$.
\textit{Denote by $T$ the operator taking a function $\bvarphi$ into the solution
 $\u$ of the problem} (2.14): $\u = T \bvarphi$. Then
 $$
 T: H^{1/2}(\partial \O;\C^n) \to H^1(\O;\C^n)
 \eqno(2.15)
 $$
  is a linear continuous operator, $T$ is the right inverse of
  the operator (2.12) (with $s=1$), and
  $$
  \|T\|_{H^{1/2}(\partial \O) \to H^1(\O)} \le \wt{c}_1.
  \eqno(2.16)
  $$

By the theorem about regularity of solutions (see, e.~g., [McL, Chapter 4]),
if $\bvarphi \in H^{3/2}(\partial \O;\C^n)$, then $\u \in H^2(\O;\C^n)$
and the operator
$$
T: H^{3/2}(\partial \O;\C^n) \to H^2(\O;\C^n)
\eqno(2.17)
$$
is continuous. Its norm satisfies the estimate
$$
 \|T\|_{H^{3/2}(\partial \O) \to H^2(\O)} \le \wt{c}_2,
 \eqno(2.18)
$$
where the constant $\wt{c}_2$ depends only on the domain $\O$.
The operator (2.17) is the right inverse of the operator (2.12) with $s=2$.

\smallskip\noindent\textbf{2.4. Definition of the conormal derivative.
The Neumann problem for the homogeneous equation.}
We need to give the precise definition of the conormal derivative
(corresponding to the operator $\A_\eps$).
The conormal derivative $\partial_{\bnu}^\eps \f$ is defined as an element of
$H^{-1/2}(\partial \O;\C^n)$,
definition is based on the Green formula, in which
the function $\f \in H^1(\O;\C^n)$ and the right-hand side of the equation are given.
For our purposes, it suffices to give the definition of the conormal derivative for a function
$\f \in H^1(\O;\C^n)$ in the case where $\A_\eps \f$ belongs to $L_2(\O;\C^n)$.

Thus, let $\f \in H^1(\O;\C^n)$ and $\bPhi_\eps \in L_2(\O;\C^n)$ be such that
$\A_\eps \f = \bPhi_\eps$ inside $\O$, which is understood as the integral identity
$$
\int_{\O}  \langle g^\eps b(\D)\f, b(\D)\eeta \rangle \,d\x
= \int_{\O} \langle \bPhi_\eps, \eeta \rangle\,d\x,
\quad \forall \eeta \in H^1_0(\O;\C^n).
\eqno(2.19)
$$

Let $l_\eps[\bvarphi]$ be the following antilinear functional of $\bvarphi \in H^{1/2}(\partial \O;\C^n)$:
$$
l_\eps[\bvarphi] =
\int_{\O}  \langle g^\eps b(\D)\f, b(\D)\u \rangle \,d\x - \int_{\O} \langle \bPhi_\eps, \u \rangle\,d\x,
\eqno(2.20)
$$
where $\u \in H^1(\O;\C^n)$ is a function such that $\gamma \u = \bvarphi$.
By (2.19), the right-hand side of (2.20) does not depend on the choice of an extension
$\u$ of the function $\bvarphi$; therefore, definition of
the functional (2.20) is correct.
It is convenient to put $\u = T \bvarphi$, where the operator $T$ is defined in Subsection 2.3.
Then, by (1.2), (1.5), and (2.16), the functional (2.20) is continuous:
$$
\begin{aligned}
|l_\eps[\bvarphi]| &\le \|\bPhi_\eps \|_{L_2(\O)} \|\u\|_{L_2(\O)} +
\alpha_1 d \|g\|_{L_\infty}  \|\D\f\|_{L_2(\O)} \|\D\u\|_{L_2(\O)}
\cr
&\le \wt{c}_1 \left( \|\bPhi_\eps \|_{L_2(\O)} + \alpha_1 d \|g\|_{L_\infty} \|\f\|_{H^1(\O)}\right)\|\bvarphi\|_{H^{1/2}(\partial \O)}.
\end{aligned}
\eqno(2.21)
$$

Recall that the space $H^{-1/2}(\partial\O;\C^n)$ is dual to $H^{1/2}(\partial\O;\C^n)$
with respect to the $L_2(\partial\O;\C^n)$-duality.
It is common to write the meaning of the functional $\bpsi \in H^{-1/2}(\partial\O;\C^n)$
on the element $\bvarphi \in H^{1/2}(\partial\O;\C^n)$ as
$(\bpsi,\bvarphi)_{L_2(\partial\O)}$ (which is extension of the inner product in
$L_2(\partial\O;\C^n)$ for the pairs
in $H^{-1/2}(\partial\O;\C^n) \times H^{1/2}(\partial\O;\C^n)$).
Herewith,
$$
\|\bpsi\|_{H^{-1/2}(\partial \O)}
= \sup_{0\ne \bvarphi \in H^{1/2}(\partial \O)}
\frac{|(\bpsi,\bvarphi)_{L_2(\partial \O)}|}{\|\bvarphi\|_{H^{1/2}(\partial \O)}}.
\eqno(2.22)
$$

For the antilinear continuous functional $l_\eps[\bvarphi]$ on $H^{1/2}(\partial\O;\C^n)$
defined by (2.20)
there exists a unique element $\bpsi_\eps \in H^{-1/2}(\partial\O;\C^n)$ such that
$$
l_\eps[\bvarphi] = (\bpsi_\eps, \bvarphi)_{L_2(\partial\O)},\quad \bvarphi \in H^{1/2}(\partial\O;\C^n).
\eqno(2.23)
$$
\textit{By definition, we say that $\bpsi_\eps$ is the conormal derivative of a function
$\f$ and write $\bpsi_\eps = \partial^\eps_{\bnu} \f$.}
From (2.21)--(2.23) it follows that
$$
\|\partial^\eps_{\bnu} \f\|_{H^{-1/2}(\partial \O)} \le
\wt{c}_1 \left( \|\bPhi_\eps\|_{L_2(\O)} + \alpha_1 d \|g\|_{L_\infty} \|\f\|_{H^1(\O)}\right).
$$

Now we discuss the statement of
the Neumann problem for the homogeneous equation with the nonhomogeneous boundary condition.
Let $\bpsi \in H^{-1/2}(\partial\O;\C^n)$.
Consider the generalized solution $\r_\eps$ of the problem
$$
\A_\eps \r_\eps + \lambda \r_\eps =0 \ \ \text{in}\ \O;
\ \ \partial^\eps_{\bnu} \r_\eps = \bpsi \ \ \text{on}\ \partial \O.
\eqno(2.24)
$$
By definition, \textit{the generalized solution of the problem} (2.24)
\textit{is an element $\r_\eps \in H^1(\O;\C^n)$ satisfying the identity}
$$
\int_{\O} \left( \langle g^\eps b(\D)\r_\eps, b(\D)\eeta \rangle +
\lambda \langle \r_\eps, \eeta \rangle\right)\,d\x =
 (\bpsi, \gamma \eeta)_{L_2(\partial \O)},
\quad \forall \eeta \in H^1(\O;\C^n).
\eqno(2.25)
$$

The following statement is checked in a standard way.

\smallskip\noindent\textbf{Proposition 2.4.} \textit{Let} $\bpsi \in H^{-1/2}(\partial\O;\C^n)$.
\textit{Then the generalized solution} $\r_\eps \in H^1(\O;\C^n)$
\textit{of the problem} (2.24) \textit{exists, it is unique and satisfies the estimate}
$$
\| \mathbf{r}_{\varepsilon} \|_{H^1 (\mathcal{O}; \mathbb{C}^{n})} \leq
c_\lambda^{-1} \widehat{c}_1 \| \bpsi \|_{H^{-1/2} (\partial \O; \mathbb{C}^{n})},
\eqno(2.26)
$$
\textit{where the constant $c_\lambda$ is defined by} (2.7), \textit{and $\wh{c}_1$ by} (2.13).

\smallskip\noindent\textbf{Proof.} As has been already mentioned, the form in the left-hand side of (2.25)
can be viewed as an inner product in $H^1(\O;\C^n)$.
By (2.13) and (2.22), the right-hand side of (2.25) is an antilinear continuous functional of
$\eeta \in H^1(\O;\C^n)$, and we have
$$
|(\bpsi, \gamma \eeta)_{L_2(\O)}| \le \wh{c}_1 \|\bpsi\|_{H^{-1/2}(\partial \O)}
\|\eeta\|_{H^1(\O)},\quad \eeta \in H^1(\O;\C^n).
$$
Then, by the Riss theorem, the solution $\r_\eps$ exists and is unique.
By (2.7), it satisfies the estimate (2.26).
$\ \bullet$

\smallskip\noindent\textbf{2.5. The "homogenized" problem.}
In $L_2(\O;\C^n)$, consider the selfadjoint operator $\A_N^0$ generated by the quadratic form
$$
a^0_N[\u,\u]= \int_{\mathcal{O}} \left\langle g^0 b (\mathbf{D}) \mathbf{u}, b(\mathbf{D}) \u \right\rangle \,d \mathbf{x},
 \quad \u \in H^1 (\mathcal{O}; \mathbb{C}^{n}).
\eqno(2.27)
$$
Here $g^0$ is the effective matrix defined by (1.8).
By (1.15), estimates (2.4), (2.5) remain true for the form (2.27).

Suppose that $\lambda$ still satisfies (2.6).
Let $\u_0 \in H^1(\O;\C^n)$ be the generalized solution of the Neumann problem
$$
b(\D)^* g^0 b(\D) \u_0(\x) + \lambda \u_0(\x)= \FF(\x),\ \ \x \in \O;
\ \ \partial_{\bnu}^0 \u_0\vert_{\partial \O} =0,
\eqno(2.28)
$$
where $\FF \in L_2(\O;\C^n)$. Then $\u_0 = (\mathcal{A}^0_N + \lambda I)^{-1} \FF$.
Here the formal differential expression for the conormal derivative
(corresponding to $\A^0$) is of the form
$$
\partial_{\bnu}^0 \u (\x) = b^*(\bnu(\x)) g^0 b(\partial) \u(\x)=
\sum_{l,j=1}^d b_l^* g^0 b_j \nu_l(\x) \partial_j \u(\x),\ \ \x \in \partial \O.
\eqno(2.29)
$$

\textit{We say that $\u_0\in H^1(\O;\C^n)$ is the generalized solution of the problem} (2.28)
\textit{if it satisfies the integral identity}
$$
\int_\O \left( \langle g^0 b(\D)\u_0, b(\D) \eeta \rangle +
\lambda \langle \u_0, \eeta \rangle \right)\,d\x = \int_\O \langle \FF,\eeta \rangle\,d\x,
\quad \forall \eeta \in H^1(\O;\C^n).
\eqno(2.30)
$$
The solution exists, it is unique and satisfies the estimate
$$
\| \u_0 \|_{H^1(\O)} \le c_\lambda^{-1} \|\FF \|_{L_2(\O)}.
\eqno(2.31)
$$
In operator terms, (2.31) means that
$$
\|(\A_{N}^0 + \lambda I)^{-1}\|_{L_2(\O)\to H^1(\O)} \le c_\lambda^{-1}.
\eqno(2.32)
$$

By the condition $\partial \mathcal{O} \in C^{1,1}$,
the solution $\mathbf{u}_0$ of the problem (2.28) satisfies
$
\mathbf{u}_0 \in H^2 (\mathcal{O}; \mathbb{C}^{n})
$
and
$$
\| \mathbf{u}_0 \|_{H^2 (\mathcal{O}; \mathbb{C}^{n})} \leq \wh{C}_\lambda
\| \mathbf{F} \|_{L_2 (\mathcal{O}; \mathbb{C}^{n})}.
\eqno(2.33)
$$
Here the constant $\wh{C}_\lambda$ depends only on the constants
${\mathcal C}_1$ and ${\mathcal C}_2$ from (2.2),
on $\alpha_0$, $\alpha_1$, $\|g\|_{L_\infty}$, $\|g^{-1}\|_{L_\infty}$, $\lambda$, and the domain
$\O$. To justify this fact, we note that the operator $b(\D)^* g^0 b(\D)$
is a \textit{strongly elliptic}
matrix operator with constant coefficients and refer to the theorems about regularity of solutions
for strongly elliptic systems (see, e.~g., [McL, Chapter 4]).

From what was said it follows that
the left-hand side of (2.28) belongs to $L_2(\O;\C^n)$;
equation is valid almost everywhere in $\O$.
The conormal derivative $\partial^0_{\bnu} \u_0$
is defined correctly by (2.29) as an element of $H^{1/2}(\partial \O;\C^n)$;
the boundary condition in (2.28) can be understood in the sense of the trace theorem.
Estimate (2.33) can be written as
$$
\| (\mathcal{A}^0_N + \lambda I)^{-1}\|_{L_2(\O;\C^n) \to H^2(\O;\C^n)} \le \wh{C}_\lambda.
\eqno(2.34)
$$

In what follows, it is shown that the solution $\u_\eps$ of the problem (2.8) converges in $L_2(\O;\C^n)$
to the solution $\u_0$ of the homogenized problem (2.28), as $\eps \to 0$.
We will estimate the error term $\|\u_\eps - \u_0\|_{L_2(\O)}$. Also, we will find approximation
of $\u_\eps$ in $H^1(\O;\C^n)$.

\section*{\S 3. Auxiliary statements}

This section contains various auxiliary statements needed below.

\smallskip\noindent\textbf{3.1.} Next statement is a version of the traditional lemma
used in homogenization theory (see, e.~g., [ZhKO, Chapter 1, \S 1]).

\smallskip\noindent\textbf{Lemma 3.1.} \textit{Let $f_l(\x)$, $l=1,\dots,d$, be $\Gamma$-periodic
$(n \times m)$-matrix-valued functions in $\R^d$ such that}
$$
f_l \in L_2(\Omega), \ \ \int_{\Omega} f_l(\x)\,d\x =0, \ \ l=1,\dots,d;
\quad \sum_{l=1}^d D_l f_l(\x) =0,
\eqno(3.1)
$$
\textit{where the last equation is understood in the distribution sense.
There exist $\Gamma$-periodic $(n \times m)$-matrix-valued functions $M_{lj}(\x)$ in $\R^d$,
$l,j=1,\dots,d$, such that}
$$
M_{lj} \in \wt{H}^1(\Omega),\ \ \int_\Omega M_{lj}(\x)\,d\x=0,
\ \ M_{lj}(\x) = - M_{jl}(\x),\ \ l,j=1,\dots,d,
\eqno(3.2)
$$
$$
f_l(\x) = \sum_{j=1}^d \partial_j M_{lj}(\x),\ \ l=1,\dots,d.
\eqno(3.3)
$$
\textit{We have}
$$
\|M_{lj}\|_{L_2(\Omega)} \le (2r_0)^{-1} \left(\|f_l\|_{L_2(\Omega)} +  \|f_j\|_{L_2(\Omega)}\right),
\quad l,j=1,\dots,d.
\eqno(3.4)
$$

\smallskip\noindent\textbf{Proof.}
Let $\Phi_l(\x)$, $l=1,\dots,d,$ be $\Gamma$-periodic
$(n \times m)$-matrix-valued functions in $\R^d$ such that
$$
\Delta \Phi_l(\x) = f_l(\x),\ \ \int_\Omega \Phi_l(\x)\,d\x =0,\ \ l=1,\dots,d.
\eqno(3.5)
$$
The solvability condition $\overline{f_l}=0$ for the equation
in (3.5) is satisfied.
The solution $\Phi_l$ exists, it is unique, and
$\Phi_l \in \wt{H}^2(\Omega)$, since $f_l \in L_2(\Omega)$.
We put
$$
M_{lj}(\x):= \partial_j \Phi_l(\x) - \partial_l \Phi_j(\x),\quad l,j =1,\dots,d.
\eqno(3.6)
$$
Then $M_{lj}\in \wt{H}^1(\Omega)$, $\overline{M_{lj}}=0$, and
$M_{lj}(\x)=-M_{jl}(\x)$. We have
$$
\sum_{j=1}^d \partial_j M_{lj}(\x) = \sum_{j=1}^d (\partial_j^2 \Phi_l(\x) - \partial_j\partial_l \Phi_j(\x))
= \Delta \Phi_l(\x) - \partial_l \left( \sum_{j=1}^d \partial_j \Phi_j(\x) \right).
\eqno(3.7)
$$
Note that, by (3.1) and (3.5),
$$
0 = \sum_{j=1}^d \partial_j f_j(\x) =
\sum_{j=1}^d \partial_j \Delta \Phi_j(\x)=
\Delta \left( \sum_{j=1}^d \partial_j \Phi_j(\x) \right).
$$
Thus, $\sum_{j=1}^d \partial_j \Phi_j$ is a periodic harmonic function
with zero mean value.
Hence, $\sum_{j=1}^d \partial_j \Phi_j=0$. Now from (3.5) and (3.7) it follows that
$f_l(\x) =\sum_{j=1}^d \partial_j M_{lj}(\x)$.
Relations (3.2) and (3.3) are proved.

It remains to check estimate (3.4).
Using (3.5), it is easily seen that $\Phi_l$ satisfies the estimate
$$
\|\D \Phi_l\|_{L_2(\Omega)}^2  \le \|f_l\|_{L_2(\Omega)} \|\Phi_l\|_{L_2(\Omega)}.
$$
Note that for any $\Gamma$-periodic function $\varphi \in \wt{H}^1(\Omega)$
with zero mean value we have $\|\varphi\|_{L_2(\Omega)}\le (2r_0)^{-1}\|\D \varphi\|_{L_2(\Omega)}$.
Hence,
$$
\|\D \Phi_l\|_{L_2(\Omega)} \le (2r_0)^{-1} \|f_l\|_{L_2(\Omega)}, \quad l=1,\dots,d.
$$
Together with (3.6) this implies (3.4). $\ \bullet$

\smallskip\noindent\textbf{3.2. Estimates for integrals over
the neighborhood of the boundary.}
The following statement is rather standard (see [PSu2, Lemma 5.1]).

\smallskip\noindent\textbf{Lemma 3.2.}
\textit{Let} $\mathcal{O} \subset \mathbb{R}^d$ \textit{be a bounded domain
with the boundary of class} $C^1$.
\textit{Denote}
$B_{\varepsilon} = \left\{ \mathbf{x} \in \mathcal{O}: \textrm{\upshape dist\itshape}
\,\{ \mathbf{x}, \partial \mathcal{O} \} < \varepsilon \right\}$.
\textit{Suppose that the number} $\eps_0 \in (0,1]$ \textit{is such that} $B_{\eps_0}$
\textit{can be covered by a finite number of open sets
admitting diffeomorphisms of class} $C^1$
\textit{rectifying the boundary} $\partial \O$.
\textit{Then for any function}
$u \in H^1 (\mathcal{O})$ \textit{we have}
$$
\int_{B_{\varepsilon}} |u|^{2} d \mathbf{x}
\leq \beta \varepsilon \| u \|_{H^1 (\mathcal{O})} \| u \|_{L_2 (\mathcal{O})},
\quad 0< \eps \leq \eps_0.
$$
\textit{The constant} $\beta=\beta(\O)$ \textit{depends only on the domain} $\mathcal{O}$.

\smallskip

The following statement is an analog of Lemma 2.6 from [ZhPas];
in the present form it was proved in [PSu2, Lemma 5.3].

\smallskip\noindent\textbf{Lemma 3.3.}
\textit{Let} $\mathcal{O} \subset \mathbb{R}^d$ \textit{be a bounded domain
with the boundary of class} $C^1$.
\textit{Denote} $\left(\partial \mathcal{O} \right)_{\varepsilon} =
\left\{ \mathbf{x} \in \mathbb{R}^d: \textrm{dist} \,
\{ \mathbf{x}, \partial \mathcal{O} \} < \varepsilon \right\}$.
\textit{Suppose that the number} $\eps_0 \in (0,1]$ \textit{is such that} $(\partial \O)_{\eps_0}$
\textit{can be covered by a finite number of open sets admitting
diffeomorphisms of class} $C^1$ \textit{rectifying the boundary} $\partial \O$.
\textit{Let} $S_\eps$ \textit{be the operator} (1.17).
\textit{Let} $f(\x)$ \textit{be a $\Gamma$-periodic function in} $\R^d$
\textit{such that} $f \in L_2(\Omega)$.
\textit{Then for any function} $\u \in H^1(\R^d;\C^m)$ \textit{we have}
$$
\begin{aligned}
\int_{\left( \partial \mathcal{O} \right)_{\varepsilon}} |f^\eps|^2
|S_\eps \u|^{2} \,d \mathbf{x}
\leq \beta_* \varepsilon |\Omega|^{-1}\|f\|^2_{L_2(\Omega)} \| \u \|_{H^1 (\R^d;\C^m)} \| \u \|_{L_2 (\R^d;\C^m)},
\\
0< \eps \leq \eps_1,
\end{aligned}
$$
\textit{where} $\eps_1 = \eps_0 (1+ r_1)^{-1}$,
$\beta_*= \beta^0 (1+ r_1)$, \textit{the constant} $\beta^0$
\textit{depends only on the domain $\mathcal{O}$, and} $2r_1 =\text{diam}\,\Omega$.

\smallskip
In what follows, the number $\eps_0$ is chosen as in Lemma 3.3. Then $\eps_0$
also satisfies the assumptions of Lemma 3.2.

\smallskip\noindent\textbf{3.3. The Dirichlet problem.}
We need the properties of solutions of the Dirichlet problem
for the equation $-\Delta \s + \s = \f$ in the domain $\O$.
Recall that the space $H^{-1}(\O;\C^n)$ is defined as the space dual to
$H^1_0(\O;\C^n)$ with respect to the $L_2(\O;\C^n)$-duality.
If $\f \in H^{-1}(\O;\C^n)$ and $\eeta \in H^1_0(\O;\C^n)$, the symbol $(\f,\eeta)_{L_2(\O)}$
is understood as the meaning of the functional $\f$ on the element $\eeta$.
Herewith,
$$
|(\f,\eeta)_{L_2(\O)}| \le \|\f\|_{H^{-1}(\O)} \|\eeta\|_{H^1(\O)}.
\eqno(3.8)
$$
The following lemma is checked in a standard way.

 \smallskip\noindent\textbf{Lemma 3.4.} \textit{Let
 $\f \in H^{-1}(\O;\C^n)$, and let $\s \in H^1_0(\O;\C^n)$ be the generalized solution
 of the Dirichlet problem
$$
-\Delta \s + \s = \f\ \ \text{in}\ \O,\quad \s\vert_{\partial \O}=0.
\eqno(3.9)
$$
 Then the following "energy estimate" is valid}:
 $$
 \|\s\|_{H^1(\O)} \le \|\f\|_{H^{-1}(\O)}.
 \eqno(3.10)
 $$

 \smallskip\noindent\textbf{Proof.}
 By definition of the generalized solution of the problem (3.9), we have
 $$
 \int_\O \left( \sum_{l=1}^d \langle D_l \s, D_l \eeta \rangle + \langle \s,\eeta \rangle \right)\,d\x =
 \int_{\O} \langle \f, \eeta \rangle\,d\x,\quad \forall \eeta \in H^1_0(\O;\C^n).
 \eqno(3.11)
 $$
 Substituting $\eeta =\s$ in (3.11) and using
 (3.8), we obtain (3.10). $\ \bullet$

 \smallskip
 Next statement concerns the Dirichlet problem for the homogeneous equation;
 it is deduced from Lemma 3.4.

 \smallskip\noindent\textbf{Lemma 3.5.}
 \textit{Let $\bphi \in H^1(\O;\C^n)$, and let $\h\in H^1(\O;\C^n)$ be the solution of the problem}
 $$
 -\Delta \h +\h =0\ \ \text{in}\ \O;\quad
 \h\vert_{\partial \O} = \bphi\vert_{\partial \O}.
 \eqno(3.12)
 $$
\textit{Then we have}
$$
\|\h\|_{H^1(\O)} \le 2 \|\bphi\|_{H^1(\O)}.
\eqno(3.13)
$$

\smallskip\noindent\textbf{Proof.}
 By (3.12), the function $\h - \bphi \in H^1_0(\O;\C^n)$ is the solution of the problem
 $$
 (-\Delta +I)(\h -\bphi) = \f \ \ \text{in}\ \O;\quad
 (\h-\bphi)\vert_{\partial \O} = 0,
 $$
where $\f = -(-\Delta +I)\bphi$.
Since
$(\f,\eeta)_{L_2(\O)} = -(\bphi,\eeta)_{H^1(\O)},$
$\eeta \in H^1_0(\O;\C^n)$,
then $\f \in H^{-1}(\O;\C^n)$, and we have
$$
\|\f\|_{H^{-1}(\O)} = \sup_{0\ne \eeta \in H^1_0(\O)}
\frac{|(\bphi,\eeta)_{H^1(\O)}|}{\|\eeta\|_{H^1(\O)}} = \|\bphi\|_{H^1(\O)}.
\eqno(3.14)
$$
By Lemma 3.4 and (3.14),
$$
\|\h - \bphi\|_{H^1(\O)}
\le \|\f\|_{H^{-1}(\O)} = \|\bphi\|_{H^1(\O)}.
$$
This implies (3.13). $\ \bullet$

\section*{\S 4. Homogenization of the Neumann problem: main results}

\smallskip\noindent\textbf{4.1. Formulations of the main results.}
Now we formulate main results of the paper.
We start with approximation of the operator $(\mathcal{A}_{N,\varepsilon} + \lambda I)^{-1}$
in the $L_2$-operator norm for small $\eps$.

\smallskip\noindent\textbf{Theorem 4.1.} \textit{Suppose that} $\mathcal{O}\subset \R^d$
\textit{is a bounded domain with the boundary of class} $C^{1,1}$, \textit{and the matrix} $g(\x)$
\textit{and DO} $b(\D)$ \textit{satisfy the assumptions of Subsection} 1.2.
\textit{Suppose that Condition} 2.1
\textit{is satisfied. Let $\lambda$ be subject to the restriction} (2.6).
 \textit{Let} $\mathbf{u}_{\varepsilon}$ \textit{be the solution of the problem} (2.8),
 \textit{and let} $\mathbf{u}_0$ \textit{be the solution of the problem} (2.28)
 \textit{with} $\FF \in L_2(\O;\C^n)$.
\textit{Then there exists a number} $\eps_1 \in (0,1]$
\textit{depending on the domain $\O$ and the lattice $\Gamma$ such that}
$$
\| \mathbf{u}_{\varepsilon} - \mathbf{u}_0
\|_{L_2 (\mathcal{O}; \mathbb{C}^{n})} \leq {\mathcal C}_0(\lambda) \varepsilon
\| \mathbf{F} \|_{L_2 (\mathcal{O}; \mathbb{C}^{n})},\ \ 0< \eps \leq \eps_1,
\eqno(4.1)
$$
\textit{or, in operator terms,}
$$
\| (\mathcal{A}_{N,\varepsilon} + \lambda I)^{-1} - (\mathcal{A}^0_N +\lambda I)^{-1}
 \|_{L_2 (\mathcal{O}; \C^n) \to L_2 (\mathcal{O}; \mathbb{C}^{n})}
 \leq {\mathcal C}_0(\lambda) {\varepsilon},\ \ 0< \eps \leq \eps_1.
 \eqno(4.2)
$$
\textit{The constant ${\mathcal C}_0(\lambda)$ depends on}
$m$, $d$, $\alpha_0$, $\alpha_1$, $\| g \|_{L_\infty}$,
$\| g^{-1} \|_{L_\infty}$, $\lambda$, \textit{the parameters of the lattice $\Gamma$,
the constants ${\mathcal C}_1$ and ${\mathcal C}_2$ from the inequality} (2.2),
\textit{and the domain} $\mathcal{O}$.

\smallskip

Our second result concerns approximation of the resolvent of
$\A_{N,\eps}$ in the norm of operators acting from $L_2(\O;\C^n)$ to $H^1(\O;\C^n)$.
In order to formulate this result, it is necessary to introduce the corrector.

We fix a linear continuous extension operator
$$
P_\O: H^2(\O;\C^n) \to H^2(\R^d;\C^n),
\eqno(4.3)
$$
and put $\wt{\u}_0 = P_\O \u_0$. Then
$$
\|\wt{\u}_0 \|_{H^2(\R^d;\C^n)} \le C_\O \| \u_0 \|_{H^2(\O;\C^n)},
\eqno(4.4)
$$
where $C_\O$ is the norm of the operator (4.3).

Let $S_\eps$ be the smoothing Steklov operator defined by
(1.17). By $R_\O$ we denote the operator of restriction of functions
in $\R^d$ onto the domain $\O$. The corrector for the Neumann problem is given by
$$
K_{N,\lambda}(\eps) = R_\O [\Lambda^\eps] S_\eps b(\D) P_\O (\A_N^0 + \lambda I)^{-1}.
\eqno(4.5)
$$
The operator $b(\D) P_\O (\A_N^0+ \lambda I)^{-1}$ is a continuous mapping of $L_2(\O;\C^n)$
to $H^1(\R^d;\C^m)$.
As was mentioned in Subsection 1.6, the operator $[\Lambda^\eps] S_\eps$
is continuous from $H^1(\R^d;\C^m)$ to $H^1(\R^d;\C^n)$.
Hence, the operator (4.5) is continuous from $L_2(\O;\C^n)$ to $H^1(\O;\C^n)$.

The first order approximation to the solution $\u_\eps$ is given by
$$
\v_\eps = (\A_N^0 + \lambda I)^{-1}\FF + \eps K_{N,\lambda}(\eps)\FF.
\eqno(4.6)
$$

Consider the following function in $\R^d$:
$$
\wt{\v}_\eps = \wt{\u}_0 + \eps \Lambda^\eps S_\eps (b(\D) \wt{\u}_0) =
\wt{\u}_0 + \eps K_\lambda(\eps) \wt{\u}_0,
\eqno(4.7)
$$
Here $K_\lambda(\eps)$ is the operator (1.25). Then
$$
\v_\eps = \wt{\v}_\eps \vert_{\O}.
\eqno(4.8)
$$

\smallskip\noindent\textbf{Theorem 4.2.} \textit{Suppose that the assumptions of Theorem} 4.1
\textit{are satisfied. Let} $\mathbf{v}_{\varepsilon}$
\textit{be the function defined by} (4.5), (4.6).
\textit{Then there exists a number} $\eps_1 \in (0,1]$
\textit{depending on the domain $\O$ and the lattice $\Gamma$ such that}
$$
\| \mathbf{u}_{\varepsilon} - \mathbf{v}_{\varepsilon}
\|_{H^1 (\mathcal{O}; \mathbb{C}^{n})} \leq {\mathcal C}(\lambda) \varepsilon^{1/2}
\| \mathbf{F} \|_{L_2 (\mathcal{O}; \mathbb{C}^{n})},\quad 0< \eps \leq \eps_1,
\eqno(4.9)
$$
\textit{or, in operator terms,}
$$
\| (\mathcal{A}_{N,\varepsilon}+\lambda I)^{-1} - (\mathcal{A}^0_N + \lambda I)^{-1} -
\varepsilon K_{N,\lambda} (\varepsilon) \|_{L_2 (\mathcal{O}; \mathbb{C}^{n}) \to
H^1 (\mathcal{O}; \mathbb{C}^{n})} \leq {\mathcal C}(\lambda) {\varepsilon}^{1/2}.
\eqno(4.10)
$$
\textit{For the flux} $\p_\eps:= g^\eps b(\D)\u_\eps$ \textit{we have}
$$
\| \p_\eps - \wt{g}^\eps S_\eps b(\D) \wt{\u}_0 \|_{L_2(\O;\C^m)}
\leq {\mathcal C}'(\lambda) \eps^{1/2} \| \mathbf{F} \|_{L_2 (\mathcal{O}; \mathbb{C}^{n})},\ \ 0< \eps \leq \eps_1,
\eqno(4.11)
$$
\textit{where} $\wt{g}(\x)$ \textit{is the matrix} (1.9).
\textit{The constants} ${\mathcal C}(\lambda)$ \textit{and} ${\mathcal C}'(\lambda)$
\textit{depend on} $m$, $d$, $\alpha_0$, $\alpha_1$, $\| g \|_{L_\infty}$,
$\| g^{-1} \|_{L_\infty}$, $\lambda$, \textit{the parameters of the lattice} $\Gamma$,
\textit{the constants ${\mathcal C}_1$ and ${\mathcal C}_2$ from the inequality} (2.2),
\textit{and the domain} $\mathcal{O}$.

\smallskip
Now we distinguish the case where the corrector is equal to zero.
Next statement follows from Theorem 4.2 and Proposition 1.2.

\smallskip\noindent\textbf{Proposition 4.3.}
\textit{Let} $\mathbf{u}_{\varepsilon}$ \textit{be the solution of the problem} (2.8), \textit{and let}
$\mathbf{u}_0$ \textit{be the solution of the problem} (2.28).
 \textit{If} $g^0 = \overline{g}$,
\textit{i.~e., relations} (1.13) \textit{are satisfied, then} $\Lambda=0$,
$K_{N,\lambda}(\eps)=0$, \textit{and we have}
$$
\| \mathbf{u}_{\varepsilon} - \mathbf{u}_0
\|_{H^1 (\mathcal{O}; \mathbb{C}^{n})} \leq {\mathcal C}(\lambda) \varepsilon^{1/2}
\| \mathbf{F} \|_{L_2 (\mathcal{O}; \mathbb{C}^{n})},\ \ 0< \eps \leq \eps_1.
$$

\smallskip\noindent\textbf{4.2. The first part of the proof:
introduction of the boundary layer correction term.}
The proof of Theorems 4.1 and 4.2 relies on application of the results
for the homogenization problem in $\R^d$ (Theorems 1.6, 1.7, 1.8)
and on estimation of the boundary layer correction term.

Consider the function $\wt{\u}_0 = P_\O \u_0 \in H^2(\R^d;\C^n)$. Clearly,
$$
\wt{\FF}:= \mathcal{A}^0 \widetilde{\mathbf{u}}_0 + \widetilde{\mathbf{u}}_0
\in L_2 (\mathbb{R}^d; \mathbb{C}^{n}).
$$
Herewith, $\wt{\FF} \vert_{\O} = \FF$.
Using the Fourier transformation and (1.4), (1.15), we obtain
$$
\begin{aligned}
&\| \wt{\FF} \|^2_{L_2 (\mathbb{R}^d)}
= \int_{\mathbb{R}^d} \left| (b (\boldsymbol{\xi})^*
g^0 b (\boldsymbol{\xi}) + \mathbf{1}) \widehat{\mathbf{u}}_0 (\boldsymbol{\xi}) \right|^2
\, d \boldsymbol{\xi}
\\
&\leq \ \int_{\mathbb{R}^d} (\alpha_1 |g^0| |\boldsymbol{\xi}|^2 + 1)^2
| \widehat{\mathbf{u}}_0 (\boldsymbol{\xi}) |^2 \, d \boldsymbol{\xi}
\leq \left( \textrm{max} \left\{ \alpha_1 \| g \|_{L_\infty}, 1 \right\} \right)^2
\| \widetilde{\mathbf{u}}_0 \|^2_{H^2 (\mathbb{R}^d)}.
\end{aligned}
\eqno(4.12)
$$
Here $\widehat{\mathbf{u}}_0 (\boldsymbol{\xi})$ is the Fourier-image of the function
$\widetilde{\mathbf{u}}_0 (\mathbf{x})$.
Combining (4.12) and (4.4), (2.33), we arrive at
$$
\begin{aligned}
\| \wt{\FF} \|_{L_2(\R^d)} &\le {\mathcal C}_3(\lambda) \|\FF\|_{L_2(\O)},
\cr
{\mathcal C}_3(\lambda) &= \wh{C}_\lambda C_\O \max \{ \alpha_1 \|g\|_{L_\infty}, 1\}.
\end{aligned}
\eqno(4.13)
$$

Let $\wt{\u}_{\varepsilon}\in H^1(\R^d;\C^n)$ be the generalized solution of the equation
$$
\mathcal{A}_{\varepsilon} \wt{\u}_{\varepsilon} + \lambda \wt{\u}_{\varepsilon}
= \wt{\FF}.
\eqno(4.14)
$$
Applying Theorems 1.6, 1.7, and 1.8, we obtain
$$
\| \wt{\u}_{\varepsilon} - \widetilde{\mathbf{u}}_0
\|_{L_2 (\mathbb{R}^d)}  \leq  C_1(\lambda) \varepsilon
\| \wt{\FF} \|_{L_2 (\mathbb{R}^d)}, \ \ \eps >0,
\eqno(4.15)
$$
$$
\| \wt{\u}_{\varepsilon} - \wt{\mathbf{v}}_{\varepsilon}
\|_{H^1 (\mathbb{R}^d)} \leq C_2(\lambda) \varepsilon \| \wt{\FF}
\|_{L_2 (\mathbb{R}^d)}, \ \ \eps >0,
\eqno(4.16)
$$
$$
\| g^\eps b(\D)\wt{\u}_{\varepsilon} - \wt{g}^\eps S_\eps b(\D)\wt{\u}_0
\|_{L_2 (\mathbb{R}^d)} \leq C_5(\lambda) \varepsilon \| \wt{\FF}
\|_{L_2 (\mathbb{R}^d)}, \ \ \eps >0.
\eqno(4.17)
$$
Here $\wt{\v}_\eps$ is defined by (4.7).

Consider the difference $\u_\eps - \wt{\u}_\eps$ in the domain $\O$.
By (2.8) and (4.14), this difference is the solution of the Neumann problem for the
homogeneous equation:
$$
\begin{aligned}
&(\A_\eps + \lambda I) (\u_\eps - \wt{\u}_\eps) =0 \ \ \text{in} \ \O;
\cr
&\partial^\eps_{\bnu}(\u_\eps - \wt{\u}_\eps) = - \bpsi_\eps\ \ \text{on}\ \partial \O,
\end{aligned}
\eqno(4.18)
$$
where $\bpsi_\eps := \partial^\eps_{\bnu} \wt{\u}_\eps$.
Since $\wt{\u}_\eps\vert_\O \in H^1(\O;\C^n)$ and
$\A_\eps \wt{\u}_\eps = \FF - \lambda \wt{\u}_\eps \in L_2(\O;\C^n)$,
we can apply the definition of the conormal derivative from Subsection
2.4. According to (2.20) and (2.23), the element
$\bpsi_\eps \in H^{-1/2}(\partial \O;\C^n)$ is defined by
$$
\begin{aligned}
&(\bpsi_\eps, \bvarphi)_{L_2(\partial \O)} =
\int_{\O} \left( \langle g^\eps b(\D) \wt{\u}_\eps, b(\D) T \bvarphi\rangle
+ \lambda \langle \wt{\u}_\eps, T \bvarphi \rangle\right) \,d\x -
\int_{\O} \langle \FF, T \bvarphi\rangle \,d\x,
\cr
&\bvarphi \in H^{1/2}(\partial\O;\C^n).
\end{aligned}
\eqno(4.19)
$$
Here $T$ is the operator (2.15) defined in Subsection 2.3.

We wish to approximate the difference $\wt{\u}_\eps - \u_\eps$ by a function
$\w_\eps$, which is solution of the Neumann problem for the homogeneous equation
with a simpler given function $\brho_\eps$ in the boundary condition.
The choice of $\brho_\eps$ is prompted by approximation of the solution (4.15)
and approximation of the flux (4.17).
We define $\brho_\eps \in H^{-1/2}(\partial \O;\C^n)$ by the relation
$$
\begin{aligned}
&(\brho_\eps, \bvarphi)_{L_2(\partial \O)} =
\int_{\O} \left( \langle \wt{g}^\eps S_\eps b(\D) \wt{\u}_0, b(\D) T \bvarphi\rangle
+ \lambda \langle {\u}_0, T \bvarphi \rangle
- \langle \FF, T \bvarphi\rangle \right) \,d\x,
\cr
&\bvarphi \in H^{1/2}(\partial\O;\C^n).
\end{aligned}
\eqno(4.20)
$$
Clearly, the right-hand side of (4.20) is an antilinear continuous functional of
$\bvarphi \in H^{1/2}(\partial\O;\C^n)$, so that definition (4.20) is correct.
Let us estimate the norm $\| \bpsi_\eps - \brho_\eps \|_{H^{-1/2}(\partial \O)}$.
From (4.19) and (4.20) it follows that
$$
\begin{aligned}
&(\bpsi_\eps - \brho_\eps, \bvarphi)_{L_2(\partial \O)}
= \int_{\O}  \langle g^\eps b(\D) \wt{\u}_\eps - \wt{g}^\eps S_\eps b(\D) \wt{\u}_0,
b(\D) T \bvarphi\rangle \,d\x
\cr
&+ \lambda \int_{\O} \langle \wt{\u}_\eps - {\u}_0, T \bvarphi \rangle \,d\x,
\quad \bvarphi \in H^{1/2}(\partial\O;\C^n).
\end{aligned}
$$
Together with (4.15) and (4.17) this implies
$$
\begin{aligned}
&|(\bpsi_\eps - \brho_\eps, \bvarphi)_{L_2(\partial \O)}|
\cr
&\le \eps \|\wt{\FF}\|_{L_2(\R^d)} \left(C_5(\lambda) \|b(\D) T \bvarphi\|_{L_2(\O)}
+ \lambda C_1(\lambda) \|T\bvarphi\|_{L_2(\O)}\right),
\cr
&\bvarphi \in H^{1/2}(\partial\O;\C^n),\quad \eps >0.
\end{aligned}
\eqno(4.21)
$$
Combining (4.21) with (1.2), (1.5), (2.16), and (4.13), we obtain
$$
\begin{aligned}
\|(\bpsi_\eps - \brho_\eps, \bvarphi)_{L_2(\partial \O)}|
\le {\mathcal C}_4(\lambda) \eps \|{\FF}\|_{L_2(\O)} \|\bvarphi\|_{H^{1/2}(\partial \O)},
\cr
 \bvarphi \in H^{1/2}(\partial\O;\C^n),\quad \eps>0,
\end{aligned}
$$
where
${\mathcal C}_4(\lambda) =
{\mathcal C}_3(\lambda) \wt{c}_1 \left(C_5(\lambda) (d \alpha_1)^{1/2} + \lambda C_1(\lambda) \right)$.
From here and (2.22) it follows that
$$
\| \bpsi_\eps - \brho_\eps \|_{H^{-1/2}(\partial \O)}
\le {\mathcal C}_4(\lambda) \eps \|{\FF}\|_{L_2(\O)},\quad \eps >0.
\eqno(4.22)
$$

Now we introduce the correction term $\w_\eps \in H^1(\O;\C^n)$ which is the generalized
solution of the Neumann problem
$$
\begin{aligned}
&(\A_\eps + \lambda I) \w_\eps = 0 \ \ \text{in}\ \O;
\cr
&\partial^\eps_{\bnu}\w_\eps =  \brho_\eps\ \ \text{on}\ \partial \O.
\end{aligned}
\eqno(4.23)
$$

\smallskip\noindent\textbf{Lemma 4.4.} \textit{Let $\u_\eps$ be the solution of the problem} (2.8),
\textit{and let $\wt{\u}_\eps$ be the solution of the equation} (4.14). \textit{Let
$\w_\eps$ be the solution of the problem} (4.23).
\textit{Then for $\eps>0$ we have}
$$
\|\u_\eps - \wt{\u}_\eps + \w_\eps\|_{H^1(\O)} \le {\mathcal C}_5(\lambda) \eps \|\FF\|_{L_2(\O)},
\eqno(4.24)
$$
\textit{where} ${\mathcal C}_5(\lambda) = c_\lambda^{-1} \wh{c}_1 {\mathcal C}_4(\lambda)$.

\smallskip\noindent\textbf{Proof.}
By (4.18) and (4.23), the function $\u_\eps - \wt{\u}_\eps + \w_\eps$
is the solution of the problem
$$
\begin{aligned}
&(\A_\eps + \lambda I) (\u_\eps - \wt{\u}_\eps + \w_\eps) = 0 \ \ \text{in}\ \O;
\cr
&\partial^\eps_{\bnu}(\u_\eps - \wt{\u}_\eps + \w_\eps) =  \brho_\eps - \bpsi_\eps \ \
\text{on}\ \partial \O.
\end{aligned}
$$
Applying Proposition 2.4, we obtain
$$
\|\u_\eps - \wt{\u}_\eps + \w_\eps\|_{H^1(\O)} \le c_\lambda^{-1} \wh{c}_1 \|\brho_\eps - \bpsi_\eps
\|_{H^{-1/2}(\partial \O)}.
$$
Together with (4.22) this implies (4.24). $\ \bullet$

 \smallskip
 \textbf{Conclusions.} \textbf{1.} Combining (4.8), (4.13), (4.16), and (4.24), we arrive at
 $$
\|\u_\eps - {\v}_\eps + \w_\eps\|_{H^1(\O)} \le \left({\mathcal C}_5(\lambda)+
C_2(\lambda){\mathcal C}_3(\lambda)\right) \eps \|\FF\|_{L_2(\O)},\quad \eps>0.
\eqno(4.25)
 $$
Hence, it is clear that, \textit{in order to prove main estimate} (4.9)
\textit{from Theorem} 4.2,
\textit{it suffices to estimate $\|\w_\eps\|_{H^{1}(\O)}$ by} $C \eps^{1/2} \|\FF\|_{L_2(\O)}$.

\textbf{2.} Replacing the $H^1(\O;\C^n)$-norm by
the $L_2(\O;\C^n)$-norm in the left-hand side of (4.24), we obtain:
$$
\|\u_\eps - \wt{\u}_\eps + \w_\eps\|_{L_2(\O)} \le {\mathcal C}_5(\lambda) \eps \|\FF\|_{L_2(\O)},
\quad \eps>0.
\eqno(4.26)
$$
Combining (4.26) and (4.13), (4.15), we arrive at
$$
\|\u_\eps - {\u}_0 + \w_\eps\|_{L_2(\O)} \le \left({\mathcal C}_5(\lambda)+
C_1(\lambda) {\mathcal C}_3(\lambda) \right) \eps \|\FF\|_{L_2(\O)},\quad \eps>0.
\eqno(4.27)
$$
This shows that, \textit{in order to prove estimate} (4.1) \textit{from Theorem} 4.1,
\textit{it suffices to estimate $\|\w_\eps\|_{L_2(\O)}$ in terms of} $C \eps \|\FF\|_{L_2(\O)}$.

 Thus, the problem is reduced to estimation of the function $\w_\eps$
 which can be interpreted as the
 \textit{boundary layer correction term}. First we estimate the norm of $\w_\eps$ in $H^1(\O;\C^n)$,
 and so prove Theorem 4.2 (see \S 5). Next, using already proved Theorem 4.2, we will estimate
 the norm of $\w_\eps$ in $L_2(\O;\C^n)$ (see \S 6).

\section*{\S 5. Proof of Theorem 4.2}

\noindent\textbf{5.1.}
According to Subsection 2.4, the generalized solution $\w_\eps\in H^1(\O;\C^n)$ of the problem (4.23)
satisfies the identity
$$
\int_{\O} \left( \langle g^\eps b(\D)\w_\eps, b(\D) \eeta \rangle + \lambda \langle \w_\eps, \eeta
\rangle \right)\,d\x =  (\brho_\eps, \gamma \eeta)_{L_2(\partial \O)},
\  \forall \eeta \in H^1(\O;\C^n).
\eqno(5.1)
$$
Using (4.20) and (2.30), we represent the right-hand side of (5.1) as
$$
\begin{aligned}
(\brho_\eps, \gamma \eeta)_{L_2(\partial \O)} &= \int_{\O} \left(
 \langle \wt{g}^\eps S_\eps b(\D) \wt{\u}_0, b(\D) T \gamma \eeta \rangle
+ \lambda \langle {\u}_0, T \gamma \eeta \rangle  - \langle \FF, T \gamma \eeta \rangle \right) \,d\x
\cr
&=  \int_{\O}  \langle \wt{g}^\eps S_\eps b(\D) \wt{\u}_0 - {g}^0 b(\D) {\u}_0, b(\D) T \gamma \eeta \rangle.
\end{aligned}
\eqno(5.2)
$$
From (5.1) and (5.2) it follows that
$$
\intop_{\O} \left( \langle g^\eps b(\D)\w_\eps, b(\D) \eeta \rangle + \lambda \langle \w_\eps, \eeta
\rangle \right)\,d\x = {\mathcal I}_\eps[\eeta],
\quad \forall \eeta \in H^1(\O;\C^n),
\eqno(5.3)
$$
where we denote
$$
{\mathcal I}_\eps[\eeta]:=
\int_{\O}  \langle \wt{g}^\eps S_\eps b(\D) \wt{\u}_0 - {g}^0 b(\D) {\u}_0,
b(\D) T \gamma \eeta \rangle\, d\x,\quad \eeta \in H^1(\O;\C^n).
\eqno(5.4)
$$

For further estimates, it is convenient to represent the functional (5.4) as
$$
{\mathcal I}_\eps[\eeta]= {\mathcal I}^{(1)}_\eps[\eeta] + {\mathcal I}^{(2)}_\eps[\eeta],
\eqno(5.5)
$$
where
$$
{\mathcal I}^{(1)}_\eps[\eeta] =
\int_{\O}  \langle {g}^0 S_\eps b(\D) \wt{\u}_0 - {g}^0 b(\D) {\u}_0, b(\D) T \gamma \eeta \rangle \,d\x,
\eqno(5.6)
$$
$$
{\mathcal I}^{(2)}_\eps[\eeta] =
\int_{\O}  \langle (\wt{g}^\eps - g^0)S_\eps b(\D) \wt{\u}_0, b(\D) T \gamma \eeta \rangle \, d\x.
\eqno(5.7)
$$

The term (5.6) can be easily estimated with the help of Proposition 1.4 and relations
(1.2), (1.5), and (1.15):
$$
\begin{aligned}
&|{\mathcal I}^{(1)}_\eps[\eeta]| \le |g^0|
\| (S_\eps -I) b(\D)\wt{\u}_0\|_{L_2(\R^d)} \|b(\D) T \gamma \eeta\|_{L_2(\O)}
\cr
&\le
\eps \|g\|_{L_\infty} r_1 \|\D b(\D)\wt{\u}_0\|_{L_2(\R^d)} (d \alpha_1)^{1/2} \|T \gamma \eeta\|_{H^1(\O)},
\quad \eps >0.
\end{aligned}
\eqno(5.8)
$$
By (1.4), (2.33), and (4.4),
$$
\begin{aligned}
\|\D b(\D)\wt{\u}_0\|_{L_2(\R^d)} &\le \alpha_1^{1/2} \|\wt{\u}_0 \|_{H^2(\R^d)}
\cr
&\le \alpha_1^{1/2} C_\O \| {\u}_0 \|_{H^2(\O)}
\le \alpha_1^{1/2} C_\O \wh{C}_\lambda \| {\FF} \|_{L_2(\O)}.
\end{aligned}
\eqno(5.9)
$$
Now, combining (5.8), (5.9), (2.13), and (2.16), we arrive at
$$
\begin{aligned}
&|{\mathcal I}^{(1)}_\eps[\eeta]| \le {\mathcal C}_6(\lambda) \eps \|\FF\|_{L_2(\O)} \|\eeta\|_{H^1(\O)},
\quad \eeta \in H^1(\O;\C^n),\ \ \eps>0,
\cr
&{\mathcal C}_6(\lambda)= d^{1/2} r_1 \alpha_1 \|g\|_{L_\infty}
 C_\O \wh{C}_\lambda \wh{c}_1 \wt{c}_1.
\end{aligned}
\eqno(5.10)
$$

\smallskip\noindent\textbf{5.2. Analysis of the term ${\mathcal I}^{(2)}_\eps$.}
Transform the term (5.7), using (1.2):
$$
{\mathcal I}^{(2)}_\eps[\eeta]
= \int_{\O} \sum_{l=1}^d \langle f_l^\eps S_\eps b(\D) \wt{\u}_0, D_l T \gamma \eeta \rangle \, d\x,
\eqno(5.11)
$$
where we denote
$$
f_l(\x) := b_l^* (\wt{g}(\x) - g^0) = b_l^* (g(\x)(b(\D)\Lambda(\x) + \1_m) - g^0),\ \ l=1,\dots,d.
\eqno(5.12)
$$
Then $f_l(\x)$ are $\Gamma$-periodic $(n\times m)$-matrix-valued functions,
$f_l \in L_2(\Omega)$, and, by the definition (1.8) of the effective matrix,
we have $\overline{f_l}=0$. Finally, from the equation (1.7) for $\Lambda$ it follows that
$$
\sum_{l=1}^d D_l f_l(\x) = b(\D)^* g(\x) (b(\D)\Lambda(\x) + \1_m)=0.
$$
Thus, the set of functions $f_l(\x)$, $l=1,\dots,d,$ satisfy the assumptions of Lemma 3.1.
By this lemma, there exist matrix-valued functions $M_{lj}(\x)$, $l,j=1,\dots,d$, such that
relations (3.2), (3.3) are valid. Then
$$
f_l^\eps(\x) = \eps \sum_{j=1}^d \partial_j M_{lj}^\eps(\x).
\eqno(5.13)
$$

We will need estimates for the $L_2(\Omega)$-norms of the functions (5.12).
From (1.5), taking into account that $\wt{g}-g^0$ is the orthogonal projection of
$\wt{g}$ onto the space orthogonal to constants, we obtain
$$
\|f_l\|_{L_2(\Omega)} \le \alpha_1^{1/2} \|\wt{g}-g^0\|_{L_2(\Omega)}
\le \alpha_1^{1/2} \|\wt{g}\|_{L_2(\Omega)}.
$$
Combining this with (1.2), (1.5), (1.9), and (1.10), we have
$$
\|f_l\|_{L_2(\Omega)} \le \alpha_1^{1/2} \|g\|_{L_\infty} (|\Omega|^{1/2}+ (d\alpha_1)^{1/2} \|\D \Lambda\|_{L_2(\Omega)})
\le  |\Omega|^{1/2} {\mathfrak C},\  l=1,\dots,d,
\eqno(5.14)
$$
where ${\mathfrak C} =
\alpha_1^{1/2} \|g\|_{L_\infty} (1+ (d m)^{1/2} \alpha_1^{1/2}
\alpha_0^{-1/2} \|g\|_{L_\infty}^{1/2} \|g^{-1}\|_{L_\infty}^{1/2})$.
Hence, by (3.4),
$$
\|M_{lj}\|_{L_2(\Omega)} \le r_0^{-1} |\Omega|^{1/2} {\mathfrak C},\quad l,j=1,\dots,d.
\eqno(5.15)
$$

By (5.13),
$$
\begin{aligned}
&f_l^\eps S_\eps b(\D) \wt{\u}_0 = \eps
\sum_{j=1}^d \left(\partial_j M_{lj}^\eps \right) S_\eps b(\D) \wt{\u}_0
\cr
&= \eps \sum_{j=1}^d \partial_j \left(M_{lj}^\eps S_\eps b(\D) \wt{\u}_0\right) -
\eps \sum_{j=1}^d M_{lj}^\eps S_\eps b(\D)\partial_j \wt{\u}_0.
\end{aligned}
\eqno(5.16)
$$
According to (5.16), the term (5.11) can be represented as
$$
{\mathcal I}^{(2)}_\eps[\eeta] = \wt{\mathcal I}^{(2)}_\eps[\eeta] + \wh{\mathcal I}^{(2)}_\eps[\eeta],
\eqno(5.17)
$$
where
$$
\wt{\mathcal I}^{(2)}_\eps[\eeta]=
\eps \int_{\O} \sum_{l,j=1}^d
\langle \partial_j \left(M_{lj}^\eps S_\eps b(\D) \wt{\u}_0 \right),
D_l T \gamma \eeta \rangle \, d\x,
\eqno(5.18)
$$
$$
\wh{\mathcal I}^{(2)}_\eps[\eeta]=
- \eps \int_{\O} \sum_{l,j=1}^d
\langle M_{lj}^\eps S_\eps b(\D) \partial_j \wt{\u}_0,
D_l T \gamma \eeta \rangle \, d\x.
\eqno(5.19)
$$

The term (5.19) can be easily estimated by using Proposition 1.5:
$$
|\wh{\mathcal I}^{(2)}_\eps[\eeta]| \le \eps \sum_{l,j=1}^d |\Omega|^{-1/2} \|M_{lj}\|_{L_2(\Omega)}
\|b(\D) \partial_j \wt{\u}_0\|_{L_2(\R^d)} \|D_l T\gamma \eeta\|_{L_2(\O)}.
$$
Together with (2.13), (2.16), (5.9), and (5.15) this implies:
$$
\begin{aligned}
&|\wh{\mathcal I}^{(2)}_\eps[\eeta]| \le {\mathcal C}_7(\lambda) \eps \|\FF\|_{L_2(\O)} \|\eeta\|_{H^1(\O)},
\quad \eeta \in H^1(\O;\C^n),\ \ \eps>0,
\cr
&{\mathcal C}_7(\lambda) = d r_0^{-1} {\mathfrak C}
\alpha_1^{1/2} C_\O \wh{C}_\lambda  \wt{c}_1 \wh{c}_1.
\end{aligned}
\eqno(5.20)
$$

\smallskip\noindent\textbf{5.3. Estimate for the term $\wt{\mathcal I}^{(2)}_\eps$.}
Recall the notation $(\partial \O)_\eps$ for the $\eps$-neighborhood of
the boundary $\partial \O$ in $\R^d$.
Assume that $0< \eps \le \eps_1$, where the number $\eps_1$ was defined in Lemma 3.3.
We fix a smooth cut-off function $\theta_\eps(\x)$ in $\R^d$ such that
$$
\begin{aligned}
&\theta_\eps \in C_0^\infty(\R^d),\ \ \text{supp}\, \theta_\eps \subset (\partial \O)_\eps,
\ \ 0 \le \theta_\eps(\x) \le 1,
\cr
&\theta_\eps(\x)\vert_{\partial \O}=1,\ \ \eps |\nabla \theta_\eps(\x)| \le \kappa = \text{const}.
\end{aligned}
\eqno(5.21)
$$
Obviously, such a function exists. The constant $\kappa$ depends only on the domain $\O$.

Note that
$$
\int_{\O} \sum_{l,j=1}^d
\langle \partial_j \left(M_{lj}^\eps (1-\theta_\eps) S_\eps b(\D) \wt{\u}_0 \right),
D_l T \gamma \eeta \rangle \, d\x =0,\quad \eeta \in H^1(\O;\C^n).
\eqno(5.22)
$$
Indeed, since the left-hand side of (5.22) is an antilinear continuous functional of
$\eeta \in H^1(\O;\C^n)$, it suffices to prove this identity for some dense set in $H^1(\O;\C^n)$.
Assume that $\eeta \in H^2(\O;\C^n)$, and, integrating by parts, rewrite the left-hand side of
(5.22) as
$$
- i \int_{\O}
\langle (1-\theta_\eps) S_\eps b(\D) \wt{\u}_0,
\sum_{l,j=1}^d
(M_{lj}^\eps)^*\partial_j \partial_l T \gamma \eeta \rangle \, d\x.
$$
This expression is equal to zero since
$\sum_{l,j=1}^d (M_{lj}^\eps)^*\partial_j \partial_l =0$ due to
$M_{lj}= - M_{jl}$ (see (3.2)).

By (5.22), the term (5.18) takes the form
$$
\wt{\mathcal I}^{(2)}_\eps[\eeta]=
\eps \int_{\O} \sum_{l,j=1}^d
\langle \partial_j \left( \theta_\eps M_{lj}^\eps S_\eps b(\D) \wt{\u}_0 \right),
D_l T \gamma \eeta \rangle \, d\x.
\eqno(5.23)
$$

Consider the following expression
$$
\begin{aligned}
&\eps \sum_{j=1}^d \partial_j \left( \theta_\eps M_{lj}^\eps S_\eps b(\D) \wt{\u}_0 \right)
=
\eps \theta_\eps \sum_{j=1}^d M_{lj}^\eps S_\eps (b(\D) \partial_j \wt{\u}_0)
\cr
&+
\eps \sum_{j=1}^d (\partial_j \theta_\eps) M_{lj}^\eps S_\eps b(\D) \wt{\u}_0
+ \theta_\eps \sum_{j=1}^d (\partial_j M_{lj})^\eps S_\eps b(\D) \wt{\u}_0.
\end{aligned}
$$
Take into account that $\sum_{j=1}^d \partial_j M_{lj} = f_l$ (see (3.3)).
Then
$$
\eps \left\|\sum_{j=1}^d \partial_j \left( \theta_\eps M_{lj}^\eps S_\eps b(\D) \wt{\u}_0 \right)\right\|_{L_2(\O)}
\le J_l^{(1)}(\eps) + J_l^{(2)}(\eps) + J_l^{(3)}(\eps),
\eqno(5.24)
$$
where
$$
J_l^{(1)}(\eps) = \eps \sum_{j=1}^d \|\theta_\eps M_{lj}^\eps S_\eps (b(\D) \partial_j \wt{\u}_0)\|_{L_2(\O)},
\eqno(5.25)
$$
$$
J_l^{(2)}(\eps) = \eps \sum_{j=1}^d \| (\partial_j \theta_\eps) M_{lj}^\eps S_\eps b(\D) \wt{\u}_0\|_{L_2(\O)},
\eqno(5.26)
$$
$$
J_l^{(3)}(\eps) = \|\theta_\eps f_{l}^\eps S_\eps b(\D) \wt{\u}_0\|_{L_2(\O)}.
\eqno(5.27)
$$

In order to estimate the term (5.25), apply (5.21), Proposition 1.5, and relations (5.9), (5.15):
$$
J_l^{(1)}(\eps)
\le \eps \sum_{j=1}^d |\Omega|^{-1/2} \|M_{lj}\|_{L_2(\Omega)}
\| b(\D) \partial_j \wt{\u}_0 \|_{L_2(\R^d)}
\le {\mathfrak C}^{(1)}(\lambda) \eps \|\FF\|_{L_2(\O)},
\eqno(5.28)
$$
where ${\mathfrak C}^{(1)}(\lambda)=
r_0^{-1} {\mathfrak C} d^{1/2} \alpha_1^{1/2} C_\O \wh{C}_\lambda$.

The term (5.26) is estimated with the help of (5.21) and Lemma 3.3:
$$
\begin{aligned}
&(J_l^{(2)}(\eps))^2 \le d \kappa^2 \sum_{j=1}^d \int_{(\partial \O)_\eps} |M^\eps_{lj}|^2 |S_\eps b(\D)\wt{\u}_0|^2\,d\x
\cr
&\le
\eps d \kappa^2 \beta_* |\Omega|^{-1} \sum_{j=1}^d \|M_{lj}\|^2_{L_2(\Omega)}
\| S_\eps b(\D) \wt{\u}_0 \|^2_{H^1(\R^d)},\quad 0< \eps \le \eps_1.
\end{aligned}
\eqno(5.29)
$$
From (1.19), (1.4), (2.33), and (4.4) it follows that
$$
\begin{aligned}
&\| S_\eps b(\D) \wt{\u}_0 \|_{H^1(\R^d)} \le \| b(\D) \wt{\u}_0 \|_{H^1(\R^d)}
\cr
&\le \alpha_1^{1/2}\|\wt{\u}_0\|_{H^2(\R^d)}
\le \alpha_1^{1/2} C_\O \wh{C}_\lambda  \|{\FF}\|_{L_2(\O)}.
\end{aligned}
\eqno(5.30)
$$
Combining (5.29), (5.30), and (5.15), we obtain
$$
J_l^{(2)}(\eps)\le {\mathfrak C}^{(2)}(\lambda) \eps^{1/2} \| \FF \|_{L_2(\O)},
\ \ 0< \eps \le\eps_1,
\eqno(5.31)
$$
where
${\mathfrak C}^{(2)}(\lambda) = \kappa d \beta_*^{1/2} r_0^{-1}{\mathfrak C}
\alpha_1^{1/2} C_\O \wh{C}_\lambda$.

The term (5.27) is estimated in a similar way by using (5.21) and Lemma 3.3:
$$
\begin{aligned}
&(J_l^{(3)}(\eps))^2 \le \int_{(\partial \O)_\eps} |f^\eps_{l}|^2 |S_\eps b(\D)\wt{\u}_0|^2\,d\x
\cr
&\le
\eps \beta_* |\Omega|^{-1} \|f_{l}\|^2_{L_2(\Omega)}
\| S_\eps b(\D) \wt{\u}_0 \|^2_{H^1(\R^d)}, \quad 0< \eps \le \eps_1.
\end{aligned}
\eqno(5.32)
$$
From (5.14), (5.30), and (5.32) it follows that
$$
J_l^{(3)}(\eps) \le {\mathfrak C}^{(3)}(\lambda) \eps^{1/2} \| \FF \|_{L_2(\O)},
\ \ 0< \eps \le\eps_1,
\eqno(5.33)
$$
where
${\mathfrak C}^{(3)}(\lambda) = \beta^{1/2}_*  {\mathfrak C}
\alpha_1^{1/2} C_\O \wh{C}_\lambda$.

Finally, relations (5.24), (5.28), (5.31), and (5.33) imply that
$$
\eps \left\|\sum_{j=1}^d \partial_j \left( \theta_\eps M_{lj}^\eps S_\eps b(\D) \wt{\u}_0 \right)\right\|_{L_2(\O)}
\le \wt{\mathfrak C}(\lambda)
\eps^{1/2} \|\FF \|_{L_2(\O)},\quad 0< \eps \le \eps_1,
\eqno(5.34)
$$
where $\wt{\mathfrak C}(\lambda)={\mathfrak C}^{(1)}(\lambda) + {\mathfrak C}^{(2)}(\lambda) + {\mathfrak C}^{(3)}(\lambda)$.
Combining (2.13), (2.16), (5.23), and (5.34), we arrive at
$$
|\wt{\mathcal I}^{(2)}_\eps[\eeta]| \le {\mathcal C}_8(\lambda) \eps^{1/2} \|\FF\|_{L_2(\O)} \|\eeta\|_{H^1(\O)},
\quad \eeta \in H^1(\O;\C^n),\ \ 0 < \eps \le \eps_1,
\eqno(5.35)
$$
where ${\mathcal C}_8(\lambda) =\wt{\mathfrak C}(\lambda) d^{1/2} \wt{c}_1 \wh{c}_1.$

\smallskip\noindent\textbf{5.4. Completion of the proof of Theorem 4.2.}
Now relations (5.5), (5.10), (5.17), (5.20), and (5.35) imply that
$$
|{\mathcal I}_\eps[\eeta]| \le {\mathcal C}_9(\lambda) \eps^{1/2} \|\FF\|_{L_2(\O)} \|\eeta\|_{H^1(\O)},
\quad \eeta \in H^1(\O;\C^n),\ \ 0 < \eps \le \eps_1,
\eqno(5.36)
$$
where ${\mathcal C}_9(\lambda) = {\mathcal C}_6(\lambda) + {\mathcal C}_7(\lambda) + {\mathcal C}_8(\lambda)$.

From (5.36) and (5.3) with $\eeta =\w_\eps$ it follows that
$$
\intop_\O \left( \langle g^\eps b(\D)\w_\eps, b(\D)\w_\eps\rangle +
\lambda |\w_\eps|^2 \right)\,d\x \le {\mathcal C}_9(\lambda) \eps^{1/2} \|\FF\|_{L_2(\O)} \|\w_\eps \|_{H^1(\O)}
$$
for $0< \eps \le \eps_1$. Together with (2.7) this yields the desired estimate
$$
\|\w_\eps\|_{H^1(\O)} \le c_\lambda^{-1}{\mathcal C}_9(\lambda) \eps^{1/2} \|\FF\|_{L_2(\O)},
\quad 0< \eps \le \eps_1.
\eqno(5.37)
$$

As a result, (4.25) and (5.37) imply (4.9) with the constant
$$
{\mathcal C}(\lambda) = {\mathcal C}_5(\lambda) + C_2(\lambda) {\mathcal C}_3(\lambda)
+  c_\lambda^{-1}{\mathcal C}_9(\lambda).
$$

It remains to obtain approximation (4.11) for the flux.
By (4.9), taking (1.2) and (1.5) into account, we have
$$
\|\p_\eps - g^\eps b(\D)\v_\eps\|_{L_2(\O)} \le \|g\|_{L_\infty} (d \alpha_1)^{1/2} {\mathcal C}(\lambda)
\eps^{1/2} \| \FF \|_{L_2(\O)},\quad 0 < \eps \le \eps_1.
\eqno(5.38)
$$
Let $\wt{\v}_\eps$ be defined by (4.7). Note that, by (4.8),
$$
\| g^\eps b(\D)\v_\eps - \wt{g}^\eps S_\eps b(\D) \wt{\u}_0 \|_{L_2(\O)} \le
\| g^\eps b(\D)\wt{\v}_\eps - \wt{g}^\eps S_\eps b(\D) \wt{\u}_0 \|_{L_2(\R^d)}.
\eqno(5.39)
$$
Let us use the proof of Theorem 1.8.
By (1.47),
$$
\| g^\eps b(\D)\wt{\v}_\eps - \wt{g}^\eps S_\eps b(\D) \wt{\u}_0 \|_{L_2(\R^d)}
\le C_6 \eps \| \D b(\D) \wt{\u}_0 \|_{L_2(\R^d)},\quad \eps>0.
\eqno(5.40)
$$
From (5.39), (5.40), and (5.9) we obtain
$$
\| g^\eps b(\D)\v_\eps - \wt{g}^\eps S_\eps b(\D) \wt{\u}_0 \|_{L_2(\O)} \le
C_6 \alpha_1^{1/2} C_\O \wh{C}_\lambda \eps \|\FF\|_{L_2(\O)},\ \ \eps>0.
\eqno(5.41)
$$
Finally, (5.38) and (5.41) imply the required estimate (4.11) with the constant
${\mathcal C}'(\lambda) = \|g\|_{L_\infty} (d \alpha_1)^{1/2} {\mathcal C}(\lambda)
+ C_6 \alpha_1^{1/2} C_\O \wh{C}_\lambda$. This completes the proof of Theorem 4.2. $\ \bullet$

\section*{\S 6. Proof of Theorem 4.1}

\noindent\textbf{6.1.} As was mentioned in the end of \S 4,
for the proof of Theorem 4.1 we need to estimate the norm
$\|\w_\eps\|_{L_2(\O)}$ by $C \eps \|\FF\|_{L_2(\O)}$.
Comparing (5.5), (5.10), (5.17), and (5.20), we see that
$$
\begin{aligned}
|{\mathcal I}_\eps[\eeta]| \le ({\mathcal C}_6(\lambda) + {\mathcal C}_7(\lambda))
\eps \|\FF\|_{L_2(\O)} \|\eeta\|_{H^1(\O)} +
|\wt{\mathcal I}^{(2)}_\eps[\eeta]|,
\cr \eeta \in H^1(\O;\C^n),\quad \eps >0.
\end{aligned}
\eqno(6.1)
$$
In the identity (5.3) which is valid for any function $\eeta \in H^1(\O;\C^n)$,
now we choose $\eeta$ in a special way. Namely, let $\bPhi \in L_2(\O;\C^n)$,
and let $\eeta_\eps \in H^1(\O;\C^n)$
be the generalized solution of the Neumann problem
$$
b(\D)^* g^\eps b(\D) \eeta_\eps + \lambda \eeta_\eps = \bPhi \ \text{in}\ \O;
\ \ \partial_{\bnu}^\eps \eeta_\eps \vert_{\partial \O}=0,
\eqno(6.2)
$$
i.~e., $\eeta_\eps = (\A_{N,\eps}+ \lambda I)^{-1}\bPhi$.
By (2.11),
$$
\|\eeta_\eps \|_{H^1(\O)} \le c_\lambda^{-1} \|\bPhi\|_{L_2(\O)}.
\eqno(6.3)
$$
By definition of the generalized solution of the problem (6.2), we have
$$
\int_\O
\left(\langle g^\eps b(\D) \w_\eps, b(\D) \eeta_\eps \rangle + \lambda \langle \w_\eps, \eeta_\eps \rangle\right) \, d\x = \int_\O \langle \w_\eps, \bPhi \rangle\,d\x.
\eqno(6.4)
$$
Combining (6.4) and (5.3), we obtain
$$
\int_\O \langle \w_\eps, \bPhi \rangle\,d\x =
{\mathcal I}_\eps[\eeta_\eps],\quad \bPhi \in L_2(\O;\C^n).
\eqno(6.5)
$$
Now (6.1), (6.3), and (6.5) imply that
$$
\begin{aligned}
|(\w_\eps, \bPhi)_{L_2(\O)}|\le ({\mathcal C}_6(\lambda) + {\mathcal C}_7(\lambda))
c_\lambda^{-1} \eps \|\FF\|_{L_2(\O)} \|\bPhi\|_{L_2(\O)} +
|\wt{\mathcal I}^{(2)}_\eps[\eeta_\eps]|,
\cr
\bPhi \in L_2(\O;\C^n),\quad \eps >0.
\end{aligned}
\eqno(6.6)
$$

\smallskip\noindent\textbf{6.2. Analysis of the term $\wt{\mathcal I}^{(2)}_\eps[\eeta_\eps]$.}
Clearly, the required estimate for $\|\w_\eps\|_{L_2(\O)}$ will be proved,
if we succeed to estimate the last term in the right-hand side of
(6.6) by $C \eps \|\FF\|_{L_2(\O)} \|\bPhi\|_{L_2(\O)}$.

We apply the (already proved) Theorem 4.2 in order to approximate the solution $\eeta_\eps$
of the problem (6.2) in the $H^1(\O;\C^n)$-norm.
Denote by $\eeta_0$ the solution of the corresponding homogenized problem:
$\eeta_0 = (\A^0_N + \lambda I)^{-1}\bPhi$.
Let $\wt{\eeta}_0 = P_\O \eeta_0$.
Note that, by (2.34),
$$
\|\wt{\eeta}_0\|_{H^2(\R^d)} \le C_\O \|\eeta_0 \|_{H^2(\O)} \le C_\O \wh{C}_\lambda \|\bPhi\|_{L_2(\O)}.
\eqno(6.7)
$$
From (4.10) it follows that
$$
\| \eeta_\eps - \eeta_0 - \eps \Lambda^\eps S_\eps b(\D) \wt{\eeta}_0 \|_{H^1(\O)}
\le {\mathcal C}(\lambda) \eps^{1/2} \|\bPhi \|_{L_2(\O)}, \quad 0< \eps \le \eps_1.
\eqno(6.8)
$$

Represent $\wt{\mathcal I}^{(2)}_\eps[\eeta_\eps]$ (see (5.23)) as the sum of three terms
$$
\wt{\mathcal I}^{(2)}_\eps[\eeta_\eps]= \wt{\mathcal I}^{(2)}_\eps[(\A_{N,\eps}+\lambda I)^{-1}\bPhi]=
 {\mathcal J}_\eps^{(1)}[\bPhi]
+ {\mathcal J}_\eps^{(2)}[\bPhi] + {\mathcal J}_\eps^{(3)}[\bPhi],
\eqno(6.9)
$$
where
$$
{\mathcal J}_\eps^{(1)}[\bPhi]=
\eps \int_{\O} \sum_{l,j=1}^d
\langle \partial_j \left( \theta_\eps M_{lj}^\eps S_\eps b(\D) \wt{\u}_0 \right),
D_l T \gamma (\eeta_\eps - \eeta_0 - \eps \Lambda^\eps S_\eps b(\D) \wt{\eeta}_0) \rangle \, d\x,
\eqno(6.10)
$$
$$
{\mathcal J}_\eps^{(2)}[\bPhi]=
\eps \int_{\O} \sum_{l,j=1}^d
\langle \partial_j \left( \theta_\eps M_{lj}^\eps S_\eps b(\D) \wt{\u}_0 \right),
D_l T \gamma \eeta_0 \rangle \, d\x,
\eqno(6.11)
$$
$$
{\mathcal J}_\eps^{(3)}[\bPhi]=
\eps \int_{\O} \sum_{l,j=1}^d
\langle \partial_j \left( \theta_\eps M_{lj}^\eps S_\eps b(\D) \wt{\u}_0 \right),
D_l T \gamma (\eps \Lambda^\eps S_\eps b(\D) \wt{\eeta}_0) \rangle \, d\x.
\eqno(6.12)
$$

The term (6.10) is easily estimated by using (6.8), (2.13), (2.16), and (5.34):
$$
\begin{aligned}
&|{\mathcal J}_\eps^{(1)}[\bPhi]| \le \eps^{1/2}\wt{\mathfrak C}(\lambda)
\|\FF\|_{L_2(\O)} d^{1/2} \wt{c}_1 \wh{c}_1
\| \eeta_\eps - \eeta_0 - \eps \Lambda^\eps S_\eps b(\D) \wt{\eeta}_0 \|_{H^1(\O)}
\cr
&\le
{\mathcal C}^{(1)}(\lambda) \eps \|\FF\|_{L_2(\O)} \|\bPhi\|_{L_2(\O)},
\quad
\bPhi \in L_2(\O;\C^n),\ 0< \eps \le \eps_1,
\end{aligned}
\eqno(6.13)
$$
where ${\mathcal C}^{(1)}(\lambda) = \wt{\mathfrak C}(\lambda)d^{1/2} \wt{c}_1 \wh{c}_1{\mathcal C}(\lambda).$

In order to estimate the term (6.11), we apply (5.21) and (5.34):
$$
|{\mathcal J}_\eps^{(2)}[\bPhi]| \le \eps^{1/2}\wt{\mathfrak C}(\lambda)
\|\FF\|_{L_2(\O)} \sum_{l=1}^d \|D_l T\gamma \eeta_0\|_{L_2(B_\eps)}.
\eqno(6.14)
$$
According to Lemma 3.2,
$$
\|D_l T\gamma \eeta_0\|_{L_2(B_\eps)}^2 \le \beta \eps \|D_l T \gamma \eeta_0 \|_{H^1(\O)}
\|D_l T \gamma \eeta_0 \|_{L_2(\O)},\quad 0< \eps \le \eps_0.
$$
Hence,
$$
\sum_{l=1}^d \|D_l T\gamma \eeta_0\|_{L_2(B_\eps)}
\le \eps^{1/2} (\beta d)^{1/2} \|T\gamma \eeta_0\|^{1/2}_{H^2(\O)}
\|T\gamma \eeta_0\|^{1/2}_{H^1(\O)},\quad 0< \eps \le \eps_0.
$$
Applying estimates (2.13), (2.16), (2.18), (2.32), and (2.34), we obtain
$$
\begin{aligned}
&\sum_{l=1}^d \|D_l T\gamma \eeta_0\|_{L_2(B_\eps)}
\le \eps^{1/2} (\beta d)^{1/2} (\wt{c}_2 \wh{c}_2 \wt{c}_1 \wh{c}_1)^{1/2}
\|\eeta_0\|^{1/2}_{H^2(\O)}
\|\eeta_0\|^{1/2}_{H^1(\O)}
\cr
&\le
\eps^{1/2} (\beta d)^{1/2} (\wt{c}_2 \wh{c}_2 \wt{c}_1 \wh{c}_1)^{1/2}
\wh{C}_\lambda^{1/2} c_\lambda^{-1/2} \|\bPhi\|_{L_2(\O)},\quad 0< \eps \le \eps_0.
\end{aligned}
\eqno(6.15)
$$
Now (6.14) and (6.15) imply that
$$
|{\mathcal J}_\eps^{(2)}[\bPhi]| \le {\mathcal C}^{(2)}(\lambda) \eps \|\FF\|_{L_2(\O)} \|\bPhi\|_{L_2(\O)},
\quad \bPhi \in L_2(\O;\C^n),\ \ 0 < \eps \le \eps_0,
\eqno(6.16)
$$
где ${\mathcal C}^{(2)}(\lambda) = \wt{\mathfrak C}(\lambda) (\beta d)^{1/2} (\wt{c}_2 \wh{c}_2 \wt{c}_1 \wh{c}_1)^{1/2}
\wh{C}_\lambda^{1/2} c_\lambda^{-1/2}$.

\smallskip\noindent\textbf{6.3. Estimate of the term ${\mathcal J}_\eps^{(3)}[\bPhi]$.}
It remains to estimate the term (6.12). Denote
$$
\z_\eps := T\gamma (\eps \Lambda^\eps S_\eps b(\D) \wt{\eeta}_0).
\eqno(6.17)
$$
By (5.34), the term (6.12) satisfies
$$
|{\mathcal J}^{(3)}_\eps[\bPhi]| \le \eps^{1/2} \wt{\mathfrak C}(\lambda) \|\FF\|_{L_2(\O)}
\sum_{l=1}^d \|D_l \z_\eps \|_{L_2(\O)},\quad 0< \eps \le \eps_1.
\eqno(6.18)
$$

By definitions of the operators $\gamma$ and $T$ (see Subsection 2.3),
the function (6.17) is the generalized solution of the Dirichlet problem
$$
\begin{aligned}
&-\Delta \z_\eps + \z_\eps =0 \ \ \text{in}\ \O;
\cr
&\z_\eps\vert_{\partial \O} = \eps \Lambda^\eps S_\eps b(\D) \wt{\eeta}_0 \vert_{\partial \O}.
\end{aligned}
\eqno(6.19)
$$

Now we estimate the norm of $\z_\eps$ in $H^1(\O;\C^n)$.
Consider the following function in $\R^d$:
$$
{\boldsymbol{\phi}}_\eps(\x)= \eps \theta_\eps(\x) \Lambda^\eps(\x) (S_\eps b(\D) \wt{\eeta}_0)(\x).
\eqno(6.20)
$$
Taking (5.21) into account, we can rewrite the problem (6.19)
in a different way: $-\Delta \z_\eps + \z_\eps =0$ in $\O$ and
$\z_\eps\vert_{\partial \O} = {\bphi}_\eps \vert_{\partial \O}$.
By Lemma 3.5,
$$
\|\z_\eps\|_{H^1(\O)} \le 2 \|\bphi_\eps\|_{H^1(\O)}.
\eqno(6.21)
$$
Thus, the question is reduced to the following statement.

 \smallskip\noindent\textbf{Lemma 6.1.} \textit{Let $\bphi_\eps$ be the function} (6.20).
\textit{Then we have}
$$
\| \bphi_\eps\|_{H^1(\O)} \le {\mathcal C}_{10}(\lambda) \eps^{1/2} \|\bPhi\|_{L_2(\O)},
\quad 0< \eps \le \eps_1.
\eqno(6.22)
$$

\smallskip\noindent\textbf{Proof.} The proof is similar to that of Lemma 7.4 from [PSu2].

First, we estimate the $L_2(\O;\C^n)$-norm
of the function (6.20).
By (5.21) and Proposition 1.5, we have
$$
\| {\boldsymbol{\phi}}_{\varepsilon} \|_{L_2 (\mathcal{O})} \leq
\varepsilon  \| \Lambda^{\varepsilon} S_\eps b (\mathbf{D}) \wt{\eeta}_0 \|_{L_2 (\R^d)}
\leq \varepsilon |\Omega|^{-1/2} \|\Lambda\|_{L_2(\Omega)} \| b(\D)\wt{\eeta}_0 \|_{L_2(\R^d)}.
\eqno(6.23)
$$
From (1.4) and (6.7) it follows that
$$
\| b(\D)\wt{\eeta}_0 \|_{L_2(\R^d)} \le \alpha_1^{1/2}
\|\wt{\eeta}_0\|_{H^1(\R^d)} \le \alpha_1^{1/2} C_\O \wh{C}_\lambda \|{\bPhi}\|_{L_2(\O)}.
\eqno(6.24)
$$
Now relations (6.23), (6.24), and (1.11) imply that
$$
\|\bphi_\eps \|^2_{L_2(\O)} \le \gamma_0(\lambda) \eps^2 \|\bPhi\|^2_{L_2(\O)},
\eqno(6.25)
$$
where
$\gamma_0(\lambda) = m (2r_0)^{-2} \alpha_1 \alpha_0^{-1}
\|g\|_{L_\infty} \|g^{-1}\|_{L_\infty} C_\O^2 \wh{C}^2_\lambda$.

Consider the derivatives
$$
\begin{aligned}
&\partial_j {\boldsymbol{\phi}}_{\varepsilon} =
\varepsilon (\partial_j \theta_{\varepsilon})
\Lambda^{\varepsilon} S_\eps b (\mathbf{D}) \widetilde{\eeta}_0
+ \theta_{\varepsilon}
(\partial_j \Lambda)^\eps S_\eps b (\mathbf{D}) \widetilde{\eeta}_0
\cr
&+ \varepsilon \theta_{\varepsilon} \Lambda^{\varepsilon}
 S_\eps b (\mathbf{D}) \partial_j \widetilde{\eeta}_0,\ \ j=1,\dots,d.
\end{aligned}
$$
Then
$$
\begin{aligned}
&\| \D {\boldsymbol{\phi}}_{\varepsilon}\|^2_{L_2(\O)} \leq
3 \eps^2 \int_\O |\nabla \theta_\eps |^2 |\Lambda^\eps S_\eps b(\D)\wt{\eeta}_0 |^2\,d\x
\\
&+
3 \int_\O |(\D \Lambda)^\eps|^2 |\theta_\eps S_\eps b(\D)\wt{\eeta}_0 |^2\,d\x
+3 \eps^2 \sum_{j=1}^d \int_\O |\theta_\eps|^2 |\Lambda^\eps S_\eps b(\D) D_j \wt{\eeta}_0|^2\,d\x.
\end{aligned}
\eqno(6.26)
$$
Denote the summands in the right-hand side of (6.26) by
${\mathfrak A}_1(\eps)$, ${\mathfrak A}_2(\eps)$, and ${\mathfrak A}_3(\eps)$, respectively.

The easiest task is to estimate
${\mathfrak A}_3(\eps)$. By (5.21) and Proposition 1.5,
$$
\begin{aligned}
{\mathfrak A}_3(\eps) &\leq
3 \eps^2 \sum_{j=1}^d \left\| \Lambda^{\varepsilon}
S_\eps b (\mathbf{D}) D_j \wt{\eeta}_0 \right\|^2_{L_2(\R^d)}
\cr
&\leq 3 \eps^2 |\Omega|^{-1} \|\Lambda\|^2_{L_2(\Omega)}
 \| \D b(\D) \wt{\eeta}_0\|^2_{L_2(\R^d)}.
\end{aligned}
\eqno(6.27)
$$
Similarly to (5.9), by (6.7),
$$
 \| \D b(\D) \wt{\eeta}_0\|^2_{L_2(\R^d)}\le
\alpha_1 C_\O^2 \wh{C}^2_\lambda \|\bPhi\|^2_{L_2(\O)}.
\eqno(6.28)
$$
From (6.27), (6.28), and (1.11) it follows that
$$
{\mathfrak A}_3(\eps) \leq \gamma_3(\lambda) \eps^2 \|\bPhi\|^2_{L_2(\O)},
\eqno(6.29)
$$
where $\gamma_3 (\lambda) = 3 m (2r_0)^{-2} \alpha_1 \alpha_0^{-1} \|g\|_{L_\infty} \|g^{-1}\|_{L_\infty}
C_\O^2 \wh{C}^2_\lambda$.

The first term in the right-hand side of (6.26) can be estimated by using (5.21) and Lemma 3.3.
For $0< \eps \leq \eps_1$ we have
$$
\begin{aligned}
&{\mathfrak A}_1(\eps) \leq 3 \kappa^2 \int_{(\partial \O)_\eps} |\Lambda^\eps
S_\eps b (\mathbf{D}) \widetilde{\eeta}_0 |^2\,d\x
\cr
&\leq 3 \kappa^2 \beta_* \eps |\Omega|^{-1}\|\Lambda\|^2_{L_2(\Omega)}
\| b (\mathbf{D}) \wt{\eeta}_0 \|^2_{H^1(\R^d)}.
\end{aligned}
\eqno(6.30)
$$
Note that, by (1.4) and (6.7),
$$
\| b (\mathbf{D}) \wt{\eeta}_0 \|_{H^1(\R^d)}
\le  \alpha_1^{1/2} \|\wt{\eeta}_0\|_{H^2(\R^d)}
\le \alpha_1^{1/2} C_\O \wh{C}_\lambda  \|\bPhi\|_{L_2(\O)}.
\eqno(6.31)
$$
As a result, (6.30), (6.31), and (1.11) yield
$$
{\mathfrak A}_1(\eps) \leq \gamma_1(\lambda) \eps \|\bPhi\|^2_{L_2(\O)},\quad 0< \eps \le \eps_1,
\eqno(6.32)
$$
where
$\gamma_1(\lambda) = 3 \kappa^2 \beta_* m (2r_0)^{-2} \alpha_1 \alpha_0^{-1} \|g\|_{L_\infty} \|g^{-1}\|_{L_\infty}
C_\O^2 \wh{C}^2_\lambda$.

It remains to consider the second term in the right-hand side of (6.26).
By (5.21),
$$
{\mathfrak A}_2(\eps)
\leq 3 \int_{(\partial \O)_\eps} |(\D \Lambda)^\eps|^2 |S_\eps b(\D)\wt{\eeta}_0|^2\,d\x.
$$
Applying Lemma 3.3, for any $0< \eps \leq \eps_1$ we obtain
$$
{\mathfrak A}_2(\eps) \leq 3 \beta_* \eps |\Omega|^{-1} \|\D \Lambda\|^2_{L_2(\Omega)}
\| b (\mathbf{D}) \wt{\eeta}_0 \|^2_{H^1(\R^d)}.
$$
Together with (1.10) and (6.31) this implies that
$$
{\mathfrak A}_2(\eps) \leq {\gamma}_2(\lambda) \eps \| \bPhi \|^2_{L_2(\O)},\ \ 0< \eps \leq \eps_1,
\eqno(6.33)
$$
where
${\gamma}_2(\lambda) = 3 \beta_* m \alpha_1
\alpha_0^{-1} \|g\|_{L_\infty}\|g^{-1}\|_{L_\infty} C_\O^2 \wh{C}^2_\lambda$.

Finally, from (6.26), (6.29), (6.32), and (6.33) it follows that
$$
\begin{aligned}
\| \D \bphi_\eps \|^2_{L_2(\O)} &\le
{\mathfrak A}_1(\eps) + {\mathfrak A}_2(\eps) + {\mathfrak A}_3(\eps)
\cr
&\le (\gamma_1(\lambda) + \gamma_2(\lambda) + \gamma_3(\lambda)) \eps \|\bPhi\|^2_{L_2(\O)},
\quad 0< \eps \le \eps_1.
\end{aligned}
$$
Together with (6.25) this yields (6.22) with the constant
$$
{\mathcal C}_{10}(\lambda)= \left(\gamma_0(\lambda) + \gamma_1(\lambda) +
\gamma_2(\lambda) + \gamma_3(\lambda)\right)^{1/2}. \ \ \bullet
$$

\smallskip
From (6.21) and Lemma 6.1 it follows that
$$
\|\z_\eps\|_{H^1(\O)} \le 2 {\mathcal C}_{10}(\lambda) \eps^{1/2} \|\bPhi\|_{L_2(\O)},\ \ 0< \eps \le \eps_1.
\eqno(6.34)
$$
Now (6.18) and (6.34) imply that
$$
|{\mathcal J}_\eps^{(3)}[\bPhi]| \le {\mathcal C}^{(3)}(\lambda) \eps \|\FF\|_{L_2(\O)} \| \bPhi \|_{L_2(\O)},
\quad \bPhi \in L_2(\O;\C^n),\quad 0 < \eps \le \eps_1,
\eqno(6.35)
$$
where ${\mathcal C}^{(3)}(\lambda) = 2 \wt{\mathfrak C}(\lambda) d^{1/2} {\mathcal C}_{10}(\lambda)$.

\smallskip\noindent\textbf{6.4. Completion of the proof of Theorem 4.1.}
Let us summarize the results. From (6.9), (6.13), (6.16), and (6.35) it follows that
$$
\begin{aligned}
|\wt{\mathcal I}_\eps^{(2)}[\eeta_\eps]| \le ({\mathcal C}^{(1)}(\lambda)+{\mathcal C}^{(2)}(\lambda)
+ {\mathcal C}^{(3)}(\lambda)) \eps \|\FF\|_{L_2(\O)} \| \bPhi \|_{L_2(\O)},
\cr \bPhi \in L_2(\O;\C^n),\quad 0 < \eps \le \eps_1.
\end{aligned}
$$
Together with (6.6) this yields
$$
|(\w_\eps,\bPhi)_{L_2(\O)}| \le {\mathcal C}_{11}(\lambda) \eps \|\FF\|_{L_2(\O)} \|\bPhi\|_{L_2(\O)},
\quad \bPhi \in L_2(\O;\C^n),\quad 0 < \eps \le \eps_1,
$$
where ${\mathcal C}_{11}(\lambda) = ({\mathcal C}_6(\lambda) + {\mathcal C}_7(\lambda))c_\lambda^{-1}
+{\mathcal C}^{(1)}(\lambda)+{\mathcal C}^{(2)}(\lambda)
+ {\mathcal C}^{(3)}(\lambda)$.
The desired estimate follows:
$$
\|\w_\eps\|_{L_2(\O)} \le {\mathcal C}_{11}(\lambda) \eps \|\FF\|_{L_2(\O)},
\quad 0 < \eps \le \eps_1.
\eqno(6.36)
$$

Finally, relations (4.27) and (6.36) imply (4.1) with the constant
${\mathcal C}_0(\lambda)= {\mathcal C}_5(\lambda) + C_1(\lambda) {\mathcal C}_3(\lambda)
+ {\mathcal C}_{11}(\lambda)$. Theorem 4.1 is proved. $\ \bullet$

\section*{\S 7. Homogenization of the Neumann problem in the case $\Lambda \in L_\infty$}

\noindent\textbf{7.1.} As for the homogenization problem in $\R^d$ (see Subsection 1.7),
the smoothing operator $S_\eps$ in the corrector (4.5)
can be removed under additional assumptions on the matrix $\Lambda(\x)$.

Assume that Condition 1.10 is satisfied, i.~e., $\Lambda \in L_\infty$. Denote
$$
K^0_{N,\lambda}(\eps) = [\Lambda^\eps] b(\D) (\A_{N}^0 + \lambda I)^{-1}.
\eqno(7.1)
$$
By (2.34), the operator $b(\D) (\A_{N}^0 + \lambda I)^{-1}$ is a continuous mapping of
$L_2(\O;\C^n)$ into $H^1(\O;\C^m)$. Under Condition 1.10 the operator $[\Lambda^\eps]$
is continuous from $H^1(\O;\C^m)$ to
$H^1(\O;\C^n)$ which easily follows from Proposition 1.11 and existence of a linear continuous
extension operator from $H^1(\O;\C^m)$ to $H^1(\R^d;\C^m)$.
Hence, the corrector (7.1) is a continuous mapping of
$L_2(\O;\C^n)$ to $H^1(\O;\C^n)$.

Instead of (4.6), we consider another approximation
of the solution $\u_\eps$ of the problem (2.8):
$$
\check{\v}_\eps = (\A_{N}^0 + \lambda I)^{-1}\FF + \eps K^0_{N,\lambda}(\eps) \FF =
\u_0 + \eps \Lambda^\eps b(\D) \u_0.
\eqno(7.2)
$$

\smallskip\noindent\textbf{Theorem 7.1.} \textit{Suppose that the assumptions
of Theorem} 4.1
\textit{and Condition} 1.10 \textit{are satisfied. Let $\check{\v}_\eps$ be defined by} (7.2).
\textit{Then there exists a number $\eps_1 \in (0,1]$ depending on the domain $\O$
and the lattice $\Gamma$ such that}
$$
\| \u_\eps - \check{\v}_\eps\|_{H^1(\O;\C^n)} \le \check{\mathcal C}(\lambda) \eps^{1/2}
\| \FF \|_{L_2(\O;\C^n)},\quad 0 < \eps \le \eps_1,
\eqno(7.3)
$$
\textit{or, in operator terms,}
$$
\| (\mathcal{A}_{N,\varepsilon}+\lambda I)^{-1} - (\mathcal{A}^0_N + \lambda I)^{-1} -
\varepsilon K^0_{N,\lambda} (\varepsilon) \|_{L_2 (\mathcal{O}; \mathbb{C}^{n}) \to
H^1 (\mathcal{O}; \mathbb{C}^{n})} \leq \check{\mathcal C}(\lambda) {\varepsilon}^{1/2}.
$$
\textit{For the flux} $\p_\eps= g^\eps b(\D)\u_\eps$ \textit{we have}
$$
\| \p_\eps - \wt{g}^\eps b(\D) {\u}_0 \|_{L_2(\O;\C^m)}
\leq \check{\mathcal C}'(\lambda) \eps^{1/2} \| \mathbf{F} \|_{L_2 (\mathcal{O}; \mathbb{C}^{n})},\ \ 0< \eps \leq \eps_1,
\eqno(7.4)
$$
\textit{where} $\wt{g}(\x)$ \textit{is the matrix} (1.9).
\textit{The constants $\check{\mathcal C}(\lambda)$ and $\check{\mathcal C}'(\lambda)$ depend on}
$m$, $d$, $\alpha_0$, $\alpha_1$, $\| g \|_{L_\infty}$,
$\| g^{-1} \|_{L_\infty}$, $\lambda$, \textit{the parameters of the lattice} $\Gamma$,
\textit{the constants ${\mathcal C}_1$ and ${\mathcal C}_2$ from the inequality} (2.2),
\textit{the norm $\|\Lambda\|_{L_\infty}$, and the domain} $\mathcal{O}$.

\smallskip
\noindent\textbf{Proof.} According to (4.5), (4.6), and (7.2), in the domain $\O$ we have
$$
\check{\v}_\eps - \v_\eps = \eps \Lambda^\eps (I - S_\eps) b(\D)\wt{\u}_0.
\eqno(7.5)
$$
Let us estimate the norm of the function (7.5) in $H^1(\O;\C^n)$.
Since the right-hand side of (7.5) is defined in the whole $\R^d$,
it suffices to estimate its norm in $H^1(\R^d;\C^n)$.
We apply the proof of Theorem 1.12. By (1.52) and (1.58),
$$
\eps \|\Lambda^\eps (I - S_\eps) b(\D)\wt{\u}_0 \|_{L_2(\R^d)}\le
2 \eps \|\Lambda\|_{L_\infty}\| b(\D)\wt{\u}_0 \|_{L_2(\R^d)},
\eqno(7.6)
$$
$$
\begin{aligned}
&\sum_{l=1}^d \| \partial_l (\eps \Lambda^\eps (I - S_\eps) b(\D)\wt{\u}_0) \|^2_{L_2(\R^d)}
\cr
&\le
2 \eps^2 \left( 4 (1+\beta_2) \|\Lambda\|_{L_\infty}^2 + \beta_1 r_1^2 \right)
\|\D b(\D)\wt{\u}_0 \|^2_{L_2(\R^d)}.
\end{aligned}
\eqno(7.7)
$$
Together with (5.9), (5.30), and (7.5) this implies that
$$
\| \check{\v}_\eps - \v_\eps\|_{H^1(\O)}
= \eps \|\Lambda^\eps (I - S_\eps) b(\D)\wt{\u}_0 \|_{H^1(\R^d)}
\le {\mathcal C}_{12}(\lambda) \eps \| {\FF}\|_{L_2(\O)},\quad \eps >0,
\eqno(7.8)
$$
where
${\mathcal C}_{12}(\lambda) = \alpha_1^{1/2} C_\O \wh{C}_\lambda
\left( 2 (3+2\beta_2)^{1/2} \|\Lambda\|_{L_\infty} + (2\beta_1)^{1/2} r_1 \right)$.

As a result, combining (4.9) and (7.8), we arrive at the estimate (7.3) with the constant
$\check{\mathcal C}(\lambda) = {\mathcal C}(\lambda) + {\mathcal C}_{12}(\lambda)$.

It remains to check (7.4). From (7.3), taking (1.2) and (1.5) into account, we obtain
$$
\| \p_\eps - g^\eps b(\D) \check{\v}_\eps \|_{L_2(\O)}
\le \eps^{1/2}  d^{1/2}\alpha^{1/2}_1 \|g\|_{L_\infty}\check{\mathcal C}(\lambda) \|\FF\|_{L_2(\O)},
\quad 0< \eps \le \eps_1.
\eqno(7.9)
$$
By (1.2), (1.9), and (7.2),
$$
g^\eps b(\D) \check{\v}_\eps = \wt{g}^\eps b(\D) \u_0 +
\eps \sum_{l=1}^d g^\eps b_l \Lambda^\eps b(\D) D_l \u_0.
\eqno(7.10)
$$
Using Condition 1.10 and relations (1.2), (1.5), and (2.33),
we estimate the norm of the second term:
$$
\begin{aligned}
&\eps \left\| \sum_{l=1}^d g^\eps b_l \Lambda^\eps b(\D) D_l \u_0 \right\|_{L_2(\O)}\le
\eps \alpha_1 \|g\|_{L_\infty} \|\Lambda \|_{L_\infty} \sum_{l,j=1}^d \|D_j D_l \u_0\|_{L_2(\O)}
\cr
&\le  \eps d \alpha_1 \|g\|_{L_\infty} \|\Lambda \|_{L_\infty} \wh{C}_\lambda \|\FF\|_{L_2(\O)}.
\end{aligned}
\eqno(7.11)
$$
Now (7.9)--(7.11) imply estimate (7.4) with the constant
$\check{\mathcal C}'(\lambda)
= d^{1/2}\alpha^{1/2}_1 \|g\|_{L_\infty}\check{\mathcal C}(\lambda)+
d \alpha_1 \|g\|_{L_\infty} \|\Lambda \|_{L_\infty} \wh{C}_\lambda$. $\ \bullet$

\smallskip
Recall that some sufficient conditions
which ensure Condition 1.10 are given in Proposition 1.13.

\smallskip\noindent\textbf{7.2. Special case.}
We distinguish the special case where $g^0 = \underline{g}$.
Then $\wt{g}(\x)= g^0 = \underline{g}$ and, by Proposition 1.13, Condition 1.10 is satisfied.
Applying the statement of Theorem 7.1 concerning the flux, we arrive at the following statement.

\smallskip\noindent\textbf{Proposition 7.2.}
\textit{Let} $\mathbf{u}_{\varepsilon}$ \textit{be the solution of the problem} (2.8),
\textit{and let} $\mathbf{u}_0$ \textit{be the solution of the problem} (2.28).
\textit{Let} $\p_\eps =g^\eps b(\D)\u_\eps$.
\textit{Assume that}
$g^0 = \underline{g}$, \textit{i.~e., relations} (1.14)
\textit{are satisfied. Then we have}
$$
\| \p_\eps - {g}^0 b(\D) {\u}_0 \|_{L_2(\O;\C^m)}
\leq \check{\mathcal C}'(\lambda) \eps^{1/2} \| \mathbf{F} \|_{L_2 (\mathcal{O}; \mathbb{C}^{n})},\ \ 0< \eps \leq \eps_1.
$$

\section*{\S 8. Approximation of solutions in a strictly interior subdomain}

\noindent\textbf{8.1.}
Let $\O'$ be a strictly interior subdomain of the domain $\O$.
Applying Theorem 4.1 and the results for the problem in $\R^d$
(Theorems 1.7 and 1.8), it is easy to obtain a sharp order error estimate
for approximation of the solution of the problem (2.8)
in $H^1(\O';\C^n)$.

\smallskip\noindent\textbf{Theorem 8.1.} \textit{Suppose that the assumptions of Theorem} 4.2
\textit{are satisfied. Let $\O'$ be a strictly interior subdomain of the domain $\O$. Let}
$\delta:=\text{dist}\, \{\O';\partial \O\}$.
\textit{Then there exists a number $\eps_1 \in (0,1]$ depending on the domain $\O$
and the lattice $\Gamma$ such that}
$$
\| \u_\eps - {\v}_\eps\|_{H^1(\O';\C^n)} \le {\mathfrak C}_{\delta}(\lambda) \eps
\| \FF \|_{L_2(\O;\C^n)},\quad 0 < \eps \le \eps_1,
\eqno(8.1)
$$
\textit{or, in operator terms,}
$$
\| (\A_{N,\eps} +\lambda I)^{-1} - (\A_N^0 + \lambda I)^{-1} - \eps K_{N,\lambda}(\eps)
\|_{L_2(\O;\C^n) \to H^1(\O';\C^n)} \le
{\mathfrak C}_{\delta}(\lambda) \eps.
\eqno(8.2)
$$
\textit{For the flux} $\p_\eps= g^\eps b(\D)\u_\eps$ \textit{we have}
$$
\| \p_\eps - \wt{g}^\eps S_\eps b(\D) \wt{\u}_0 \|_{L_2(\O';\C^m)}
\leq {\mathfrak C}'_{\delta}(\lambda) \eps
\| \mathbf{F} \|_{L_2 (\mathcal{O}; \mathbb{C}^{n})},\ \ 0< \eps \leq \eps_1,
\eqno(8.3)
$$
\textit{where} $\wt{g}(\x)$ \textit{is the matrix} (1.9).
\textit{The constants in estimates are given by} ${\mathfrak C}_{\delta}(\lambda)=
{\mathfrak C}_1(\lambda) \delta^{-1} + {\mathfrak C}_2(\lambda)$,
${\mathfrak C}_{\delta}'(\lambda)= {\mathfrak C}'_1(\lambda) \delta^{-1} + {\mathfrak C}'_2(\lambda)$,
 \textit{where ${\mathfrak C}_1(\lambda)$, ${\mathfrak C}_2(\lambda)$,
 ${\mathfrak C}'_1(\lambda)$, and ${\mathfrak C}'_2(\lambda)$ depend on}
$m$, $d$, $\alpha_0$, $\alpha_1$, $\| g \|_{L_\infty}$,
$\| g^{-1} \|_{L_\infty}$, $\lambda$, \textit{the parameters of the lattice} $\Gamma$,
\textit{the constants ${\mathcal C}_1$ and ${\mathcal C}_2$ from} (2.2),
\textit{and the domain} $\mathcal{O}$.

\smallskip\noindent\textbf{Proof.}
We fix a smooth cut-off function $\zeta(\x)$ such that
$$
\begin{aligned}
&\zeta \in C_0^\infty(\O),\ \ 0\le \zeta(\x) \le 1,\cr
&\zeta(\x)=1\  \text{for}\ \x \in \O';\ \ |\nabla \zeta(\x)| \le \kappa' \delta^{-1}.
\end{aligned}
\eqno(8.4)
$$
Here the constant $\kappa'$ depends only on the domain $\O$.
Let $\u_\eps$ be the solution of the problem (2.8), and let $\wt{\u}_\eps$ be the
solution of the equation (4.14).
Then $(\A_\eps + \lambda I)(\u_\eps - \wt{\u}_\eps)=0$ in the domain $\O$.
Hence,
$$
\int_{\O} \left( \langle g^\eps b(\D)(\u_\eps - \wt{\u}_\eps), b(\D)\eeta \rangle
+ \lambda \langle \u_\eps - \wt{\u}_\eps, \eeta \rangle\right)\,d\x =0,
\quad \forall \eeta \in H^1_0(\O;\C^n).
\eqno(8.5)
$$
Substitute $\eeta= \zeta^2 (\u_\eps - \wt{\u}_\eps)$ in (8.5).
By (1.2),
$$
b(\D)\left(\zeta^2 (\u_\eps - \wt{\u}_\eps)\right) =
\zeta b(\D)\left(\zeta (\u_\eps - \wt{\u}_\eps)\right) +
\sum_{l=1}^d b_l (D_l\zeta) \zeta (\u_\eps - \wt{\u}_\eps),
\eqno(8.6)
$$
$$
\zeta b(\D) (\u_\eps - \wt{\u}_\eps) =
 b(\D)\left(\zeta (\u_\eps - \wt{\u}_\eps)\right)
- \sum_{l=1}^d b_l (D_l\zeta) (\u_\eps - \wt{\u}_\eps).
\eqno(8.7)
$$
From (8.5)--(8.7) it follows that
$$
\begin{aligned}
J(\eps)&:= \int_{\O} \left( \langle g^\eps b(\D)\left(\zeta(\u_\eps - \wt{\u}_\eps)\right),
b(\D) \left(\zeta(\u_\eps - \wt{\u}_\eps)\right)\rangle
+ \lambda \zeta^2 |\u_\eps - \wt{\u}_\eps|^2 \right)\,d\x
\cr
&= J_1(\eps) +J_2(\eps),
\end{aligned}
\eqno(8.8)
$$
where
$$
\begin{aligned}
&J_1(\eps) = -
 \int_{\O}  \langle g^\eps  b(\D)\left(\zeta(\u_\eps - \wt{\u}_\eps)\right),
\sum_{l=1}^d b_l (D_l \zeta) (\u_\eps - \wt{\u}_\eps)\rangle \,d\x
\cr
&+ \int_{\O}  \langle g^\eps \sum_{j=1}^d b_j (D_j \zeta) (\u_\eps - \wt{\u}_\eps),
\sum_{l=1}^d b_l (D_l \zeta) (\u_\eps - \wt{\u}_\eps)\rangle \,d\x,
\end{aligned}
\eqno(8.9)
$$
$$
J_2(\eps) =
\int_{\O} \langle g^\eps \sum_{l=1}^d b_l (D_l \zeta)(\u_\eps - \wt{\u}_\eps),
b(\D) \left( \zeta (\u_\eps - \wt{\u}_\eps)\right) \rangle \,d\x.
\eqno(8.10)
$$
The term (8.10) is estimated with the help of (1.5) and (8.4):
$$
\begin{aligned}
&|J_2(\eps)| \le \|(g^\eps)^{1/2}b(\D) \left(\zeta (\u_\eps - \wt{\u}_\eps) \right) \|_{L_2(\O)}
\|(g^\eps)^{1/2} \sum_{l=1}^d b_l (D_l \zeta)(\u_\eps - \wt{\u}_\eps)\|_{L_2(\O)}
\cr
&\le \frac{1}{4} J(\eps) + \|g\|_{L_\infty} \alpha_1 d (\kappa')^2 \delta^{-2} \|\u_\eps - \wt{\u}_\eps\|^2_{L_2(\O)}.
\end{aligned}
\eqno(8.11)
$$
Similarly for the term (8.9) we obtain:
$$
|J_1(\eps)| \le
\frac{1}{4} J(\eps) + 2 \|g\|_{L_\infty} \alpha_1 d (\kappa')^2 \delta^{-2} \|\u_\eps - \wt{\u}_\eps\|^2_{L_2(\O)}.
\eqno(8.12)
$$
Now (8.8), (8.11), and (8.12) imply that
$$
J(\eps) \le 6 \|g\|_{L_\infty} \alpha_1 d (\kappa')^2 \delta^{-2} \|\u_\eps - \wt{\u}_\eps\|^2_{L_2(\O)}.
\eqno(8.13)
$$

Extending the function
$\bvarphi_\eps:= \zeta (\u_\eps - \wt{\u}_\eps)$ by zero to $\R^d \setminus \O$, note that
$J(\eps) = a_\eps[\bvarphi_\eps,\bvarphi_\eps] + \lambda \|\bvarphi_\eps\|^2_{L_2(\R^d)}$
and apply the lower inequality (1.6). We obtain
$$
J(\eps) \ge \alpha_0 \|g^{-1}\|^{-1}_{L_\infty} \|\D \bvarphi_\eps \|_{L_2(\R^d)}^2 + \lambda \|\bvarphi_\eps\|^2_{L_2(\R^d)},
$$
whence
$$
\| \zeta (\u_\eps - \wt{\u}_\eps)\|^2_{H^1(\O)} \le
\max \{ \alpha_0^{-1}\|g^{-1}\|_{L_\infty},\lambda^{-1}\} J(\eps).
\eqno(8.14)
$$

By (4.26) and (6.36),
$$
\|\u_\eps - \wt{\u}_\eps\|_{L_2(\O)} \le ({\mathcal C}_5(\lambda)+ {\mathcal C}_{11}(\lambda))\eps \| \FF\|_{L_2(\O)},\ \ 0< \eps \le \eps_1.
\eqno(8.15)
$$
Combining (8.13)--(8.15) and (8.4), we arrive at
$$
\| \u_\eps - \wt{\u}_\eps\|_{H^1(\O')} \le
\| \zeta (\u_\eps - \wt{\u}_\eps)\|_{H^1(\O)} \le {\mathfrak C}_1(\lambda) \delta^{-1}\eps \|\FF\|_{L_2(\O)},
\ \ 0< \eps \le \eps_1,
\eqno(8.16)
$$
where
$$
{\mathfrak C}_1(\lambda)= \max\{\alpha_0^{-1/2}\|g^{-1}\|^{1/2}_{L_\infty},\lambda^{-1/2}\} \|g\|_{L_\infty}^{1/2} (6 \alpha_1 d)^{1/2} \kappa'
({\mathcal C}_5(\lambda)+ {\mathcal C}_{11}(\lambda)).
\eqno(8.17)
$$

Taking (4.8), (4.13), and (4.16) into account,
we check that
$$
\|\wt{\u}_\eps - \v_\eps\|_{H^1(\O)} \le C_2(\lambda) {\mathcal C}_3(\lambda) \eps \| \FF \|_{L_2(\O)},\ \ \eps>0.
\eqno(8.18)
$$
Now relations (8.16) and (8.18) imply (8.1) with the constant
${\mathfrak C}_\delta(\lambda) = {\mathfrak C}_1(\lambda) \delta^{-1} + {\mathfrak C}_2(\lambda)$,
where ${\mathfrak C}_1(\lambda)$ is defined by (8.17), and
${\mathfrak C}_2(\lambda)= C_2(\lambda) {\mathcal C}_3(\lambda)$.

It remains to check (8.3). From (8.16) and (1.2), (1.5) it follows that
$$
\| \p_\eps - g^\eps b(\D) \wt{\u}_\eps \|_{L_2(\O')} \le
\|g\|_{L_\infty} \alpha_1^{1/2} d^{1/2}{\mathfrak C}_1(\lambda) \delta^{-1} \eps \|\FF\|_{L_2(\O)},\ \ 0< \eps \le \eps_1.
\eqno(8.19)
$$
By (4.13) and (4.17), we obtain
$$
\|g^\eps b(\D) \wt{\u}_\eps - \wt{g}^\eps S_\eps b(\D) \wt{\u}_0\|_{L_2(\O)}
\le C_5(\lambda) {\mathcal C}_3(\lambda) \eps \|\FF\|_{L_2(\O)}, \ \ \eps>0.
\eqno(8.20)
$$
Now relations (8.19) and (8.20) imply (8.3) with the constant
${\mathfrak C}'_\delta(\lambda) = {\mathfrak C}'_1(\lambda) \delta^{-1} + {\mathfrak C}'_2(\lambda)$,
where
$$
{\mathfrak C}'_1(\lambda)= \|g\|_{L_\infty} \alpha_1^{1/2} d^{1/2}{\mathfrak C}_1(\lambda),\quad
{\mathfrak C}'_2(\lambda) = C_5(\lambda) {\mathcal C}_3(\lambda). \ \ \bullet
\eqno(8.21)
$$

\noindent\textbf{8.2. The case where $\Lambda \in L_\infty$.}
Similarly, under Condition 1.10,
Theorems 4.1 and 1.12 imply the following result.

\smallskip\noindent\textbf{Theorem 8.2.} \textit{Suppose that the assumptions of Theorem} 7.1
\textit{are satisfied. Let $\O'$ be a strictly interior subdomain of the domain $\O$. Denote}
$\delta:=\text{dist}\, \{\O';\partial \O\}$.
\textit{Then there exists a number $\eps_1 \in (0,1]$ depending on the domain $\O$ and the
lattice $\Gamma$ such that}
$$
\| \u_\eps - \check{\v}_\eps\|_{H^1(\O';\C^n)} \le \check{\mathfrak C}_{\delta}(\lambda) \eps
\| \FF \|_{L_2(\O;\C^n)},\quad 0 < \eps \le \eps_1,
\eqno(8.22)
$$
\textit{or, in operator terms},
$$
\| (\A_{N,\eps} +\lambda I)^{-1} - (\A_N^0 + \lambda I)^{-1} - \eps K^0_{N,\lambda}(\eps)
\|_{L_2(\O;\C^n) \to H^1(\O';\C^n)} \le
\check{\mathfrak C}_{\delta}(\lambda) \eps.
$$
\textit{For the flux} $\p_\eps= g^\eps b(\D)\u_\eps$ \textit{we have}
$$
\| \p_\eps - \wt{g}^\eps b(\D) {\u}_0 \|_{L_2(\O';\C^m)}
\leq \check{\mathfrak C}'_{\delta}(\lambda) \eps
\| \mathbf{F} \|_{L_2 (\mathcal{O}; \mathbb{C}^{n})},\ \ 0< \eps \leq \eps_1,
\eqno(8.23)
$$
\textit{where} $\wt{g}(\x)$ \textit{is the matrix} (1.9).
\textit{The constants in estimates are given by} $\check{\mathfrak C}_{\delta}(\lambda)=
{\mathfrak C}_1(\lambda) \delta^{-1} + \check{\mathfrak C}_2(\lambda)$,
$\check{\mathfrak C}_{\delta}'(\lambda)= {\mathfrak C}'_1(\lambda) \delta^{-1} + \check{\mathfrak C}'_2(\lambda)$,
 \textit{where ${\mathfrak C}_1(\lambda)$, $\check{\mathfrak C}_2(\lambda)$,
 ${\mathfrak C}'_1(\lambda)$, and $\check{\mathfrak C}'_2(\lambda)$ depend on}
$m$, $d$, $\alpha_0$, $\alpha_1$, $\| g \|_{L_\infty}$,
$\| g^{-1} \|_{L_\infty}$, $\lambda$, \textit{the parameters of the lattice} $\Gamma$,
\textit{the constants ${\mathcal C}_1$ and ${\mathcal C}_2$ from} (2.2),
\textit{the norm $\|\Lambda\|_{L_\infty}$, and the domain} $\mathcal{O}$.

\smallskip\noindent\textbf{Proof.}
The inequality (8.16) remains true. Instead of (8.18), now we use
the next estimate which follows from Theorem 1.12
and (4.13):
$$
\| \wt{\u}_\eps - \check{\v}_\eps \|_{H^1(\O)} \le C_7(\lambda) \eps \|\wt{\FF}\|_{L_2(\R^d)}
\le C_7(\lambda) {\mathcal C}_3(\lambda) \eps \| {\FF}\|_{L_2(\O)},\ \ \eps >0.
\eqno(8.24)
$$
Then (8.22) is a direct consequence of (8.16) and (8.24), herewith,
$\check{\mathfrak C}_{\delta}(\lambda)=
{\mathfrak C}_1(\lambda) \delta^{-1} + \check{\mathfrak C}_2(\lambda)$,
where ${\mathfrak C}_1(\lambda)$ is defined by (8.17) and $\check{\mathfrak C}_2(\lambda)=
C_7(\lambda) {\mathcal C}_3(\lambda)$.

In order to prove (8.23), we use (8.19) and the next estimate which follows from
Theorem 1.12 and (4.13):
$$
\begin{aligned}
 \| g^\eps b(\D) \wt{\u}_\eps - \wt{g}^\eps b(\D) {\u}_0 \|_{L_2(\O)}
 &\le C_8(\lambda) \eps \|\wt{\FF} \|_{L_2(\R^d)}
 \cr
 &\le C_8(\lambda) {\mathcal C}_3(\lambda) \eps \| {\FF} \|_{L_2(\O)},\ \ \eps>0.
\end{aligned}
$$
As a result, we obtain (8.23) with the constant
$\check{\mathfrak C}_{\delta}'(\lambda)=
{\mathfrak C}'_1(\lambda) \delta^{-1} + \check{\mathfrak C}'_2(\lambda)$,
where ${\mathfrak C}'_1(\lambda)$ is defined by (8.21) and
$\check{\mathfrak C}'_2(\lambda)= C_8(\lambda) {\mathcal C}_3(\lambda)$. $\ \bullet$

\smallskip\noindent\textbf{Remark 8.3.}
Theorems 8.1 and 8.2 can be applied also in the case
where $\delta$ depends on $\eps$ and a strictly interior subdomain
$\O_\eps'$ approaches the boundary $\partial \O$.
For instance, if $\delta(\eps)=O(\eps^\alpha)$, then the error estimate in (8.1) or in (8.22)
is of order $O(\eps^{1-\alpha})$. Clearly, in the case where $\alpha<1/2$
the error estimate in a strictly interior subdomain
$\O'_\eps$ is better in order than in the whole domain $\O$.
If $\alpha \ge 1/2$, then Theorems 4.2 and 7.1 give better results (in the whole domain $\O$)
and it does not make sense to apply Theorems 8.1 and 8.2.

\section*{\S 9. Homogenization of the Neumann problem in the case $\lambda =0$}

\noindent\textbf{9.1. Statement of the problem.}
For applications, it is interesting to study the Neumann problem
for the equation $\A_{N,\eps} \u_\eps =\FF$, i.~e., in
the case where $\lambda =0$. In this case, the right-hand side of the equation
must be subject to the solvability conditions, and the solution
is subject to the additional condition of orthogonality to the kernel of the operator.
In the present section, we deduce the results for this case
from the results of \S 4 and \S 7.

We denote
$$
Z = \{\z \in H^1(\O;\C^n):\ b(\D)\z =0\}.
\eqno(9.1)
$$
Then $Z$ is a (closed) subspace of $H^1(\O;\C^n)$.
Note that for sure $Z$ contains the $n$-dimensional subspace
$\{\u \in H^1(\O;\C^n):\ \u(\x)=\c \in \C^n\}$ consisting of constant vector-valued functions.
From (2.2) it follows that
$$
\| \D \z\|_{L_2(\O)}^2 \le {\mathcal C}_1^{-1} {\mathcal C}_2 \|\z\|^2_{L_2(\O)},\quad \z \in Z.
\eqno(9.2)
$$
Due to compactness of the embedding $H^1(\O;\C^n) \subset L_2(\O;\C^n)$, inequality of the form (9.2)
can be satisfied only on a finitedimensional subspace in $H^1(\O;\C^n)$.
Consequently, $Z$ is finitedimensional.

Obviously, the finitedimensional space (9.1) is also a (closed) subspace of $L_2(\O;\C^n)$.
We denote ${\mathcal H}(\O):= L_2(\O;\C^n) \ominus Z$,
$H^1_\perp(\O;\C^n) = H^1(\O;\C^n) \cap {\mathcal H}(\O)$. In other words,
$$
H^1_\perp(\O;\C^n)=
\{\u \in H^1(\O;\C^n):\  (\u,\z)_{L_2(\O)}=0, \  \forall \z \in Z \}.
\eqno(9.3)
$$
Clearly, $H^1_\perp(\O;\C^n)$ is a (closed) subspace of $H^1(\O;\C^n)$.

\smallskip\noindent\textbf{Proposition 9.1.} \textit{The form $\|b(\D)\u\|_{L_2(\O)}$
defines a norm in the space} (9.3) \textit{equivalent to the standard $H^1$-norm.}

\smallskip\noindent\textbf{Proof.}
The estimate $\|b(\D)\u\|^2_{L_2(\O)} \le \alpha_1 d \|\u\|^2_{H^1(\O)}$
obviously follows from (1.2) and (1.5) and is valid for all functions $\u \in H^1(\O;\C^n)$.

Let us prove the opposite estimate: there exists a constant $\wt{\mathcal C}_1 >0$ such that
$$
\|\u\|^2_{H^1(\O)} \le \wt{\mathcal C}_1 \|b(\D) \u \|^2_{L_2(\O)},\quad \u \in H^1_\perp(\O;\C^n).
\eqno(9.4)
$$
Assume the opposite. Suppose that for any $k\in \N$ there exists a function
$\u_k \in H^1_\perp(\O;\C^n)$ such that
$$
\| \u_k \|^2_{H^1(\O)} > k \| b(\D) \u_k \|^2_{L_2(\O)}.
$$
Then for the sequence $\v_k = \|\u_k\|^{-1}_{H^1(\O)}\u_k$ we have
$\|\v_k\|_{H^1(\O)}=1$ and $\|b(\D) \v_k\|_{L_2(\O)} < 1/\sqrt{k} \to 0$, as $k\to \infty$.
Due to compactness of the embedding $H^1(\O;\C^n) \subset L_2(\O;\C^n)$
there exists a subsequence $\{\v_{k_j}\}$ converging in
$L_2(\O;\C^n)$. Applying the inequality (2.2) with $\u = \v_{k_j} - \v_{k_l}$
and using convergence of $\v_{k_j}$ in $L_2(\O)$ and convergence of $b(\D)\v_{k_j}$ in $L_2(\O)$,
we conclude that $\{\v_{k_j}\}$ is a fundamental sequence
in $H^1(\O;\C^n)$. Consequently, it tends to some function $\v_* \in H^1_\perp(\O;\C^n)$.
Then $\|\v_*\|_{H^1(\O)}=1$, and, from the other side, $b(\D)\v_*=0$, i.~e., $\v_* \in Z$.
Since $Z \cap H^1_\perp(\O;\C^n) = \{0\}$, we arrive at a contradiction. $\ \bullet$

\smallskip
Let $\A_{N,\eps}$ be the operator corresponding to the form (2.3).
Obviously, $\Ker \A_{N,\eps}=Z$.
Therefore, the orthogonal decomposition \hbox{$L_2(\O;\C^n)= Z \oplus {\mathcal H}(\O)$}
reduces the operator $\A_{N,\eps}$. Denote by $\B_{N,\eps}$ the selfadjoint operator in
${\mathcal H}(\O)$ which is the part of $\A_{N,\eps}$ in ${\mathcal H}(\O)$.
Equivalently, the operator $\B_{N,\eps}$ is the selfadjoint operator in
${\mathcal H}(\O)$ generated by the quadratic form
$$
b_{N,\eps}[\u,\u] = \int_\O \langle g^\eps(\x) b(\D)\u, b(\D)\u \rangle d\x,
\quad \u \in H^1_\perp(\O;\C^n).
$$
This form is closed and, by (9.4), positive definite. We have
$$
\|g^{-1}\|^{-1}_{L_\infty} \wt{\mathcal C}_1^{-1}
\|\u\|^2_{H^1(\O)} \le b_{N,\eps}[\u,\u] \le \|g\|_{L_\infty} d \alpha_1 \|\u\|^2_{H^1(\O)},
\quad \u \in H^1_\perp(\O;\C^n).
\eqno(9.5)
$$

\textit{Our goal in this section} is to approximate the inverse operator $\B_{N,\eps}^{-1}$
for small $\eps$. In terms of solutions, we study the behavior of the generalized solution
$\u_\eps \in H^1_\perp(\O;\C^n)$ of the Neumann problem
$$
\begin{aligned}
&b(\D)^* g^\eps(\x) b(\D) \u_\eps(\x) = \FF(\x), \ \x \in \O;
\cr
&\partial_{\bnu}^\eps \u_\eps \vert_{\partial \O}=0;\ \ (\u_\eps,\z)_{L_2(\O)}=0,\ \forall \z \in Z,
\end{aligned}
\eqno(9.6)
$$
where $\FF \in {\mathcal H}(\O)$. Then $\u_\eps = \B_{N,\eps}^{-1} \FF$.

By definition, we say that
$\u_\eps \in H^1_\perp(\O;\C^n)$ is the generalized solution of the problem (9.6)
if it satisfies the identity
$$
\int_\O \langle g^\eps b(\D) \u_\eps, b(\D)\eeta \rangle \,d\x
= \int_\O \langle  \FF, \eeta \rangle \,d\x,\quad \forall \eeta \in H^1_\perp(\O;\C^n).
\eqno(9.7)
$$
By (9.5), the form $b_{N,\eps}[\u,\v]$ can be viewed as an inner product in the space
$H^1_\perp(\O;\C^n)$. The right-hand side of (9.7) is an antilinear continuous functional of
$\eeta \in H^1_\perp(\O;\C^n)$. By the Riss theorem, the solution $\u_\eps$ exists and is unique.
It satisfies the estimate
$$
\|\u_\eps\|_{H^1(\O)} \le \|g^{-1}\|_{L_\infty} \wt{\mathcal C}_1 \|\FF\|_{L_2(\O)}.
\eqno(9.8)
$$

\smallskip\noindent\textbf{9.2. The homogenized problem.}
Let $\A^0_{N}$ be the operator generated by the form (2.27). Obviously, $\Ker \A^0_{N}=Z$.
The decomposition \hbox{$L_2(\O;\C^n)= Z \oplus {\mathcal H}(\O)$}
reduces the operator $\A_{N}^0$. Denote by $\B^0_{N}$ the selfadjoint operator in
${\mathcal H}(\O)$ which is the part of $\A^0_{N}$ in ${\mathcal H}(\O)$.
In other words, $\B^0_{N}$ is the selfadjoint operator in
${\mathcal H}(\O)$ generated by the quadratic form
$$
b^0_{N}[\u,\u] = \int_\O \langle g^0 b(\D)\u, b(\D)\u \rangle d\x,
\quad \u \in H^1_\perp(\O;\C^n).
\eqno(9.9)
$$
By (1.15), the form (9.9) satisfies twosided estimates of the form (9.5)
with the middle part replaced by the form (9.9); the constants remain the same.
Then the operator $\B_N^0$ is positive definite:
$$
\B_N^0 \ge \|g^{-1}\|^{-1}_{L_\infty} \wt{\mathcal C}_1^{-1} I_{{\mathcal H}(\O)}.
\eqno(9.10)
$$

Let $\u_0 \in H^1_\perp(\O;\C^n)$ be the generalized solution of the Neumann problem
$$
\begin{aligned}
&b(\D)^* g^0 b(\D) \u_0(\x) = \FF(\x),\ \ \x \in \O;
\cr
&\partial_{\bnu}^0 \u_0\vert_{\partial \O}=0;\ \ (\u_0,\z)_{L_2(\O)}=0,\ \forall \z \in Z,
\end{aligned}
\eqno(9.11)
$$
where $\FF \in {\mathcal H}(\O)$. Then $\u_0 = (\B_N^0)^{-1}\FF$.

The solution $\u_0$ is defined similarly to (9.7). It exists, is unique, and
satisfies the estimate
$$
\|\u_0\|_{H^1(\O)} \le \|g^{-1}\|_{L_\infty} \wt{\mathcal C}_1 \|\FF\|_{L_2(\O)}.
\eqno(9.12)
$$
By the theorem about regularity of solutions, we have $\u_0 \in H^2(\O;\C^n)$ and
$$
\| \u_0 \|_{H^2(\O)} \le \wh{C} \|\FF\|_{L_2(\O)}.
\eqno(9.13)
$$
Note that the estimate (9.13) can be derived from (2.33)
and the fact that $\u_0$ is the solution of the problem (2.28)
with the right-hand side $\FF+\lambda \u_0$. Then, by (2.33) and (9.12),
$$
\| \u_0 \|_{H^2(\O)} \le \wh{C}_\lambda \|\FF + \lambda \u_0\|_{L_2(\O)}
\le (\wh{C}_\lambda + \lambda \|g^{-1}\|_{L_\infty} \wt{\mathcal C}_1) \|\FF\|_{L_2(\O)}.
$$
Fixing $\lambda = \lambda_0 = {\mathcal C}_2 \|g^{-1}\|^{-1}_{L_\infty} +1$,
we obtain (9.13) with the constant
$\wh{C} = \wh{C}_{\lambda_0} + \lambda_0 \|g^{-1}\|_{L_\infty} \wt{\mathcal C}_1$.

\smallskip\noindent\textbf{9.3. Approximation of solutions in $L_2(\O;\C^n)$.}
Next theorem is an analog of Theorem 4.1 for the case where $\lambda=0$.

\smallskip\noindent\textbf{Theorem 9.2.} \textit{Suppose that} $\mathcal{O}\subset \R^d$
\textit{is a bounded domain with the boundary of class} $C^{1,1}$. \textit{Suppose that
the matrix} $g(\x)$
\textit{and DO} $b(\D)$ \textit{satisfy the assumptions of Subsection} 1.2.
\textit{Suppose that Condition} 2.1
 \textit{is satisfied. Let} $\mathbf{u}_{\varepsilon}$ \textit{be the solution of the problem}
 (9.6) \textit{and let} $\mathbf{u}_0$ \textit{be the solution of the problem} (9.11)
 \textit{with} $\FF \in {\mathcal H}(\O)$.
\textit{Then there exists a number} $\eps_1 \in (0,1]$
\textit{depending on the domain $\O$ and the lattice $\Gamma$ such that}
$$
\| \mathbf{u}_{\varepsilon} - \mathbf{u}_0
\|_{L_2 (\mathcal{O}; \mathbb{C}^{n})} \leq \wt{\mathcal C}_0 \varepsilon
\| \mathbf{F} \|_{L_2 (\mathcal{O}; \mathbb{C}^{n})},\ \ 0< \eps \leq \eps_1,
\eqno(9.14)
$$
\textit{or, in operator terms,}
$$
\| \mathcal{B}_{N,\varepsilon}^{-1} - (\mathcal{B}^0_N)^{-1}
 \|_{\mathcal{H}(\O) \to \mathcal{H}(\O)}
 \leq \wt{\mathcal C}_0 {\varepsilon},\ \ 0< \eps \leq \eps_1.
$$
\textit{The constant $\wt{\mathcal C}_0$ depends on}
$m$, $d$, $\alpha_0$, $\alpha_1$, $\| g \|_{L_\infty}$,
$\| g^{-1} \|_{L_\infty}$, \textit{the parameters of the lattice $\Gamma$,
the constants ${\mathcal C}_1$ and ${\mathcal C}_2$ from the inequality} (2.2),
\textit{the constant $\wt{\mathcal C}_1$ from} (9.4),
\textit{and the domain} $\mathcal{O}$.

\smallskip
\noindent\textbf{Proof.} We apply Theorem 4.1.
Let $\lambda$ be subject to the restriction (2.6).
Since $\B_{N,\eps}$ is the part of $\A_{N,\eps}$ in ${\mathcal H}(\O)$
and $\B_N^0$ is the part of $\A_N^0$ in the same subspace, it follows from (4.2) that
$$
\| (\mathcal{B}_{N,\varepsilon} +\lambda I)^{-1} - (\mathcal{B}^0_N + \lambda I)^{-1}
 \|_{\mathcal{H}(\O) \to \mathcal{H}(\O)}
 \leq {\mathcal C}_0(\lambda) {\varepsilon},\ \ 0< \eps \leq \eps_1.
 \eqno(9.15)
$$
Obviously, the solution $\u_\eps$ of the problem (9.6) also satisfies the equation
$(\B_{N,\eps} + \lambda I)\u_\eps =\FF_\eps$, where $\FF_\eps = \FF + \lambda \u_\eps \in {\mathcal H}(\O)$.
Consequently, \hbox{$\u_\eps = (\B_{N,\eps} +\lambda I)^{-1} \FF_\eps$}.
By (9.8), we have
$$
\|\FF_\eps \|_{L_2(\O)} \le \|\FF\|_{L_2(\O)} + \lambda \|\u_\eps\|_{L_2(\O)}
\le  (1+ \lambda \|g^{-1}\|_{L_\infty} \wt{\mathcal C}_1)\|\FF\|_{L_2(\O)}.
\eqno(9.16)
$$
We put $\u_{0,\eps}:= (\B^0_N + \lambda I)^{-1}\FF_\eps$.
By (9.15),
$$
\| \u_\eps - \u_{0,\eps} \|_{L_2(\O)} \le {\mathcal C}_0(\lambda) {\varepsilon} \|\FF_\eps\|_{L_2(\O)},
\quad 0< \eps \le \eps_1.
\eqno(9.17)
$$

From the other side, the solution $\u_0$ of the problem (9.11)
also satisfies the equation
$(\B_N^0 +\lambda I)\u_0 = \FF_0$, where $\FF_0 = \FF + \lambda \u_0 \in {\mathcal H}(\O)$.
Therefore, $\u_0 = (\B_N^0 +\lambda I)^{-1}\FF_0$.
We have:
$$
\u_{0,\eps} -\u_0 = (\B_N^0 +\lambda I)^{-1}(\FF_\eps - \FF_0)=
\lambda (\B_N^0 +\lambda I)^{-1} (\u_\eps - \u_0).
\eqno(9.18)
$$
Consequently,
$$
\u_\eps - \u_0 = \u_\eps - \u_{0,\eps} +
\lambda (\B_N^0 +\lambda I)^{-1} (\u_\eps - \u_0),
$$
whence
$$
\left(I - \lambda (\B_N^0 +\lambda I)^{-1}\right)(\u_\eps - \u_0) = \u_\eps - \u_{0,\eps}.
\eqno(9.19)
$$
By (9.10),
$$
\|\lambda (\B_N^0 +\lambda I)^{-1}\|_{{\mathcal H}(\O) \to {\mathcal H}(\O)}
\le \lambda \left(\lambda + \|g^{-1}\|^{-1}_{L_\infty} \wt{\mathcal C}_1^{-1}\right)^{-1} < 1,
$$
and so the operator $I - \lambda (\B_N^0 +\lambda I)^{-1}$ is invertible, and the norm of the
inverse operator satisfies the estimate
$$
\| \left(I - \lambda (\B_N^0 +\lambda I)^{-1}\right)^{-1}\|_{{\mathcal H}(\O) \to {\mathcal H}(\O)}
\le 1 + \lambda \| g^{-1}\|_{L_\infty} \wt{\mathcal C}_1.
\eqno(9.20)
$$
As a result, (9.19) implies that
$$
\u_\eps - \u_0
= \left(I - \lambda (\B_N^0 +\lambda I)^{-1}\right)^{-1}
(\u_\eps - \u_{0,\eps}).
\eqno(9.21)
$$
Combining (9.16), (9.17), (9.20), and (9.21), we obtain
$$
\|\u_\eps - \u_0\|_{L_2(\O)} \le (1 + \lambda \| g^{-1}\|_{L_\infty} \wt{\mathcal C}_1)^2
{\mathcal C}_0(\lambda) \eps \|\FF\|_{L_2(\O)},\quad 0< \eps \le \eps_1.
\eqno(9.22)
$$
Estimate (9.22) is proved for any $\lambda$ satisfying (2.6).
We can fix $\lambda$, for instance, putting
$\lambda = \lambda_0:= 1+ {\mathcal C}_2 \|g^{-1}\|^{-1}_{L_\infty}$.
Then we obtain (9.14) with the constant
$\wt{\mathcal C}_0 = (1 + \lambda_0 \| g^{-1}\|_{L_\infty} \wt{\mathcal C}_1)^2
{\mathcal C}_0(\lambda_0)$. $\ \bullet$

\smallskip\noindent\textbf{9.4. Approximation of solutions in $H^1(\O;\C^n)$.}
The corrector for the problem (9.6) is defined similarly to (4.5):
$$
K_N(\eps) = R_\O [\Lambda^\eps] S_\eps b(\D) P_\O (\B_N^0)^{-1}.
\eqno(9.23)
$$
The operator (9.23) is a continuous mapping of ${\mathcal H}(\O)$ into $H^1(\O;\C^n)$.
The first order approximation to the solution $\u_\eps$ of the problem (9.6) is given by
 $$
 \v_\eps = (\B_N^0)^{-1} \FF + \eps K_N(\eps) \FF=
 \u_0 + \eps \Lambda^\eps S_\eps b(\D) \wt{\u}_0,
 \eqno(9.24)
 $$
where $\wt{\u}_0 = P_\O \u_0$. Note that $\v_\eps$ belongs to $H^1(\O;\C^n)$ but
(in general) does not belong to $H^1_\perp(\O;\C^n)$.

The following result is an analog of
Theorem 4.2 in the case where $\lambda =0$.

\smallskip\noindent\textbf{Theorem 9.3.} \textit{Suppose that the assumptions of Theorem} 9.2
\textit{are satisfied. Let} $\mathbf{v}_{\varepsilon}$
\textit{be defined by} (9.24).
\textit{Then there exists a number} $\eps_1 \in (0,1]$
\textit{depending on the domain $\O$ and the lattice $\Gamma$ such that}
$$
\| \mathbf{u}_{\varepsilon} - \mathbf{v}_{\varepsilon}
\|_{H^1 (\mathcal{O}; \mathbb{C}^{n})} \leq \wt{\mathcal C} \varepsilon^{1/2}
\| \mathbf{F} \|_{L_2 (\mathcal{O}; \mathbb{C}^{n})},\quad 0< \eps \leq \eps_1,
\eqno(9.25)
$$
\textit{or, in operator terms,}
$$
\| \mathcal{B}_{N,\varepsilon}^{-1} - (\mathcal{B}^0_N)^{-1} -
\varepsilon K_{N} (\varepsilon) \|_{\mathcal{H} (\mathcal{O}) \to
H^1 (\mathcal{O}; \mathbb{C}^{n})} \leq \wt{\mathcal C} {\varepsilon}^{1/2}.
$$
\textit{For the flux} $\p_\eps:= g^\eps b(\D)\u_\eps$ \textit{we have}
$$
\| \p_\eps - \wt{g}^\eps S_\eps b(\D) \wt{\u}_0 \|_{L_2(\O;\C^m)}
\leq \wt{\mathcal C}' \eps^{1/2} \| \mathbf{F} \|_{L_2 (\mathcal{O}; \mathbb{C}^{n})},\ \ 0< \eps \leq \eps_1,
\eqno(9.26)
$$
\textit{where} $\wt{g}(\x)$ \textit{is the matrix} (1.9).
\textit{The constants $\wt{\mathcal C}$ and $\wt{\mathcal C}'$ depend only on}
$m$, $d$, $\alpha_0$, $\alpha_1$, $\| g \|_{L_\infty}$,
$\| g^{-1} \|_{L_\infty}$, \textit{the parameters of the lattice} $\Gamma$,
\textit{the constants ${\mathcal C}_1$ and ${\mathcal C}_2$ from} (2.2),
\textit{the constant $\wt{\mathcal C}_1$ from} (9.4),
\textit{and the domain} $\mathcal{O}$.

\smallskip
 \smallskip\noindent\textbf{Proof.} We apply Theorem 4.2.
Suppose that $\lambda$ satisfies condition (2.6).
Denote by $\mathcal P$ the orthogonal projection of
$L_2(\O;\C^n)$ onto the subspace ${\mathcal H}(\O)$.
Multiplying operators under the norm sign in (4.10) by the projection
$\mathcal P$ from the right, we arrive at
$$
\begin{aligned}
\|&(\B_{N,\eps} + \lambda I)^{-1} - (\B_{N}^0 + \lambda I)^{-1}
- \eps \Lambda^\eps S_\eps b(\D) P_\O (\B_{N}^0 + \lambda I)^{-1} \|_{{\mathcal H}(\O)\to H^1(\O)}
\cr
&\le {\mathcal C}(\lambda) \eps^{1/2},\quad 0< \eps \le \eps_1.
\end{aligned}
\eqno(9.27)
$$

As in the proof of Theorem 9.2, we use that
$\u_\eps = (\B_{N,\eps}+ \lambda I)^{-1}\FF_\eps$ and consider the function
$\u_{0,\eps} = (\B_{N}^0 + \lambda I)^{-1} \FF_\eps$.
By (9.27) and (9.16),
$$
\begin{aligned}
&\| \u_\eps - \u_{0,\eps} - \eps \Lambda^\eps S_\eps b(\D) \wt{\u}_{0,\eps} \|_{H^1(\O)}
\le {\mathcal C}(\lambda) \eps^{1/2} \|\FF_\eps\|_{L_2(\O)}
\cr
&\le {\mathcal C}(\lambda) (1+\lambda \|g^{-1}\|_{L_\infty} \wt{\mathcal C}_1) \eps^{1/2} \|\FF \|_{L_2(\O)},
\quad 0< \eps \le \eps_1,
\end{aligned}
\eqno(9.28)
$$
where $\wt{\u}_{0,\eps} = P_\O \u_{0,\eps}$.

Note that from (2.32) and (2.34), by multiplying operators under the norm sign
by the projection $\mathcal P$ from the right, it follows that
$$
\|(\B_{N}^0 + \lambda I)^{-1}\|_{{\mathcal H}(\O) \to H^1(\O)} \le c_\lambda^{-1},
\eqno(9.29)
$$
$$
\|(\B_{N}^0 + \lambda I)^{-1}\|_{{\mathcal H}(\O) \to H^2(\O)} \le \wh{C}_\lambda.
\eqno(9.30)
$$

Combining (9.18), (9.14), and (9.29), we obtain
$$
\begin{aligned}
&\| \u_{0,\eps} - \u_0 \|_{H^1(\O)} = \lambda \|(\B_{N}^0 + \lambda I)^{-1}(\u_\eps - \u_0)\|_{H^1(\O)}
\cr
&\le \lambda c_\lambda^{-1} \|\u_\eps - \u_0\|_{L_2(\O)}
\le \lambda c_\lambda^{-1} \wt{\mathcal C}_0 \eps \|\FF\|_{L_2(\O)},
\quad 0< \eps \le \eps_1.
\end{aligned}
\eqno(9.31)
$$

Now we consider the function $\eps \Lambda^\eps S_\eps b(\D) (\wt{\u}_{0,\eps} - \wt{\u}_0)$.
It is required to estimate its norm in $H^1(\O;\C^n)$.
Since this function is defined on the whole $\R^d$,
we will estimate its norm in $H^1(\R^d;\C^n)$.
By Proposition 1.5 and (1.4),
$$
\|\eps \Lambda^\eps S_\eps b(\D) (\wt{\u}_{0,\eps} - \wt{\u}_0)\|_{L_2(\R^d)}
\le \eps |\Omega|^{-1/2} \|\Lambda\|_{L_2(\Omega)} \alpha_1^{1/2}
\|\wt{\u}_{0,\eps} - \wt{\u}_0\|_{H^1(\R^d)}.
\eqno(9.32)
$$

Consider the derivatives
$$
\begin{aligned}
\partial_j \left(\eps \Lambda^\eps S_\eps b(\D) (\wt{\u}_{0,\eps} - \wt{\u}_0)\right)
&= (\partial_j \Lambda)^\eps S_\eps  b(\D) (\wt{\u}_{0,\eps} - \wt{\u}_0)
\cr
&+ \eps \Lambda^\eps S_\eps b(\D) \partial_j (\wt{\u}_{0,\eps} - \wt{\u}_0).
\end{aligned}
$$
By Proposition 1.5 and (1.4), we have:
$$
\begin{aligned}
&\sum_{j=1}^d
\|\partial_j (\eps \Lambda^\eps S_\eps b(\D) (\wt{\u}_{0,\eps} - \wt{\u}_0))\|^2_{L_2(\R^d)}
\cr
&\le 2 |\Omega|^{-1} \|\D \Lambda\|_{L_2(\Omega)}^2 \alpha_1 \|\wt{\u}_{0,\eps} - \wt{\u}_0\|_{H^1(\R^d)}^2
\cr
&+ 2\eps^2 |\Omega|^{-1} \| \Lambda\|_{L_2(\Omega)}^2 \alpha_1 \|\wt{\u}_{0,\eps} - \wt{\u}_0\|_{H^2(\R^d)}^2.
\end{aligned}
$$
Together with (9.32) this implies that
$$
\begin{aligned}
&\|\eps \Lambda^\eps S_\eps b(\D) (\wt{\u}_{0,\eps} - \wt{\u}_0) \|^2_{H^1(\O)}
\cr
&\le |\Omega|^{-1} \alpha_1 (3\eps^2 \|\Lambda\|^2_{L_2(\Omega)}+ 2\|\D\Lambda\|^2_{L_2(\Omega)})
\|\wt{\u}_{0,\eps} - \wt{\u}_0\|_{H^2(\R^d)}^2.
\end{aligned}
\eqno(9.33)
$$

From (9.14), (9.18), and (9.30) it follows that
$$
\begin{aligned}
&\|\wt{\u}_{0,\eps} - \wt{\u}_0\|_{H^2(\R^d)} \le
{\mathcal C}_\O \|{\u}_{0,\eps} - {\u}_0\|_{H^2(\O)}
\cr
&\le {\mathcal C}_\O \lambda \wh{C}_\lambda
\|\u_\eps - \u_0\|_{L_2(\O)}
\le \eps {\mathcal C}_\O \lambda \wh{C}_\lambda \wt{\mathcal C}_0 \|\FF\|_{L_2(\O)},
\quad 0< \eps \le \eps_1.
\end{aligned}
\eqno(9.34)
$$
Combining (9.33), (9.34), (1.10), and (1.11), we obtain
$$
\|\eps \Lambda^\eps S_\eps b(\D) (\wt{\u}_{0,\eps} - \wt{\u}_0)\|_{H^1(\O)}
\le \check{C}_\lambda \eps \|\FF\|_{L_2(\O)},\quad 0< \eps \le \eps_1,
\eqno(9.35)
$$
 where $\check{C}_\lambda= m^{1/2}\alpha_1^{1/2} \alpha_0^{-1/2}
 \|g\|^{1/2}_{L_\infty} \|g^{-1}\|^{1/2}_{L_\infty} \left( 2+ 3(2r_0)^{-2} \right)^{1/2}
 {\mathcal C}_\O \lambda \wh{C}_\lambda \wt{\mathcal C}_0$.

As a result, relations (9.28), (9.31), and (9.35) imply that
$$
\| \u_\eps - \v_\eps\|_{H^1(\O)} \le \eps^{1/2}
\left( {\mathcal C}(\lambda)(1+ \lambda \|g^{-1}\|_{L_\infty}
\wt{\mathcal C}_1) + \lambda c_\lambda^{-1} \wt{\mathcal C}_0+ \check{C}_\lambda\right) \|\FF\|_{L_2(\O)}
\eqno(9.36)
$$
for $0< \eps \le \eps_1$.
Estimate (9.36) is proved for any $\lambda$ satisfying (2.6).
Putting $\lambda = \lambda_0:= 1+ {\mathcal C}_2 \|g^{-1}\|^{-1}_{L_\infty}$, we obtain
estimate (9.25) with the constant
 $\wt{\mathcal C}= {\mathcal C}(\lambda_0)(1+ \lambda_0 \|g^{-1}\|_{L_\infty}
\wt{\mathcal C}_1) + \lambda_0 c_{\lambda_0}^{-1} \wt{\mathcal C}_0+ \check{C}_{\lambda_0}$.

It remains to prove (9.26). The arguments are similar to the proof of (4.11) (see Subsection 5.4).
By (9.25), (1.2), and (1.5), we have
$$
\| \p_\eps - g^\eps b(\D) \v_\eps \|_{L_2(\O)}
\le \|g\|_{L_\infty} (d \alpha_1)^{1/2} \wt{\mathcal C} \eps^{1/2} \|\FF\|_{L_2(\O)},
\quad 0< \eps \le \eps_1.
\eqno(9.37)
$$
Similarly to (5.40),
$$
\| g^\eps b(\D) \v_\eps - \wt{g}^\eps S_\eps b(\D)\wt{\u}_0 \|_{L_2(\O)}
\le C_6 \eps \|\D b(\D) \wt{\u}_0\|_{L_2(\R^d)}.
\eqno(9.38)
$$
By (1.4) and (9.13),
$$
\begin{aligned}
&\|\D b(\D) \wt{\u}_0\|_{L_2(\R^d)} \le \alpha_1^{1/2} \| \wt{\u}_0 \|_{H^2(\R^d)}
\cr
&\le \alpha_1^{1/2} C_\O \| {\u}_0 \|_{H^2(\O)} \le
\alpha_1^{1/2} C_\O \wh{C} \| \FF \|_{L_2(\O)}.
\end{aligned}
\eqno(9.39)
$$
As a result, (9.37)--(9.39) imply (9.26)
with the constant $\wt{\mathcal C}' = \|g\|_{L_\infty} (d\alpha_1)^{1/2} \wt{\mathcal C} +
C_6 \alpha_1^{1/2} C_\O \wh{C}$. $\ \bullet$

\smallskip\noindent\textbf{9.5. Results in the case where $\Lambda \in L_\infty$.}
Under Condition 1.10, instead of the corrector (9.23), it is possible to use
a simpler corrector
$$
K^0_N(\eps) = [\Lambda^\eps] b(\D) (\B_N^0)^{-1},
$$
which is a continuous mapping of ${\mathcal H}(\O)$ into $H^1(\O;\C^n)$.
Instead of (9.24), we consider another approximation
of the solution $\u_\eps$ of the problem (9.6):
$$
\check{\v}_\eps = (\B_N^0)^{-1} \FF + \eps K_N^0(\eps)\FF =
\u_0 + \eps \Lambda^\eps b(\D) \u_0.
\eqno(9.40)
$$

Next result is an analog of Theorem 7.1 in the case where $\lambda =0$.

\smallskip\noindent\textbf{Theorem 9.4.} \textit{Suppose that the assumptions of Theorem} 9.2
\textit{and Condition} 1.10 \textit{are satisfied. Let $\check{\v}_\eps$ be defined by} (9.40).
\textit{Then there exists a number $\eps_1 \in (0,1]$
depending on the domain $\O$ and the lattice $\Gamma$ such that}
$$
\| \u_\eps - \check{\v}_\eps\|_{H^1(\O;\C^n)} \le \check{\mathcal C} \eps^{1/2}
\| \FF \|_{L_2(\O;\C^n)},\quad 0 < \eps \le \eps_1,
\eqno(9.41)
$$
\textit{or, in operator terms,}
$$
\| \mathcal{B}_{N,\varepsilon}^{-1} - (\mathcal{B}^0_N)^{-1} -
\varepsilon K^0_{N} (\varepsilon) \|_{{\mathcal H}(\mathcal{O}) \to
H^1 (\mathcal{O}; \mathbb{C}^{n})} \leq \check{\mathcal C} {\varepsilon}^{1/2}.
$$
\textit{For the flux} $\p_\eps= g^\eps b(\D)\u_\eps$ \textit{we have}
$$
\| \p_\eps - \wt{g}^\eps b(\D) {\u}_0 \|_{L_2(\O;\C^m)}
\leq \check{\mathcal C}' \eps^{1/2} \| \mathbf{F} \|_{L_2 (\mathcal{O}; \mathbb{C}^{n})},\ \ 0< \eps \leq \eps_1,
\eqno(9.42)
$$
\textit{where} $\wt{g}(\x)$ \textit{is the matrix} (1.9).
\textit{The constants $\check{\mathcal C}$ and $\check{\mathcal C}'$ depend on}
$m$, $d$, $\alpha_0$, $\alpha_1$, $\| g \|_{L_\infty}$,
$\| g^{-1} \|_{L_\infty}$, \textit{the parameters of the lattice} $\Gamma$,
\textit{the constants ${\mathcal C}_1$ and ${\mathcal C}_2$ from the inequality} (2.2),
\textit{the constant $\wt{\mathcal C}_1$ from} (9.4),
\textit{the norm $\|\Lambda\|_{L_\infty}$, and the domain} $\mathcal{O}$.

\smallskip
\noindent\textbf{Proof.} The proof is similar to that of Theorem 7.1.
It is required to estimate the $H^1$-norm of the function
$\check{\v}_\eps - \v_\eps = \eps \Lambda^\eps (I-S_\eps) b(\D) \wt{\u}_0$.

From (1.4) and (9.13) it follows that
$$
\| b(\D) \wt{\u}_0 \|_{L_2(\R^d)} \le \alpha_1^{1/2} \|\wt{\u}_0\|_{H^1(\R^d)}
\le \alpha_1^{1/2} C_\O \| {\u}_0\|_{H^2(\O)} \le
\alpha_1^{1/2} C_\O \wh{C} \| \FF \|_{L_2(\O)}.
\eqno(9.43)
$$
Estimates (7.6) and (7.7) can be applied.
Together with (9.39) and (9.43) they yield
$$
\| \check{\v}_\eps - \v_\eps \|_{H^1(\O)}
= \| \eps \Lambda^\eps (I-S_\eps) b(\D) \wt{\u}_0 \|_{H^1(\O)}
\le {\mathcal C}_{12}^0 \eps \|\FF\|_{L_2(\O)},
\eqno(9.44)
$$
where ${\mathcal C}_{12}^0 = \alpha_1^{1/2}
C_\O \wh{C} \left( 2(3+2\beta_2)^{1/2} \|\Lambda\|_{L_\infty} +(2\beta_1)^{1/2}r_1\right)$.

Now (9.25) and (9.44) imply (9.41) with the constant
$\check{\mathcal C} = \wt{\mathcal C} + {\mathcal C}^0_{12}$.

It remains to check (9.42). From (9.41), (1.2), and (1.5) it follows that
$$
\| \p_\eps - g^\eps b(\D) \check{\v}_\eps \|_{L_2(\O)}
\le \eps^{1/2} d^{1/2} \alpha_1^{1/2} \|g\|_{L_\infty} \check{\mathcal C} \| \FF \|_{L_2(\O)},
\quad 0< \eps \le \eps_1.
\eqno(9.45)
$$
Similarly to (7.10) and (7.11), taking (9.13) into account, we have
$$
\begin{aligned}
&\| g^\eps b(\D) \check{\v}_\eps - \wt{g}^\eps b(\D) \u_0 \|_{L_2(\O)}
\le
\eps d \alpha_1 \|g\|_{L_\infty} \|\Lambda\|_{L_\infty}
\|\u_0\|_{H^2(\O)}
\cr
&\le \eps d \alpha_1 \|g\|_{L_\infty} \|\Lambda\|_{L_\infty} \wh{C} \| \FF \|_{L_2(\O)}.
\end{aligned}
\eqno(9.46)
$$
Relations (9.45) and (9.46) imply (9.42) with the constant
$\check{\mathcal C}' = d^{1/2} \alpha_1^{1/2} \|g\|_{L_\infty} \check{\mathcal C} +
d \alpha_1 \|g\|_{L_\infty} \|\Lambda\|_{L_\infty} \wh{C}$. $\ \bullet$

\smallskip
\noindent\textbf{9.6. Special cases.} Now we distinguish special cases.
Next statement follows from Theorem 9.3 and Proposition 1.2.

\smallskip
\noindent\textbf{Proposition 9.5.} \textit{Let $\u_\eps$ be the solution of the problem} (9.6),
\textit{and let $\u_0$ be the solution of the problem} (9.11).
\textit{If $g^0 = \overline{g}$, i.~e., relations} (1.13) \textit{are satisfied, then
$\Lambda=0$, $K_N(\eps)=0$, and we have}
$$
\|\u_\eps - \u_0 \|_{H^1(\O;\C^n)} \le \wt{\mathcal C} \eps^{1/2} \|\FF\|_{L_2(\O;\C^n)},\quad 0< \eps \le \eps_1.
$$

\smallskip
Applying Propositions 1.3, 1.13
and the statement of Theorem 9.4 concerning the flux, we obtain the following statement.

\smallskip
\noindent\textbf{Proposition 9.6.} \textit{Let $\u_\eps$ be the solution of the problem} (9.6),
\textit{and let $\u_0$ be the solution of the problem} (9.11). \textit{Let}
$\p_\eps = g^\eps b(\D)\u_\eps$.
\textit{If $g^0 = \underline{g}$, i.~e., relations} (1.14)
\textit{are satisfied, then}
$$
\|\p_\eps - g^0 b(\D)\u_0 \|_{L_2(\O;\C^m)} \le
\check{\mathcal C}' \eps^{1/2} \|\FF\|_{L_2(\O;\C^n)},\quad 0< \eps \le \eps_1.
$$

\smallskip\noindent\textbf{9.7. Approximation of solutions in a strictly interior subdomain.}
As in \S 8, it is possible to obtain a sharp order error estimate for approximation
of the solution in $H^1(\O';\C^n)$, where $\O'$ is a strictly interior subdomain of the domain $\O$.

\smallskip\noindent\textbf{Theorem 9.7.} \textit{Suppose that the assumptions of Theorem} 9.3
\textit{are satisfied. Let $\O'$ be a strictly interior subdomain of the domain $\O$. Let}
$\delta:=\text{dist}\, \{\O';\partial \O\}$.
\textit{Then there exists a number $\eps_1 \in (0,1]$ depending on the domain
$\O$ and the lattice $\Gamma$ such that}
$$
\| \u_\eps - {\v}_\eps\|_{H^1(\O';\C^n)} \le \wt{\mathfrak C}_{\delta} \eps
\| \FF \|_{L_2(\O;\C^n)},\quad 0 < \eps \le \eps_1,
\eqno(9.47)
$$
\textit{or, in operator terms},
$$
\| \B_{N,\eps}^{-1} - (\B_N^0)^{-1} - \eps K_{N}(\eps)
\|_{\mathcal{H}(\O) \to H^1(\O';\C^n)} \le
\wt{\mathfrak C}_{\delta} \eps,\quad 0 < \eps \le \eps_1.
$$
\textit{For the flux} $\p_\eps= g^\eps b(\D)\u_\eps$ \textit{we have}
$$
\| \p_\eps - \wt{g}^\eps S_\eps b(\D) \wt{\u}_0 \|_{L_2(\O';\C^m)}
\leq \wt{\mathfrak C}'_{\delta} \eps
\| \mathbf{F} \|_{L_2 (\mathcal{O}; \mathbb{C}^{n})},\ \ 0< \eps \leq \eps_1,
\eqno(9.48)
$$
\textit{where} $\wt{g}(\x)$ \textit{is the matrix} (1.9).
\textit{The constants in estimates are given by} $\wt{\mathfrak C}_{\delta}=
\wt{\mathfrak C}_1 \delta^{-1} + \wt{\mathfrak C}_2$,
$\wt{\mathfrak C}_{\delta}' = \wt{\mathfrak C}'_1 \delta^{-1} + \wt{\mathfrak C}'_2$,
 \textit{where $\wt{\mathfrak C}_1$, $\wt{\mathfrak C}_2$,
 $\wt{\mathfrak C}'_1$, $\wt{\mathfrak C}'_2$ depend on}
$m$, $d$, $\alpha_0$, $\alpha_1$, $\| g \|_{L_\infty}$,
$\| g^{-1} \|_{L_\infty}$, \textit{the parameters of the lattice} $\Gamma$,
\textit{the constants ${\mathcal C}_1$ and ${\mathcal C}_2$ from the inequality} (2.2),
\textit{the constant $\wt{\mathcal C}_1$ from} (9.4), \textit{and the domain} $\mathcal{O}$.

\smallskip\noindent\textbf{Proof.} We apply Theorem 8.1. Assume that $\lambda$ satisfies
(2.6). Multiplying operators under the norm sign in (8.2) by the projection
$\mathcal P$ from the right, we arrive at
$$
\begin{aligned}
&\| (\B_{N,\eps}+\lambda I)^{-1} - (\B_{N}^0 +\lambda I)^{-1} - \eps \Lambda^\eps S_\eps b(\D)P_\O
(\B_{N}^0 +\lambda I)^{-1}\|_{{\mathcal H}(\O) \to H^1(\O')}
\cr
&\le {\mathfrak C}_\delta(\lambda) \eps,\quad 0< \eps \le \eps_1.
\end{aligned}
\eqno(9.49)
$$
Next, we proceed like in the proof of Theorem 9.3.
We write $\u_\eps$ as
$\u_\eps = (\B_{N,\eps}+\lambda I)^{-1} \FF_\eps$ and consider
$\u_{0,\eps} = (\B_{N}^0 +\lambda I)^{-1} \FF_\eps$. By (9.16) and (9.49),
$$
\begin{aligned}
&\| \u_\eps - \u_{0,\eps} - \eps \Lambda^\eps S_\eps b(\D) \wt{\u}_{0,\eps} \|_{H^1(\O')}
\le
{\mathfrak C}_\delta(\lambda) \eps \|\FF_\eps\|_{L_2(\O)}
\cr
&\le
{\mathfrak C}_\delta(\lambda) (1+ \lambda \|g^{-1}\|_{L_\infty} \wt{\mathcal C}_1)
\eps \|\FF \|_{L_2(\O)},\quad 0< \eps \le \eps_1.
\end{aligned}
\eqno(9.50)
$$
Now relations (9.31), (9.35), and (9.50) imply that
$$
\| \u_\eps - \v_\eps\|_{H^1(\O')}
\le \eps \left( {\mathfrak C}_\delta(\lambda)(1+ \lambda \|g^{-1}\|_{L_\infty}
\wt{\mathcal C}_1) + \lambda c_\lambda^{-1} \wt{\mathcal C}_0+ \check{C}_\lambda\right) \|\FF\|_{L_2(\O)}
\eqno(9.51)
$$
for $0< \eps \le \eps_1$. The inequality (9.51) is valid for any $\lambda$ satisfying (2.6).
Putting $\lambda = \lambda_0 = 1+ {\mathcal C}_2 \|g^{-1}\|^{-1}_{L_\infty}$,
we obtain (9.47) with the constant $\wt{\mathfrak C}_\delta = \wt{\mathfrak C}_1 \delta^{-1} +
\wt{\mathfrak C}_2$, where
$$
\wt{\mathfrak C}_1=
{\mathfrak C}_1(\lambda_0)(1+ \lambda_0 \|g^{-1}\|_{L_\infty} \wt{\mathcal C}_1),
\eqno(9.52)
$$
$\wt{\mathfrak C}_2=
{\mathfrak C}_2(\lambda_0)(1+ \lambda_0 \|g^{-1}\|_{L_\infty} \wt{\mathcal C}_1)
+ \lambda_0 c_{\lambda_0}^{-1} \wt{\mathcal C}_0+ \check{C}_{\lambda_0}$.

It remains to prove (9.48). From (9.47), (1.2), and (1.5) it follows that
$$
\| \p_\eps - g^\eps b(\D)\v_\eps \|_{L_2(\O')}
\le \|g\|_{L_\infty} d^{1/2} \alpha_1^{1/2} \wt{\mathfrak C}_\delta \eps \| \FF \|_{L_2(\O)},
\quad 0< \eps \le \eps_1.
\eqno(9.53)
$$
Combining (9.53), (9.38), and (9.39), we arrive at
$$
\| \p_\eps - \wt{g}^\eps S_\eps b(\D) \wt{\u}_0 \|_{L_2(\O')}
\le \left(\|g\|_{L_\infty} d^{1/2} \alpha_1^{1/2} \wt{\mathfrak C}_\delta
+ C_6 \alpha_1^{1/2} C_\O \wh{C} \right) \eps \| \FF \|_{L_2(\O)}
$$
for $0< \eps \le \eps_1$. This yields (9.48) with the constant
$\wt{\mathfrak C}'_\delta = \wt{\mathfrak C}'_1 \delta^{-1} + \wt{\mathfrak C}'_2$,
where
$$
\wt{\mathfrak C}'_1 = \|g\|_{L_\infty} d^{1/2} \alpha_1^{1/2} \wt{\mathfrak C}_1,\ \
\wt{\mathfrak C}'_2 = \|g\|_{L_\infty} d^{1/2} \alpha_1^{1/2} \wt{\mathfrak C}_2 +
C_6 \alpha_1^{1/2} C_\O \wh{C}.
\ \ \bullet
\eqno(9.54)
$$

\smallskip
Under Condition 1.10, we obtain the following result.

\smallskip\noindent\textbf{Theorem 9.8.} \textit{Suppose that the assumptions of Theorem} 9.4
\textit{are satisfied. Let $\O'$ be a strictly interior subdomain of the domain $\O$, and let}
$\delta:=\text{dist}\, \{\O';\partial \O\}$.
\textit{Then there exists a number $\eps_1 \in (0,1]$ depending on the domain $\O$
and the lattice $\Gamma$ such that}
$$
\| \u_\eps - \check{\v}_\eps\|_{H^1(\O';\C^n)} \le \check{\mathfrak C}_{\delta} \eps
\| \FF \|_{L_2(\O;\C^n)},\quad 0 < \eps \le \eps_1,
\eqno(9.55)
$$
\textit{or, in operator terms,}
$$
\| \B_{N,\eps}^{-1} - (\B_N^0)^{-1} - \eps K_{N}^0(\eps)
\|_{\mathcal{H}(\O) \to H^1(\O';\C^n)} \le
\check{\mathfrak C}_{\delta} \eps,\quad 0 < \eps \le \eps_1.
$$
\textit{For the flux} $\p_\eps= g^\eps b(\D)\u_\eps$ \textit{we have}
$$
\| \p_\eps - \wt{g}^\eps b(\D) {\u}_0 \|_{L_2(\O';\C^m)}
\leq \check{\mathfrak C}'_\delta \eps
\| \mathbf{F} \|_{L_2 (\mathcal{O}; \mathbb{C}^{n})},\ \ 0< \eps \leq \eps_1,
\eqno(9.56)
$$
\textit{where} $\wt{g}(\x)$ \textit{is the matrix} (1.9).
\textit{The constants in estimates are given by} $\check{\mathfrak C}_\delta =
\wt{\mathfrak C}_1 \delta^{-1} + \check{\mathfrak C}_2$,
$\check{\mathfrak C}_{\delta}' = \wt{\mathfrak C}'_1 \delta^{-1} + \check{\mathfrak C}'_2$,
 \textit{where $\wt{\mathfrak C}_1$, $\check{\mathfrak C}_2$,
 $\wt{\mathfrak C}'_1$, $\check{\mathfrak C}'_2$ depend on}
$m$, $d$, $\alpha_0$, $\alpha_1$, $\| g \|_{L_\infty}$,
$\| g^{-1} \|_{L_\infty}$, \textit{the parameters of the lattice} $\Gamma$,
\textit{the constants ${\mathcal C}_1$ and ${\mathcal C}_2$ from the inequality} (2.2),
\textit{the constant $\wt{\mathcal C}_1$ from} (9.4),
\textit{the norm $\|\Lambda\|_{L_\infty}$, and the domain} $\mathcal{O}$.

\smallskip\noindent\textbf{Proof.} We apply the proof of Theorem 9.4.
Under Condition 1.10, (9.44) is true. Together with
(9.47) it implies (9.55) with the constant
 $\check{\mathfrak C}_\delta =
\wt{\mathfrak C}_1 \delta^{-1} + \check{\mathfrak C}_2$,
where $\wt{\mathfrak C}_1$ is defined by (9.52) and
$\check{\mathfrak C}_2 = \wt{\mathfrak C}_2 + {\mathcal C}_{12}^0$.

 Let us check (9.56). From (9.55), (1.2), and (1.5) it follows that
 $$
 \| \p_\eps - g^\eps b(\D) \check{\v}_\eps \|_{L_2(\O')} \le
 \|g\|_{L_\infty} d^{1/2} \alpha_1^{1/2} \check{\mathfrak C}_\delta \eps \|\FF\|_{L_2(\O)},
 \quad 0< \eps \le \eps_1.
 \eqno(9.57)
 $$
 Under Condition 1.10, (9.46) is true.
 Now relations (9.46) and (9.57) implies (9.56) with the constant
 $\check{\mathfrak C}_{\delta}' = \wt{\mathfrak C}'_1 \delta^{-1} + \check{\mathfrak C}'_2$,
 where $\wt{\mathfrak C}'_1$ is defined by (9.54) and
 $\check{\mathfrak C}'_2 = \|g\|_{L_\infty} d^{1/2} \alpha_1^{1/2} \check{\mathfrak C}_2
 + d \alpha_1 \|g\|_{L_\infty} \|\Lambda\|_{L_\infty} \wh{C}$. $\ \bullet$

\smallskip
What was said in Remark 8.3 concerns also Theorems 9.7 and 9.8.

\end{document}